\documentclass[11pt]{article}
\usepackage{bbm}
\usepackage{mathrsfs}
\usepackage{amsfonts}
\usepackage{amssymb}
\usepackage{amssymb,amsmath,amsthm, amsfonts}

\newtheorem{definition}{Definition}
\newtheorem{theorem}{Theorem}
\newtheorem{corollary}{Corollary}
\newtheorem{lemma}{Lemma}
\newtheorem{remark}{Remark}

\newtheorem{assumption}{Assumption}
\def\qed{ \ \vrule width.2cm height.2cm depth0cm\smallskip}

\newcommand{\la}{\langle}
\newcommand{\ra}{\rangle}
\newcommand{\hP}{\hat\dbP}

\newcommand{\ba}{\begin{array}}
\newcommand{\ea}{\end{array}}
\newcommand{\be}{\begin{equation}}
\newcommand{\ee}{\end{equation}}
\newcommand{\bea}{\begin{eqnarray}}
\newcommand{\eea}{\end{eqnarray}}
\newcommand{\beaa}{\begin{eqnarray*}}
\newcommand{\eeaa}{\end{eqnarray*}}

\def\neg{\negthinspace}

\def\a{\alpha}
\def\b{\beta}

\def\d{\delta}

\def\l{\lambda}
\def\m{\mu}
\def\n{\nu}
\def\si{\sigma}

\def\f{\varphi}
\def\th{\theta}
\def\o{\omega}
\def\h{\widehat}
\def\G{\Gamma}

\def\Th{\Theta}

\def\Si{\Sigma}

\def\O{\Omega}
\def\cF{{\cal F}}
\def\cG{{\cal G}}
\def\cH{{\cal H}}

\def\cP{{\cal P}}

\def\hC{\mathbb{C}}
\def\hD{\mathbb{D}}
\def\hE{\mathbb{E}}
\def\hF{\mathbb{F}}
\def\hG{\mathbb{G}}
\def\hH{\mathbb{H}}

\def\hN{\mathbb{N}}

\def\hP{\mathbb{P}}
\def\hQ{\mathbb{Q}}
\def\hR{\mathbb{R}}

\def\hX{\mathbb{X}}

\def\sB{\mathscr{B}}

\def\sE{\mathscr{E}}

\def\sH{\mathscr{H}}

\def\scL{\mathscr{L}}
\def\sM{\mathscr{M}}

\def\sP{\mathscr{P}}

\def\sT{\mathscr{T}}
\def\sU{\mathscr{U}}

\def\no{\noindent}

\def\ss{\smallskip}
\def\ms{\medskip}
\def\bs{\bigskip}
\def\q{\quad}
\def\qq{\qquad}

\def\pa{\partial}
\def\cd{\cdot}
\def\cds{\cdots}
\def\lan{{\langle}}
\def\ran{{\rangle}}

\def\be{\begin{equation}}
\def\bel{\begin{equation}\label}
\def\ee{\end{equation}}
\def\bt{\begin{theorem}}
\def\bcd{\begin{condition}}
\def\ecd{\end{condition}}
\def\et{\end{theorem}}
\def\bc{\begin{corollary}}
\def\ec{\end{corollary}}
\def\bde{\begin{definition}}
\def\ede{\end{definition}}
\def\bl{\begin{lemma}}
\def\el{\end{lemma}}
\def\bp{\begin{proposition}}
\def\ep{\end{proposition}}
\def\br{\begin{remark}}
\def\er{\end{remark}}
\def\ba{\begin{array}}
\def\ea{\end{array}}
\def\ed{\end{document}}

\def\square#1{\vbox{\hrule\hbox{\vrule height#1%
     \kern#1\vrule}\hrule}}
\def\rectangle#1#2{\vbox{\hrule\hbox{\vrule height#1%
     \kern#2\vrule}\hrule}}
\def\qed{\hfill \vrule height7pt width3pt depth0pt}

\def\sqr#1#2{{\vcenter{\vbox{\hrule height.#2pt
              \hbox{\vrule width.#2pt height#1pt \kern#1pt \vrule width.#2pt}
              \hrule height.#2pt}}}}

\def\qed{ \hfill \vrule width.25cm height.25cm depth0cm\smallskip}
\newcommand{\dfnn}{\stackrel{\triangle}{=}}
\newcommand{\basa}{\begin{assumption}}
\newcommand{\easa}{\end{assumption}}

\newcommand{\bas}{\begin{assum}}
\newcommand{\eas}{\end{assum}}

\def\lan{\mathop{\langle}}
\def\ran{\mathop{\rangle}}

\def\limw{\mathop{\buildrel w\over\rightarrow}}

\def\pa{\partial}
\def\h{\widehat}
\def\wt{\widetilde}
 \def\cd{\cdot}
\def\cds{\cdots}

\def\dis{\displaystyle}
\def\wt{\widetilde}

\def\1{{\bf 1}}

\def\:{\!:\!}

at 9pt
\begin{document}
\newtheorem{thm}{Theorem}[section]
\newtheorem{lem}[thm]{Lemma}
\newtheorem{cor}[thm]{Corollary}
\newtheorem{prop}[thm]{Proposition}
\newtheorem{rem}[thm]{Remark}
\newtheorem{eg}[thm]{Example}
\newtheorem{defn}[thm]{Definition}
\newtheorem{assum}[thm]{Assumption}

\renewcommand {\theequation}{\arabic{section}.\arabic{equation}}
\def\thesection{\arabic{section}}

\title{\bf A Mean-field Stochastic Control Problem with Partial Observations\footnote{The authors would like to dedicate this
paper to Prof. Hans-J\"urgen Engelbert, in the occasion of his 70th birthday, for his generous guidance and inspirational
discussions throughout the past decades.  \newline\ \ \indent $^1$\!\!  Laboratoire de Math\'{e}matiques,
Universit\'{e} de Bretagne-Occidentale, F-29285 Brest Cedex, France,
email: Rainer.Buckdahn@univ-brest.fr; School of Mathematics, Shandong University, Jinan, 250100, China. This work is part of the French ANR project CAESARS (ANR-15-CE05-0024).
\newline\ \ \indent $^2$\!\! Corresponding author. School of Mathematics, Shandong University, Weihai. Weihai, 264209, China. Email: juanli@sdu.edu.cn. This author has been supported by the NSF of P.R.China (No. 11222110), Shandong Province (No. JQ201202), NSFC-RS (No. 11661130148; NA150344).
\newline\ \ \indent $^3$\!\!  Department of
Mathematics, University of Southern California, Los Angeles, 90089, USA.
Email: jinma@usc.edu. This author is supported in part by US NSF grant  \#1106853. }}
\author{Rainer Buckdahn$^1$,~~{Juan Li}$^2$,~~and~~Jin Ma$^3$}

\date{\today}
\maketitle
\vspace{-8mm}
\begin{abstract}
In this paper we are interested in a new type of {\it mean-field}, non-Markovian stochastic control problems with partial
observations. More precisely, we assume that the coefficients of the controlled dynamics depend not only on the paths of the state,
but also on the conditional law of the state, given the observation to date. Our problem is strongly motivated by the recent
study of the mean field games and the related McKean-Vlasov stochastic control problem,
but with added aspects of path-dependence and partial observation.
We shall first investigate the well-posedness of the state-observation dynamics, with combined reference probability measure arguments in nonlinear filtering theory and  the Schauder fixed point theorem. We then study the stochastic control problem
with a partially  observable system in which the conditional law appears nonlinearly in both the coefficients of the system
and cost function. As a consequence the control problem is intrinsically ``time-inconsistent", and we prove that the Pontryagin Stochastic Maximum Principle holds in this case and characterize the adjoint equations,
which turn out to be a new form of mean-field type BSDEs.
\end{abstract}

\vfill


\bs

\no{\it 2000 AMS Mathematics subject classification:} 60H10,30;
93E03,11,20.

\eject

\section{Introduction}
\label{sect-Introduction}
\setcounter{equation}{0}
In this paper we are interested in the following {\it mean-field-type} stochastic control problem,
on a given filtered probability space $(\O, \cF, \hP;\hF=\{\cF_t\}_{t\ge 0})$:
 \bea
 \label{SDE-0}
  \left\{\ba{lll}
  dX_t= \hE\{b(t,  \f_{\cd\wedge t},\hE[X_{ t}|\cG_t], u)\}|_{\f=X, u=u_t}dt+
  \hE\{\si(t,  \f_{\cd\wedge t},   \hE[X_{t}|\cG_t], u)\}|_{\f=X, u=u_t}dB_t, \\
  X_0=x,
 \ea\right.
 % \q t\ge 0,
\eea
where $B$ is an $\hF$-Brownian motion, $b$ and $\si$ are measurable functions satisfying reasonable conditions,
$\f_{\cd \wedge t}$ and $X_{\cd\wedge t}$ denote the  continuous function and process, respectively, ``stopped" at $t$;
 $\hG\dfnn\{\cG_t\}_{t\ge0}$ is a given filtration that could  involve the information  of $X$ itself, and
 $u=\{u_t:t\ge 0\}$  is the ``control process", assumed to be adapted to a filtration $\hH=\{\cH_t\}_{t\ge0}$, where
 $\cH_t\subseteq \cF^X_t\vee\cG_t$,  $t\ge 0$.   We note that if $\cG_t=\{\emptyset, \O\}$, for all $t\ge 0$ (i.e., the
 conditional expectation in (\ref{SDE-0}) becomes expectation), $\cH_t=\cF^X_t$, and coefficients are ``Markovian"
 (i.e., $\f_{\cd\wedge t}=\f_t$), then the problem becomes a  stochastic control problem with McKean-Vlasov dynamics
 and/or a Mean-field game (see, for example, \cite{CD1,CD2,CD3} in its ``forward" form, and \cite{BDL, BDLP, BLP} in
 its ``backward" form). On the other hand, when $\hG$ is a given filtration, this is the so-called {\it conditional
 mean-field SDE} (CMFSDE for short) studied in \cite{CZ}. We note that in that case the conditioning is essentially
 ``open-looped".

The  problem that this paper is particularly focusing on is when $\cG_t=\cF^Y_t$, $t\ge 0$, where $Y$ is an
``observation process" of the dynamics of $X$, i.e., the case when the pair $(X, Y)$ forms a ``close-looped"
or ``coupled" CMFSDE.
More precisely, we shall consider the following partially observed controlled dynamics (assuming $b=0$ for
notational simplicity):
\bea
\label{dynamics}
\left\{\ba{lll}
\dis  dX_t=
%\hE\{b(t,  \f_{\cd\wedge t},\hE[X_{ t}|\cF^Y_t], u)\}|_{\f=X, u=u_t}dt+
  \hE\{\si(t,  \f_{\cd\wedge t},   \hE[X_{t}|\cF^Y_t], u)\}|_{\f=X, u=u_t}dB^1_t;\ms \\
\dis dY_t=h(t, X_{ t})dt +\hat \si d B^2_t; \qq X_0=x, ~Y_0=0.
\ea\right.
\eea
Here $X$ is the  ``signal" process that can only be observed through $Y$,  $(B^1, B^2)$ is a standard Brownian motion,
and $\hat \si$ is a constant. We should
%specify the probability  $\hP^0$ to facilitate the discussion below. W
note that in SDEs (\ref{dynamics}) the conditioning filtration $\hF^Y$ now depends on $X$ itself, therefore it is much
more convoluted than the CMFSDE we have seen in the literature. Furthermore, the path-dependent nature of the
coefficients makes the SDE essentially {\it non-Markovian}.
%because of its intrinsic ``path-dependent" nature. A more precise formulation will be given in section 2 below.   But
Such form of CMFSDEs, to the best of our knowledge,
% to date the problem  of the form (\ref{SDE-0})
has not been explored fully in the literature.

% (through the observation process $Y$).

Our study of the CMFSDE (\ref{dynamics}) is strongly motivated by the following variation of the mean-field game in a finance
context, which would result in a type of  stochastic control problem involving a controlled dynamics of such a form. Consider a firm whose {\it fundamental value}, under the risk neutral measure $\hP^0$ with zero interest, evolves
 as the following SDE with ``stochastic volatility" $\si=\si(t,\o)$, $(t,\o)\in[0,\infty)\times\O$:
 \bea
 \label{fvalue}
X_t=x+\int_0^t \si(s, \cd) dB^1_s, \q t\ge 0,
\eea
where $B^1$ is the intrinsic noise from inside the firm. We assume that such fundamental value process cannot be observed directly, but can be observed through a stochastic dynamics (e.g., its stock value) via an SDE:
\bea
\label{stock}
Y_t=\int_0^t h(s, X_{s})ds+B^2_t, \q t\ge 0,
\eea
where $B^2$ is the noise from the market, which we assume is independent of $B^1$ (this is by no means
necessary, we can certainly consider the filtering problem with correlated noises).

Now let us assume that the volatility $\si$ in (\ref{fvalue})  is affected by the actions of a large number of investors,
and all can only
make decisions based on the information from the process $Y$. Therefore, similar to \cite{CD2} (or \cite{HMC}) we begin by considering
%as starting point a system of
$N$ individual investors, and assume that $i$-th investor's
%\, (1\le i\le N)$, whose
private state dynamics is of the form:
\bea
\label{playeri}
dU^i_t=\si^i(t, U^i_{\cd\wedge t}, \bar \n^N_t, \a^i_t)dB^{1,i}_t, \qq t\ge 0, \q 1\le i\le N,
\eea
where $B^{1,i}$'s are independent Brownian motions, and $\bar\n^N_t$ denotes the empirical conditional distribution of
$U=(U^1, \cds, U^N)$, given the (common) observation $Y=\{Y_t:t\ge0\}$, that is, $\bar\n^N_t\dfnn\frac1N\sum_{j=1}^N \d_{\hE[U^j_t|\cF^Y_t]}$, where $\d_x$ denotes the Dirac measure at $x$. More precisely, the notation in (\ref{playeri}) means (see, e.g., \cite{CD2}),
\bea
\label{empirical}
\si^i(t, U^i_{\cd\wedge t}, \bar \n^N_t, \a^i_t)&\dfnn&\int_{\hR}\tilde\si^i(t, U^i_{\cd\wedge t}, y, \a^i_t)\bar\n^N_t(dy) \nonumber\\
&=& \frac1N\sum_{j=1}^N\int_{\hR}\tilde\si^i(t, U^i_{\cd\wedge t}, y, \a^i_t)\d_{\hE[U^j_t|\cF^Y_t]}(dy)\\
 &=& \frac1N\sum_{j=1}^N\tilde\si^i(t, U^i_{\cd\wedge t},\hE[U^j_t|\cF^Y_t], \a^i_t). \nonumber
 \eea
Here, $\tilde\si^i$'s are the functions defined on appropriate (Euclidean) spaces.
%on both sides of the above by the

We now assume that each investor chooses an individual strategy to minimize
the cost; the cost functional of the $i$-th agent is of the form:
\bea
\label{costi}
J^i(\a^i)\dfnn\hE\Big\{\Phi^i(U^i_T)+\int_0^TL^i(t,U^i_{\cd\wedge t}, \bar\n^N_t, \a^i_t)dt\Big\},\qq 1\le i\le N,
\eea
Following the argument of Lasry and Lions \cite{LasryLions} (see also \cite{CD2, CD3, CD5, CZ, HMC}), if we assume that the game is {\it symmetric}, i.e., $\tilde\si^i=\tilde\si,\ L^{i}$ and $\Phi^{i}=\Phi$ are independent of $i$, and that the number of investors $N$ converges to $+\infty$, then under suitable technical conditions, one could
find (approximate) Nash equilibriums through a limiting dynamics, and assign
a representative investor the unified strategy $\a$, determined by a {\it conditional} McKean-Vlasov type SDE
\bea
\label{McKeanV}
dX_t=\si(t, X_{\cd\wedge t}, \m_t, \a_t)dB^1_t,  \q t\ge 0,
\eea
where $\m$ is the conditional distribution of $X_t$ given $\cF^Y_t$, and
$$ \si(t, X_{\cd\wedge t}, \m_t, u_t)\dfnn\int \si(t, X_{\cd\wedge t}, y, u_t)\m_t(dy)=\hE\{\si(t, \f_{\cd\wedge t}, \hE[X_t|\cF^Y_t], u)\}|_{\f
=X, u=u_t}.
$$
Furthermore, the value function becomes, with similar notations,
\bea
\label{cost0}
V(x)=\inf_\a J(\a)\dfnn\hE\Big\{\Phi(X_T)+\int_0^TL(t,X_{\cd\wedge t}, \m_t, \a_t)dt\Big\}.
\eea
We note that (\ref{McKeanV}) and (\ref{cost0}), together with (\ref{stock}), form a stochastic control problem involving
CMFSDE dynamics and partial observations, as we are proposing.

The main objective of this paper is two-fold: We shall first study the exact meaning as well as the well-posedness of the dynamics, and then investigate the Stochastic Maximum Principle for the corresponding stochastic control problem. For the wellposedness of (\ref{dynamics}) we shall use a scheme that combines the idea of \cite{CD1} and the techniques of nonlinear filtering, and prove the existence and uniqueness of the solution to SDE (\ref{McKeanV}) via Schauder's fixed point theorem on $\sP_2(\O)$, the space of probability measures with finite second moment, endowed with the 2-Wasserstein metric. We note that the
important elements in this argument include the so-called {\it reference probability space} that is often seen in the nonlinear
filtering theory  and the Kallianpur-Striebel formula (cf. e.g., \cite{Bens, Zeit}), which enable us to define the solution mapping.

Our next task is to prove Pontryagin's Maximum Principle for our stochastic control problem. The main idea is similar to
earlier works of the first two authors (\cite{BLP, Li}), with some significant modifications. In particular, since in the present case
the control problem can only be carried out in a weak form, due to the lack of strong solution of CMFSDE, the existence of
the common reference probability space is essential. Consequently, extra efforts are needed to overcome the complexity caused
by the change of probability measures, which, together with the path-dependent nature of the underlying dynamic system,
makes even the first order adjoint equation more complicated than the traditional ones. To the best of our knowledge, the resulting mean-field backward SDE is new.

The paper is organized as follows. In Section 2 we provide all the necessary preparations, including some known
facts of nonlinear filtering. In Sections 3 and 4 we prove the well-posedness of the partially observable dynamics. In Section 5 we
introduce the stochastic control problem, and in Section 6 we study the variational equations and give some important
estimates. Finally, in Section 7 we prove the Pontryagin maximum principle.

\section{Preliminaries}
\setcounter{equation}{0}

Throughout this paper we consider the {\it canonical  space} $(\O, \cF)$, where $\O\dfnn \hC_0([0,\infty);\hR^{2d})= \{\o\in \hC([0,\infty);\hR^{2d}): \o_0 = {\bf 0}\}$, and $\cF$ be its topological $\si$-field.
%Let $\bar B_t(\o)=\o_t$, $\o\in\O$, be the canonical process, and
Let $\hF=\{\cF_t\}_{t\ge 0}$ be the natural filtration on $\O$, that is,  for each $t\ge 0$,
$\cF_t$ is the topological $\si$-field of the space $\O_t\dfnn \{\o(\cd\wedge t): \o\in\O\}$.
%generated by $\bar B$.
%coordinate process on $\hC_0([0,\infty);\hR^{2d})$a complete probability space $(\O, \cF, \hP)$ on which is
%defined a $2d$-dimensional Brownian motion $B=\{B_t\}_{t\ge0}$. We often assume that
%$\O$ is the .
For simplicity, throughout this paper we assume $d=1$, and that all the processes are 1-dimensional, although
the higher dimensional cases can be argued similarly without substantial difficulties.
Furthermore, we let $\sP(\O)$ denote the space
of all probability measures on $(\O, \cF)$, and for each $\hP\in \sP(\O)$, we assume that $\hF$ is $\hP$-augmented so that the filtered probability space $(\O, \cF, \hP; \hF)$
% be the  Wiener measure on $(\O,\cF)$.
%so that the canonical process $\bar B=(\bar B^1, \bar B^2)$ is a $\hQ$-Brownian motion.
%We assume that
satisfies the {\it usual hypotheses}.

Next, for given $T>0$ we denote $\hC_T=\hC([0,T])$ endowed by the supremum norm $\|\cdot\|_{\hC_T}$, and  let $\sB(\hC_T)$ be its topological $\si$-field. Consider now the space of all probability measures on $(\hC_T, \sB(\hC_T))$, denoted by $\sP(\hC_T)$, and for $p\ge 1$ we let $\sP_p(\hC_T)
\subseteq \sP(\hC_T)$ be those that have finite $p$-th moment. We recall that the {\it $p$-Wasserstein metric} on
$\sP_p(\hC_T)$ is defined as a mapping
$W_p:\sP_p(\hC_T)\times\sP_p(\hC_T)\mapsto \hR_+$ such that, for all $\m, \n\in\sP_p(\hC_T)$,
\bea
\label{Wp-metric}
W_p(\m,\n) \dfnn \inf\{(\int_{\hC^2_T}\|x-y\|^p_{\hC_T}\pi(dx,dy))^{\frac1p}: \pi\in \sP_p(\hC^2_T)~\mbox{with marginals $\m$ and $\n$}\}.
\eea

In this paper we shall use the 2-Wasserstein metric $W_2$, and abbreviate $(\sP_2(\hC_T), W_2)$ by $\sP_2(\hC_T)$. Since $\hC_T$ is a separable Banach space, it is known that $\sP_2(\hC_T)$ is a separable and complete metric space. Furthermore, it is known that (cf. e.g., \cite{Villani}), for $\mu_n,\mu\in{\sP}_2(\hC_T)$,
\bea\begin{array}{lcl}
\label{Wasserstein}
\lim_{n\to \infty}W_2(\mu_n,\mu)= 0 & \Longleftrightarrow&\, \, \mu_n\limw \mu ~\mbox{ in $\sP_2(\hC_T)$ and, as } N\rightarrow +\infty,\\
& & ~ \displaystyle\sup_n \int_{\Omega} \|\f\|^2_{\hC_T}I\{\|\f\|_{\hC_T}\ge N\}\mu_n(d\f)\rightarrow 0.
\end{array}
\eea

 Next, for any $\hP\in\sP(\O)$, $p,q\ge 1$, any sub-filtration $\hG\subseteq \hF$, and any Banach space $\hX$, we denote
$L^p(\hP; \hX)$ to be all $\hX$-valued $L^p$-random variables under $\hP$. In particular, we denote by $L^p(\hP;\hR)$ to be
all real valued $L^p$-random variables under $\hP$. Further, we denote by $L^p_{\hG}(\hP; L^q([0,T]))$  the $L^p$-space of all $\hG$-adapted processes $\eta$, such that
\bea
\label{LpqPnorm}
\| \eta\|_{p,q,\hP}\dfnn \Big\{\hE^{\hP}\Big[\int_0^T|\eta_t|^qdt\Big]^{p/q}\Big\}^{1/p}<\infty.
%, \qq \forall \eta\in L^p_{\hG}(\hP;L^q([0,T])).
\eea
If $p=q$, we simply write  $L^p_\hG(\hP;[0,T])\dfnn L^p_\hG(\hP; L^p([0,T]))$.  Finally,
we define $L^{\infty-}_{\hG}(\hP;[0,T])\dfnn \bigcap_{p>1} L^p_{\hG}(\hP;[0,T])$ and $\scL^{\infty-}_\hG(\hP; \hC_T)
\dfnn \bigcap_{p>1} L^p_\hG(\hP;\hC_T)$, where $L^p_\hG(\hP; \hC_T)$ is the space of all continuous, $\hF$-adapted,
processes $\xi=\{\xi_t\}$ such that $\|\xi\|_{\hC_T}\in L^p(\hP;\hR)$.
We will often drop ``$\hP$" from the subscript/superscript when the context is clear.

%\subsection{Formulation of the control problem}

%\ms
%{\bf The State-observation Dynamics.} \,
We now give a more precise description of
%our stochastic control problem. We first reformulate
the SDEs (\ref{dynamics}), in terms of the
standard McKean-Vlasov SDE. Again we consider only the case $b=0$, and we assume further that $\hat\si=1$\ in (\ref{dynamics}) for simplicity.

We begin by introducing some notations. Let $X$ be the state process and $Y$ the observation process, defined on $(\O, \cF, \hP)$, for some $\hP\in\sP(\O)$. We denote the ``filtered" state process by $U^{X|Y}_t=\hE^\hP [X_t|\cF^Y_t]$, $t\ge0$. Since (as we show in Lemma \ref{contiU} below) the process
$U^{X|Y}$ is continuous,  we denote its law
under $\hP $ on $\hC_T$ by $\m^{X|Y}=\hP \circ [U^{X|Y}]^{-1}\in\sP(\hC_T)$. Next, let $P_t(\f)=\f(t)$,
$\f\in \hC_T$, $t\ge0$, be the projection mapping, and define $\m^{X|Y}_t=\m^{X|Y}\circ {P_t}^{-1}$. Then, for any $\f\in \hC_T$, and
$u\in \hR$,
we can write
%consider the SDE:
%we can write
$$  \hE [\si(t, \f_{\cd\wedge t}, \hE [X_{t}|\cF^Y_t], u)]=\int\si(t, \f_{\cd\wedge t}, y, u)\m^{X|Y}_t(dy)
\dfnn\si(t, \f_{\cd\wedge t}, \m^{X|Y}_t, u).
$$

We should note that since the dynamics $X$ is non-observable, the decision of the controller
can only be made based on the information observed from the process $Y$. Therefore, it is reasonable
to assume that the control process $u$ is $\hF^Y=\{\cF^Y_t\}_{t\ge 0}$ adapted (or progressively measurable).
We should remark that, for a given such control, it is by no means clear that the state-observation SDEs will
have a strong solution on a prescribed probability space, as we shall see from our well-posedness result in the next sections. We therefore consider a ``weak formulation" which we now describe. Consider the pairs $( \hP, u)$, where $\hP\in\sP(\O)$, $u \in L^2_{\hF}(\hP;[0,T])$, such that the following SDEs are well-defined:
\bea
\label{controlsys}
X_t&=&x+\int_0^t\hE^{\hP}[\si(s, \f_{\cd\wedge s}, \hE^{\hP}[X_s|\cF^Y_s], z)]\Big|_{\f=X, z=u_s}
dB^1_s\\
&=&
x+\int_0^t\int_{\hR}\si(s, X_{\cd\wedge s}, y, u_s)\m_s(dy)dB^1_s=x+\int_0^t\si(s, X_{\cd\wedge s}, \m_s, u_s)dB^1_s,
\nonumber\\
\label{observation}
Y_t&=&\int_0^th(s, X_{ s})ds + B^2_t, \qq t \ge 0,
\eea
where $(B^1,B^2)$ is a standard 2-$d$ Brownian motion under $\hP$, and $\m_t(\cd)\dfnn \hP\circ \hE^{\hP}[X_t|\cF^Y_t]^{-1}(\cd)$ is the distribution, under $\hP$, of the conditional expectation of
$X_t$, given $\cF^Y_t$. We note that we {\it do not} require that the solution to (\ref{controlsys}) and (\ref{observation}) (or probability
$\hP$ for given $u$) be unique(!).
%; and the observation process $Y$ satisfies the SDEwhere .
%We now define the set of {\it admissible controls}.
Now let $U$ be a convex subset of $\hR^k$. For simplicity, assume $k=1$.
\begin{defn}
\label{admissible}
A pair $(\hP, u)\in \sP(\O)\times L^2_{\hF}(\hP;[0,T])$ is called an ``admissible control" if

{\rm (i)} $u_t\in U$, for all $t\in[0, T]$, and $B=(B^1,B^2)$ is a $(\hF,\hP)$-Brownian motion;

\ms
{\rm (ii)} There exist processes $(X, Y)\in L^2_{\hF}(\hP; [0, T])$ satisfying SDEs (\ref{controlsys}) and (\ref{observation}); and

\ms
{\rm (iii)} $u\in L^{\infty-}_{\hF^Y}(\hP; [0,T])$.
\qed
\end{defn}
We shall denote the set of all admissible controls by $\sU_{ad}$. For simplicity, we often write $u\in\sU_{ad}$, and
denote the associated probability measure(s) $\hP$ by  $\hP^u$, for $u\in\sU_{ad}$.

\begin{rem}
\label{rem1}
{\rm As we will shall see later, under our standing assumptions to every control $u\in\sU_{ad}$ there is only one probability measure $\hP^{u}$ associated.
We should note, however,  that
% the probability measure $\hP^u$ varies with the control and,
unlike the traditional filtering problem, the main difficulty of SDE (\ref{controlsys})-(\ref{observation}) lies in
the mutual dependence between the solution pair $X^u$ and $Y$, via
%``closed-looped" nature, that is, the dependence of
the law of conditional expectation
$\m^u_t=\hP^u\circ \hE^{\hP^u}[X^u_t|\cF^Y_t]^{-1}$ in the coefficients.
% of $X^u$ knowing $Y$.
Moreover, the requirement that $u$ is $\hF^Y$-adapted
adds an additional seemingly ``circular" nature to the problem. Thus, the well-posedness of the problem is far from obvious, and will be the main subject of \S3.\qed}
\end{rem}

% which will also have the ``conditional mean-field" type.{\bf The Reference Measure.}

We note that under the weak formulation the state-observation processes $(X^u, Y)$ are often defined on different probability spaces. To facilitate our discussion we shall designate  a common space on which all the controlled dynamics can be evaluated. In light of the nonlinear filtering theory, we make the following assumption.
\begin{assum}
\label{Assump2}
There exists a probability measure $\hQ^0$ on $(\O, \cF)$, such that, under $\hQ^0$, $(B^1, Y)$ is a 2-dimensional Brownian motion, where
$Y$ is the observation process.
\qed
\end{assum}
We note that the probability measure $\hQ^0$ is commonly known as the ``reference probability measure" in nonlinear filtering theory. The existence of such measure  can be argued once the existence of the weak solution of (\ref{controlsys})-(\ref{observation}) is known.
% as we will elaborate in \S3.
%further, which allows to consider $Y$ independent of the control $u$. This will be the subject of Section 3.
%Independently of what will be done in Section under adequate assumptions on the coefficients,
%Let us assume for the moment that such a ``reference measure'' $\hQ^0$ exists, so that
%$(B^1,Y)$ is a $\hQ^0$-Brownian motion.  Then the relation between $\hQ^{0}$ and $\hP^{u}$ can be described  with the help of the
% Girsanov Theorem.
Indeed, suppose that $u\in \sU_{ad}$ and $\hP^u\in \sP(\O)$ is the associated probability such that the SDEs (\ref{controlsys}) and (\ref{observation}) have a solution $(X^u,Y)$ on $(\O, \cF, \hP^u)$.
%measure and $X=X^u$ be the corresponding state process such that
%is satisfied.
Consider the following SDE:
\bea
\label{barL}
\bar L_t=1-\int_0^t h(s, X^u_{s}) \bar L_sdB^2_s=1+\int_0^t \bar L_s dZ^u_s,
\eea
where $Z^u_t=-\int_0^t h(s, X^u_{ s})dB^2_s$. We denote its solution by $\bar L^u$. Then, under appropriate conditions on $h$,
both $Z^u$ and $\bar L^u$ are $\hP^u$-martingales, and $\bar L^u$ is the stochastic exponential:
\bea
\label{barLexp}
\bar L^u_t=\exp\Big\{Z^u_t-\frac12 \la Z^u\ra_t\Big\}=\exp\Big\{-\int_0^th(s, X^u_{s})dB^2_s-\frac12 \int_0^t
|h(s, X^u_{s})|^2ds\Big\}.
\eea
Thus, the Girsanov Theorem suggests that $d\hQ^0=\bar{L}_T^{u}d\hP^{u}$ defines a new probability measure $\hQ^0$ under which
 $(B^1,Y)$ is a Brownian motion, hence a ``reference measure".

The essence of Assumption \ref{Assump2} is, therefore, to assign a {\it prior distribution} on the observation process $Y$ {\it before} the well-posedness of the control system is established. In fact, with such an assumption one can begin by assuming that $(B^1, Y)$
is the canonical process (i.e., $(B^1_t, Y_t)(\o)=\o(t)$, $\o\in\O$) and $\hQ^0$ the Wiener measure on $(\O, \cF)$, and then
%existence of the reference measure  $\hQ^0$ is actually crucial for us to
proceed to prove the existence of the weak solution of the system (\ref{controlsys}) and (\ref{observation}). This
scheme will be carried out in details in \S3.

Continuing with our control problem,  for any $u\in\sU_{ad}$, we define the {\it cost functional} by
\bea
\label{cost}
J(t, x; u)&\dfnn& \hE^{\hQ^0}\Big\{\int_t^T f(s, X^u_{\cd\wedge s}, \m^u_s, u_s)ds+\Phi(X^u_T, \m^u_T)\Big\}\nonumber\\
&=&\hE^{\hQ^0}\Big\{\int_t^T\hE^{\hP^u}[f(s, \f_{\cd\wedge s}, \hE^{\hP^u}[X^u_s|\cF^Y_s], u)]\Big|_{\f=X^u, u=u_s}ds\\
&&+\hE^{\hP^u}[\Phi(x,\hE^{\hP^u}[X^u_T|\cF^Y_T])]\Big|_{x=X^u_T}\Big\},\nonumber
\eea
and we denote the  value function as
\bea
\label{value}
V(t, x)\dfnn \inf_{u\in\sU_{ad}} J(t, x; u).
\eea

We shall make use of the following {\it Standing Assumptions} on the coefficients.
\begin{assum}
\label{Assum1}
{\rm (i)} The mappings $(t, \f,x, y, z)\mapsto \sigma(t, \f_{\cd\wedge t},y, z)$, $h(t,x)$,  $f(t, \f_{\cd\wedge t}, y, z)$, and $\Phi(x, y)$ are bounded and continuous, for $(t,\f,x,y, z)\in [0,T]\times\hC_T\times\hR\times \hR\times U$;

\ss
{\rm (ii)} The partial derivatives $\pa_y\si$, $\pa_z\si$, $\pa_y f$, $\pa_z f$, $\pa_x h$, $\pa_x\Phi$, $\pa_y\Phi$ are bounded and  continuous, for $(\f, x, y,z)\in \hC_T\times \hR\times\hR\times U$, uniformly in $t\in[0,T]$;

\ss
{\rm (iii)} The mappings $\f\mapsto \si(t, \f_{\cd\wedge t}, y, z), f(t, \f_{\cd\wedge t}, y, z)$, as functionals from $\hC_T$ to $\hR$, are
Fr\'echet differentiable. Furthermore, there exists a family of measures $\{\ell(t, \cd)\}|_{t\in[0,T]}$, satisfying $0\le \int_0^T\ell(t,ds)\le C$, for all $t\in[0,T]$, such that both derivatives, denoted by $D_\f\si=D_\f\si(t, \f_{\cd\wedge t}, y, z)$ and $D_\f f=D_\f f(t, \f_{\cd\wedge t}, y, z)$, respectively,  satisfy
\bea
\label{DsiDf}
|D_\f\si(t, \f_{\cd\wedge t}, y, z)(\psi)|+|D_\f f(t, \f_{\cd\wedge t}, y, z)(\psi)|\le \int_0^T|\psi(s)|\ell(t, ds), \ \ \psi\in\hC_T,
\eea
uniformly in $(t, \f, y,z)$;

\ss
{\rm (iv)} The mapping $y\mapsto y\pa_y\si(t, \f_{\cd\wedge t}, y, z)$ is uniformly bounded, uniformly in $(t, \f,z)$;

\ss
{\rm (v)} The mapping $x\mapsto x\pa_x h(t, x)$ is bounded, uniformly in $(t,x)\in[0,T]\times\hR$;

\ss

{\rm (vi)} The mappings $x\mapsto x h(t,x), x^2\pa_x h(t,x)$ are bounded, uniformly in $(t,x)\in[0,T]\times\hR$.
\qed
\end{assum}

We note that some of the assumptions above are merely
technical and can be improved, but we prefer not to dwell on such technicalities and focus on the
main ideas instead.

\begin{rem}
\label{remark3}
{\rm Note that if $(t, \f, y, z)\mapsto\phi(t, \f_{\cd\wedge t}, y, z)$ is a function defined on $[0,T]\times\hC_T\times \hR\times \hR$
satisfying Assumption \ref{Assum1}-(i), (ii), then for any $\m\in\sP_2(\hC_T)$, we can define a function on the space $[0,T]\times\O\times\hC_T\times \sP_2(\hC_T)\times U$:
\bea
\label{barphi}
\bar\phi(t,\o, \f_{\cd\wedge t}, \m_t, z)\dfnn \int_\hR\phi(t, \f_{\cd\wedge t}, y, z)\m_t(dy),
\eea
where $\m_t=\m\circ P_t^{-1}$ and $P_t(\f)\dfnn\f(t)$, $(t, \f)\in [0,T]\times\hC_T$. Then, $\bar \phi$ must satisfy
the following Lipschitz condition:
\bea
\label{Lip}
|\bar\phi(t, \f^1_{\cd\wedge t}, \m^1_t, z^1)-\bar\phi(t, \f^2_{\cd\wedge t}, \m^2_t, z^2)|\le  K\Big\{\|\f^1-
\f^2\|_{\hC_t}+ W_2(\m^1, \m^2)+|z^1-z^2|\Big\},
\eea
where $\|\cd\|_{\hC_t}$ is the sup-norm on $\hC([0,t])$ and $W_2(\cd, \cd)$ is the 2-Wasserstein metric.
\qed}
\end{rem}

\begin{rem}
\label{remark6}
{\rm
The Fr\'echet derivatives $D_\f\si$ and $D_\f f$ by definition belong to $\hC_T^*\dfnn \sM[0,T]$, the
space of all finite signed Borel measures on $[0,T]$, endowed with the total variation norm $|\cd|_{TV}$ (with a slight abuse
of notation, we still denote it by $|\cd|$). Thus the Assumption \ref{Assum1}-(iii)
amounts to saying that, as measures,
\bea
\label{dominate}
|D_\f\si(t, \f_{\cd\wedge t}, y, z)(ds)|+|D_\f f(t, \f_{\cd\wedge t}, y, z)(ds)|\le \ell(t, ds), \q \forall (t, \f, y, z).
\eea
This inequality will be crucial in our discussion in Section 7.
\qed}
\end{rem}

%\subsection

To end this section we recall some basic facts in nonlinear filtering theory, adapted to our situation.
We begin by considering the inverse Girsanov kernel of $\bar L^u$ defined by (\ref{barLexp}):
%We note that in what follows we will use more often the L^u_t\dfnn
\bea
\label{Lexp}
L^u_t\dfnn [\bar L^u_t]^{-1}=\exp\Big\{\int_0^th(s, X^u_{s})dY_s-\frac12 \int_0^t
|h(s, X^u_{s})|^2ds\Big\}, \ \ t\in [0,T].
\eea
Then $L^u$ is a $\hQ^0$-martingale, $d\hP^u=L^u_Td\hQ^0$, and $L^u$
%We now start from the reference measure $\hQ^0$ . For fixed control $(\hP, u)$, we denote $L=L^u$ by (\ref{Lexp}). Then $L$
satisfies the following SDE on $(\O, \cF, \hQ^0)$:
 \bea
\label{L}
L_t=1+\int_0^t h(s, X_{s})  L_sdY_s, \qq t\in[0,T].
\eea
Let us now denote $L=L^u$ for simplicity. An important ingredient that we are going to use frequently is the SDEs known as the {\it Kushner-Stratonovic} or {\it Fujisaki-Kallianpur-Kunita} (FKK) equation for the  ``normalized conditional probability". Let us denote
\bea
\label{S}
S_t\dfnn \hE^{\hQ^0}[L_tX_t|\cF^Y_t],  \q S^{0}_t\dfnn \hE^{\hQ^0}[L_t|\cF^Y_t], \q t\ge 0.
\eea
Since under $\hQ^0$ the process $(B^1,Y)$ is a Brownian motion,  the $\si$-field $\cF^Y_{t, T}$ and $\cF^Y_t\vee \cF^{B^1}_t$
are independent, where $\cF^Y_{t, T}\dfnn \si\{Y_r-Y_t: t\le r\le T\}$.
It is standard to show that (in light of (\ref{L})) $S$ and $S^0$ satisfy the following SDEs:
\bea
\label{Sn0}
S^{0}_t=1+\int_0^t\hE^{\hQ^0}[ h(s, X_{s})L_s|\cF^Y_s]dY_s, \q t\ge 0.
\eea
and
\bea
\label{SDES}
S_t=x+\int_0^t \hE^{\hQ^0}[L_sX_sh(s, X_{s})|\cF^Y_s]dY_s, \q t\ge 0.
\eea

Furthermore,  let  $U_t\dfnn \hE^{\hP^u}[X_t|\cF^Y_t]$, $t\ge 0$. Then, by the Bayes formula (also known as the Kallianpur-Striebel formula, see, e.g., \cite{Bens}) we have
\bea
\label{Bayes}
U_t=\frac{\hE^{\hQ^0}[L_tX_t|\cF^Y_t]}{\hE^{\hQ^0}[L_t|\cF^Y_t]}=
\frac{S_t}{S^0_t}, \q t\ge 0, \q \hQ^0\mbox{-a.s.}
\eea
 A simple application of It\^o's formula and some direct computation then lead to the following FKK equation:
\bea
\label{FKK}
dU_t&=&\Big\{\hE^{\hP^u}[X_th(t, X_{t})|\cF^Y_t]-\hE^{\hP^u}[X_t|\cF^Y_t]\hE^{\hP^u}
[h(t,X_{ t})|\cF^Y_t]\Big\}dY_t\\
&&+\Big\{\hE^{\hP^u}[X_t|\cF^Y_t]\big\{\hE^{\hP^u}[h(t,X_{ t})|\cF^Y_t]\big\}^2-\hE^{\hP^u}[X_th(t,X_{ t})|\cF^Y_t]\hE^{\hP^u}[h(t,X_{t})|\cF^Y_t]\Big\}dt.\nonumber
\eea
In fact, one can easily show that
\bea
\label{ZakaivsFKK}
S_t&=&U_t\exp\Big\{\int_0^t\hE^{\hP^u}[h(s,X_{s})|\cF^Y_s]dY_s-\frac12\int_0^t\hE^{\hP^u}[h(s,X_{ s})|\cF^Y_s]^2ds\Big\}.
\eea

\section{Well-posedness of the State-Observation Dynamics}
\setcounter{equation}{0}

In this and next sections we investigate the well-posedness of the controlled state-observation system (\ref{controlsys}) and (\ref{observation}). More precisely, we shall argue that the admissible control set $\sU_{ad}$, defined by Definition \ref{admissible}, is not empty. We first note that, for a fixed $\hP\in\sP(\O)$ and
% such that $\hP^0\sim \hQ^0$ ,
$u\in L^{\infty-}_{\hF^Y}(\hP, [0,T])$, if we define
%along with a  for the associated $\hP^{u}$, which we suppose to be defined over $(\O,\hF)=(\hC_T,\sB(\hC_T))$.the weak sense.
\bea
\label{sigma}
\phi^u(t,\o, \f_{\cd\wedge t}, \m_t)\dfnn \int_\hR\phi(t, \f_{\cd\wedge t}, y, u_t(\o))\m_t(dy),
\eea
where $\phi=b, \si$, then we can write the control-observation system (\ref{controlsys}) and (\ref{observation}) as
%for the given control $u\in\sU_{ad}$,first give
a slightly more generic form (denoting $b^u=b$ and $\si^u=\si$ for simplicity):
%
%, which will include (\ref{controlsys}) and (\ref{observation})
%as special case. We note that SDE (\ref{SDE})
%is exactly  in which
%the coefficients $b$ and $\si$ should be understood as the random fields defined by
%For this
%
%We recall from  that the system (\ref{controlsys}) and (\ref{observation}) is satisfied with respect to some  .
%
%
%We shall consider the following system of SDEs:
\bea
\label{SDE}
\left\{\ba{lll}
\dis X_t=x+\int_0^tb(s,\cd, X_{\cd\wedge s}, \m^{X|Y}_s)ds+\int_0^t\si(s,\cd, X_{\cd\wedge s}, \m^{X|Y}_s)dB^1_s;\ms\\
\dis Y_t=\int_0^t h(s, X_{s})ds+B^2_t,
\ea\right. \qq t\ge 0,
\eea
where $B=(B^1,B^2)$ is a $\hP$-Brownian motion, and $\m^{X|Y}_t=\hP\circ [\hE^{\hP}[X_t|\cF^Y_t]]^{-1}$. Our task is to prove the well-posedness of SDE (\ref{SDE}) in a {\it weak} sense (i.e., including the existence of
the probability measure $\hP$(!)). In light of Remark \ref{remark3},
we shall  assume that the coefficients $b$ and $\si$
in (\ref{SDE}) satisfy the following assumptions that are slightly weaker than Assumption \ref{Assum1}, but sufficient for our
purpose in this section.

\begin{assum}
\label{Assum2}
The coefficients $b, \si: [0, T]\times \hC_T\times \sP_2(\hC_T)\mapsto \hR$ enjoy the following properties:

{\rm(i)} For fixed $(\f, \m)\in \hC_T\times \sP_2(\hC_T)$, the mapping $(t,\o)\mapsto (b, \si)(t,\o, \f, \m)$ is an $\hF$-progressively
measurable process;

{\rm(ii)} For fixed $t\in [0,T]$, and $\hQ^0$-a.e. $\o\in\O$, there exists $K>0$, independent of $(t,\o)$, such that for all
$(\f^1,\m^1), (\f^2,\m^2)\in \hC_T\times \sP_2(\hC_T)$, it holds that
\bea
\label{Lip2}
|\phi(t,\o, \f^1_{\cd\wedge t}, \m^1_t) -\phi(t,\o, \f^2_{\cd\wedge t}, \m^2_t)|\le K(\sup_{t\in[0,T]}|\f^1_t-\f^2_t|+W_2(\m^1, \m^2)),
\eea
for $\phi=b, \si$, respectively.
\qed
\end{assum}

In the rest of the section we shall still assume $b=0$, as it does not add extra difficulties.
Now assume that $(X,Y)$ satisfies (\ref{SDE}) under $\hP$, and
%For a continuous $\hF$-adapted process $X$
let us denote
$U^{X|Y}_t\dfnn  \hE^{\hP}[X_t|\cF^Y_t]$, $t\ge0$. (We note that $U^{X|Y}$ should be understood as the ``optional
projection" of $X$ onto $\hF^Y$!) We first check that $U^{X|Y}$ is indeed a continuous process.
\begin{lem}
\label{contiU}
Assume that Assumption \ref{Assum1} holds. Then  $U^{X|Y}$ admits a continuous version.
\end{lem}

{\it Proof.}  First note that $\hP\sim \hQ^0$, and $X$ has continuous paths, $\hP$-a.s. By Bayes formula (\ref{Bayes}) we can write
$ U^{X|Y}_t=\frac{\hE^{\hQ^0}[L_tX_t|\cF^Y_t]}{\hE^{\hQ^0}[L_t|\cF^Y_t]}=\frac{S_t}{S^0_t}$, where $S^0$ and $S$ satisfy
(\ref{Sn0}) and (\ref{SDES}), respectively, and $L$ satisfies (\ref{L}). Clearly, the representations (\ref{Sn0}) and (\ref{SDES}) indicate that both $S^0$ and $S$ have continuous paths, thus $U^{X|Y}$ must have a continuous version.
\qed

\ms
We now define $\m^{X|Y}(\cd)=\hP\circ [U^{X|Y}]^{-1}(\cd)$, and $\m^{X|Y}_t(\cd)=\hP\circ [U^{X|Y}_t]^{-1}(\cd)$, for any $t\ge0$.
Lemma \ref{contiU} then implies that
$\m^{X|Y}\in\sP_2(\hC_T)$, justifying the definition of SDE (\ref{SDE}). In what follows when the context is clear,  we shall omit ``$X|Y$" from the superscript.

We note that the special circular nature of SDE (\ref{SDE}) between its solution and
its law of the conditional expectation (whence the underlying probability) makes it necessary to specify the meaning of a solution.
We have the following definition.
\begin{defn}[Weak Solution]
\label{sol}
An eight-tuple $(\O, \cF, \hP, \hF, X, Y, B^1,B^2)$  is called a solution
to the filtering equation (\ref{SDE}) if

{\rm(i)}  $(\O, \cF)$ is the canonical space, $\hP\in \sP(\O)$, and $\hF$ is the canonical filtration;

{\rm(ii)} $(B^1, B^2)$ is a 2-dimensional $\hF$-Brownian motion under $\hP$;

{\rm(iii)} $(X, Y)$ is an $\hF$-adapted continuous process such that (\ref{SDE}) holds for all $t\in[0,T]$, $\hP$-almost surely.
\end{defn}

To prove the well-posedness we shall use a  generalized version of the Schauder Fixed Point Theorem (see Cauty \cite{Cauty}, or
a recent generalization
%of the proof of Schauder's conjecture in see \cite{Cauty} we refer the reader to
in \cite{Cauty2}). To this end we consider the
following subset of $\sP_2(\hC_T)$:
\bea
\label{cE}
\sE\dfnn \Big\{\m\in \sP_2(\hC_T)
%\hC([0,T], \sP(\hR))
\big|\sup_{t\in [0,T]}\int_{\hR}|y|^4\m_t(dy)<\infty\Big\}.
\eea
In the above $\m_t=\m\circ {P_t}^{-1}\in\sP_2(\hR)$, and  $P_t(\f)=\f(t)$, $\f\in \O$, is the projection mapping.
Clearly, $\sE$ is a convex subset of $\sP_2(\hC_T)$.

 We now construct a mapping $\sT:\sE\mapsto \sE$, whose fixed point, if exists, would give a solution to the SDE (\ref{SDE}).
We shall begin with %%we recall Assumption \ref{Assum2},
 the reference probability space $(\O, \cF, \hQ^0)$, thanks to Assumption \ref{Assump2}, then $(B^1, Y)$ is a $\hQ^0$-Brownian motion.
We may assume without loss of generality that $(B^1,Y)$ is the canonical process,
% on $(\O,\cF)=(\hC_T,\sB(\hC_T))$
and $\hQ^0$ is the Wiener measure.
%In particular, we have that $(B^1,Y)$ is an $(\hF,\hQ^0)$- Brownian motion.

For any $\m\in \sE$ we consider the SDE on the space $(\O, \cF, \hQ^0)$:
\bea
\label{fileq2}
X_t=x+\int_0^t \si(s, \cd, X_{\cd\wedge s}, \m_s)dB^1_s, \q t\ge 0.
%Y_t&=&\int_0^t h(s, X_{\cd\wedge s})ds +W_t. \nonumber
\eea
Note that as the distribution $\m$ is given, (\ref{fileq2}) is an ``open-loop" SDE with ``functional Lipschitz" coefficient, thanks to
Assumption \ref{Assum2}. Thus,
%under the standard assumptions on the coefficients
there exists a unique (strong) solution to (\ref{fileq2}), which we denote by $X=X^\m$.

%We should note that given the state-observation pair $(X,Y)$,  such an assumption is necessary since we do not yet know that
%$(X, Y)$ satisfies the SDE (\ref{SDE}).
Now, using  $X^\m$ we  define the process $L^\m=\{L^\m_t\}_{t\ge0}$ as in (\ref{Lexp}) on probability space  $(\O, \cF, \hQ^0)$,
 and then we define the probability $d\hP^\m\dfnn L^\m_Td\hQ^0$. By the Kallianpur-Striebel formula (\ref{Bayes}) we can define a  process
\bea
\label{KS}
U^\m_t\dfnn\hE^{\hP^\m}[X^\m_t|\cF^Y_t]=\frac{\hE^{\hQ^0}[L^\m_tX^\m_t|\cF^Y_t]}{\hE^{\hQ^0}[L^\m_t|\cF^Y_t]}=\frac{S^\m_t}{S^{\m,0}_t}, \q t\ge0,
\eea
where $S^\m_t\dfnn \hE^{\hQ^0}[L^\m_tX^\m_t|\cF^Y_T]$, $S^{\m,0}_t\dfnn \hE^{\hQ^0}[L^\m_t|\cF^Y_T]$, $t\ge 0$,
and then we denote
\bea
\label{Tnu}
\sT(\m)\dfnn \n^\m= \hP^\m\circ[U^\m]^{-1}\in\sP(\hC_T).
\eea

Our task is to show that the solution mapping $\sT: \m\mapsto \n^\m$ satisfies the desired assumptions for
Schauder's Fixed Point Theorem.

\begin{thm}
\label{compact}
The solution mapping $\sT:\sE\to \sP_2(\hC_T)$ enjoys the following properties:

%(1) $\sT$ is continuous on $\sP(\hC([0,T];\hR^d))$;

(1) $\sT(\sE)\subseteq \sE$;

(2) $\sT(\sE)$ is compact under 2-Wasserstein metric.

(3) $\sT:(\sE,W_1(\cdot,\cdot))\rightarrow (\sP_2(\hC_T),W_2(\cdot,\cdot))$ is continuous, i.e., whenever $\mu,\mu^n\in \sE,\, n\ge 1,$ is such that $W_1(\mu^n,\mu)\rightarrow 0$, we have that $W_2(\sT(\mu^n),\sT(\mu))\rightarrow 0.$
\end{thm}

We remark that an immediate consequence of (3) is that $\sT:\sE\rightarrow \sP_2(\hC_T)$ is continuous under both the 1- and the 2-Wasserstein metrics. Moreover, the compactness of $\sT(\sE)$ under the 2-Wasserstein metric stated in (2) implies that in the 1-Wasserstein metric.

{\it Proof.} (1) Given $\m\in \sE$ we need only show that
\bea
\label{moment4}
\sup_{t\in[0,T]}\int_{\hR}|y|^4\n^\m_t(dy)<\infty.
%, \qq \mbox{for all} ~t\in[0,T].
\eea
To see this we note that for $t\in[0,T]$, by Jensen's inequality,
\beaa
\int_{\hR}|y|^4\n^\m_t(dy)=\int_{\hR}|y|^4\hP^\m\circ [U^\m]^{-1}(dy)=\hE^{\hP^\m}[|\hE^{\hP^\m}[X^\m_t|\cF^Y_t]|^4]
\le \hE^{\hP^\m}[|X^\m_t|^4].
\eeaa
Since under $\hQ^0$, $B^1$ is also a Brownian motion, it is standard to argue that,  as $X^\m$ is the solution to
the SDE (\ref{fileq2}),  it holds that
\bea
\label{momentQ0}
\sup_{0\le t\le T} \hE^{\hQ^0}[|X^\m_t|^{2n}] \le C(1+|x|^{2n}), \qq \mbox{for all } n\in\hN.
\eea
%have finite moments of all orders.
Furthermore, noting that the process $L^\m$ is an
$L^2$-martingale under $\hQ^0$, we have
\beaa
\sup_{0\le t\le T}\int_{\hR^d}|y|^4\n^\m_t(dy)
%\hE^{\hP^\m}\Big[\sup_{0\le t\le T}\hE^{\hP^\m}[|X^\m_t|^4|\cF^Y_t]\Big]
&\le& \sup_{0\le t\le T}\hE^{\hP^\m}\Big[|X^\m_t|^4\Big]
=\sup_{0\le t\le T}\hE^{\hQ^0}\Big[L^\m_T|X^\m_t|^4\Big]\\
&\le&\big(\hE^{\hQ^0}[|L^\m_T|^2]\big)^{\frac12}\sup_{0\le t\le T}\hE^{\hQ^0}\Big[|X^\m_t|^8\Big]^{\frac12}
<\infty,
\eeaa
thanks to (\ref{momentQ0}). In other words, $\n^\m=\sT(\m)\in\sE$, proving (1).
%$$ X^\m_t=x+\int_0^t\int_{\hR^d}\si(s, \f_{\cd\wedge s}, y)\m_t(dy)\big|_{\f=X} dB^1_s, \q t\ge 0.$$

(2) We shall prove that for any sequence $\{\m^n_t\}\subseteq \sE$, there exists a subsequence,
denoted by $\{\m^{n}_t\}$ itself, such that  $\lim_{n\to\infty}\sT(\m^{n})= \n$  in 2-Wasserstein metric, for some
$\n\in  \sT(\sE)$.

In light of the equivalence relation (\ref{Wasserstein}), we shall first argue that the family $\{\sT(\m^n)\}_{n\ge1}$ is tight. To this end, recall that
\bea
\label{Un}
U^n_t=\hE^{\hP^n}[X^n_t|\cF^Y_t]=\frac{S^n_t}{S^{n,0}_t},
\eea
where $S^n_t\dfnn \hE^{\hQ^0}[L^n_tX^n_t|\cF^Y_t]$,
%$$
$S^{n,0}_t\dfnn \hE^{\hQ^0}[L^n_t|\cF^Y_t]$, $t\ge 0$, and $d\hP^n\dfnn L^n_Td\hQ^0$.
It then follows from the FKK  equation (\ref{FKK})  that
%% by integration by parts:
\bea
\label{FKK1}
dU^n_t&=&\big\{\hE^{\hP^n}[X^n_th(t, X^n_{t})|\cF^Y_t]-\hE^{\hP^n}[X^n_t|\cF^Y_t]\hE^{\hP^n}[h(t,X^n_{ t})
|\cF^Y_t]\big\}dY_t\\
&&+\big\{\hE^{\hP^n}[X^n_t|\cF^Y_t](\hE^{\hP^n}[h(t,X^n_{ t})|\cF^Y_t])^2-\hE^{\hP^n}[X^n_th(t, X^n_{ t})|\cF^Y_t]\hE^{\hP^n}[h(t, X^n_{ t})|\cF^Y_t]\big\}dt. \nonumber
\eea
Now denote $B^{2,n}_t\dfnn Y_t-\int_0^th(s,X^n_{\cd\wedge s})ds$. Then $(B^1, B^{2,n})$ is a 2-dimensional
standard $\hP^n$-Brownian motion. Furthermore, since $h$ is bounded, so is
$\hE^{\hP^n}[h(t, X^n_{\cd\wedge t})|\cF^Y_t]$.
We thus have the following estimate:
$$\begin{array}{lll}
\hE^{\hP^n}[|U^n_t-U^n_s|^4]&\le& C\hE^{\hP^n}\Big[\Big(\int_s^t\hE^{\hP^n}[|X^n_s|^2|\cF^Y_s]ds\Big)^2
\Big]\nonumber\\
&\le& C\hE^{\hP^n}\Big[\sup_{0\le s\le T}\big|\hE^{\hP^n}[|X^n_s|^2|\cF^Y_s]|^2\Big]|t-s|^2\nonumber\\
\end{array}$$
\bea \label{est1}
&\le& C\hE^{\hP^n}\Big[\sup_{0\le s\le T}\big|\hE^{\hP^n}[\sup_{0\le r\le T}|X^n_r|^2|\cF^Y_s]|^2\Big]|t-s|^2\\
&\le & C\hE^{\hP^n}\Big[\sup_{0\le s\le T}|X^n_s|^4\Big]|t-s|^2\le C|t-s|^2. \nonumber
\eea

Thus, as $U^n_0=x$, $n\geq 1$, the sequence of  continuous processes $\{U^n\}$ is relatively compact  (cf. e.g., Ethier-Kurtz \cite{EK}).
Therefore, the sequence of their laws $\{\sT(\m^n)\dfnn \hP^n\circ [U^n]^{-1}, n\ge1\}\subseteq \sP(\hC_T)$ is tight. Consequently,
we can find a subsequence, we may assume itself, that converges weakly to a limit $\n\in\sP_2(\hC_T)$.
Furthermore, for each $n\ge 1$, we apply the Jensen, Burkholder-Davis-Gundy, and H\"older inequalities
to get, with $\n^n\dfnn\sT(\m^n)$,
%$\sup_{n\ge 1}\int_{C_T}|y|^4\mu^n(dy)<+\infty$, in order to be able to conclude that
%$W_2(\mu_n,\mu)\le W_4(\mu_n,\nu)\rightarrow 0$, as $n\rightarrow +\infty.$
%
\bea
\label{estnu}
\int_{\hC_T}\|\f\|^4_{\hC_T}\nu^n(d\f)&=&\hE^{\hP^n}[\|U^n\|^4_{\hC_T}]=\hE^{\hP^n}[\sup_{0\le t\le T}|\hE^{\hP^n}[X_t^n
|{\cal F}_t^n]|^4]\nonumber\\
&\le&
\hE^{\hP^n}\Big[\sup_{0\le t\le T} \hE^{\hP^n}\big[\sup_{0\le r\le T}|X_r^n||{\cal F}_t^n\big]^4\Big]\\
&\le& C\Big[\hE^{\hP^n}\big[\sup_{0\le r\le T}|X_r^n|^{6}\big]\Big]^{2/3}=C\Big[\hE^{\hQ_0}\big[L_T^n\sup_{0\le r\le T}|X_r^n|^{6}\big]\Big]^{2/3}\nonumber\\
&\le& C\Big[\hE^{\hQ_0}[(L_T^n)^4]\Big]^{1/6}\Big[\hE^{\hQ_0}[\sup_{0\le r\le T}|X_r^n|^8]\Big]^{1/2}<+\infty. \nonumber
\eea
But noting that $h$ is bounded, one deduces from (\ref{momentQ0}) that
\bea
\label{uniformbdd}
\sup_{n\ge 1}\int_{\hC_T}\|\f\|^4_{\hC_T}\n^n(d\f)<\infty,
\eea
and, thus,

\centerline{$\displaystyle\sup_{n\ge 1}\int_{\hC_T}\|\f\|^2_{\hC_T}I\{|\f\|_{\hC_T}\ge N \}\n^n(d\f)\rightarrow 0,
\mbox{ as }N\rightarrow +\infty.$}

\smallskip

\noindent This, together with the fact that $\n^n=\sT(\m^n)\limw \n$, implies that $W_2(\n^n, \n)\to 0$, and $\n\in \sE$, as $n\to\infty$, where $W_2(\cd,\cd)$ is the 2-Wasserstein metric on $\sP_2(\hC_T)$. This proves (2).

(3) We now check that the mapping  $\sT:(\sE,W_1(\cdot,\cdot))\rightarrow (\sP_2(\hC_T),W_2(\cdot,\cdot))$ is continuous.  To this end,
for each $\m\in\sE$, we consider the following SDE  on the probability
space $(\O, \cF, \hQ^0)$:
\bea
\label{XBLQ0}
\left\{\ba{lll}
dX_t=\si(t, X_{\cd\wedge t}, \m_t)dB^1_t, \q & X_0=x;\\
dB^2_t=dY_t-h(t, X_{t})dt, & B^2_0=0; \\
dL_t=h(t, X_{ t})L_t dY_t, & L_0=1.
\ea\right.
\eea
Now let $\{\m^n\}\subseteq \sE$ be any sequence such that $\m^n\to \m$, as $n\to\infty$, in the 1-Wasserstein metric, and denote by $(X^n, B^{n,2}, L^n)$ the corresponding solutions to (\ref{XBLQ0}). Define
$$\si^n(t, \o_{\cd\wedge t})\dfnn \si(t, \o_{\cd\wedge t}, \m^n_t), \q (t, \o)\in [0,T]\times \O.
$$
Then by Assumption \ref{Assum2}-(ii), the $\si^n$'s are functional Lipschitz deterministic functions, with Lipschitz constant independent
of $n$. This and standard SDE arguments lead to that, as $n\to\infty$,
\bea
\label{converge}
 \hE^{\hQ^0}\Big\{\sup_{0\le t\le T}|X^n_t-X_t|^p+\sup_{0\le t\le T}|L^n_t- L_t|^p\Big\}\to 0, \q  \mbox{in $L^p(\hQ^0)$, \q
 $p\ge 1$.}
\eea
We deduce that $U^n_t=\hE^{\hP^n}[X^n_t|\cF^Y_t]=S^n_t/S^{n,0}_t$ converges in probability under $\hQ^0$ to
$\frac{\hE^{\hQ^0}[L_tX_t|\cF^Y_t]}{\hE^{\hQ^0}[L_t|\cF_t^Y]}=\hE^{\hP}[X_t|\cF^Y_t]$, where  $d\hP\dfnn L_T
d\hQ^0$.

Now for any $\psi\in \hC_b(\hR)$, letting $n\to\infty$ we have
\bea
\label{lawconverge}
\lan \psi, \sT(\m^n)_t\ran&=& \hE^{\hP^n}\big[\psi(\hE^{\hP^n}[X^n_t|\cF^Y_t])\big]= \hE^{\hQ^0}\big[L^n_T
\psi(\hE^{\hP^n}[X^n_t|\cF^Y_t])\big]\nonumber\\
&\longrightarrow& \hE^{\hQ^0}\big[L_T\psi(\hE^{\hP}[X_t|\cF^Y_t])\big]=
\hE^{\hP}\big[\psi(\hE^{\hP}[X_t|\cF^Y_t])\big]\\
&=& \lan \psi, \hP\circ [\hE^{\hP}[X_t|\cF^Y_t]]^{-1}\ran, \qq \mbox{\rm as~} n\to\infty. \nonumber
\eea
This implies that $\n_t=\hP\circ [\hE^{\hP}[X_t|\cF^Y_t]]^{-1}=\sT(\m)_t$, for all $t\in[0,T]$. With the same argument one shows that,
for any $0\le t_1<t_2<\cds <t_k<\infty$,
$$\sT(\m^n)_{t_1, \cds, t_k}\dfnn \hP\circ\big(\hE^{\hP}[X^n_{t_1}|\cF^Y_{t_1}], \cds, \hE^{\hP}[X^n_{t_k}|\cF^Y_{t_k}])^{-1} \stackrel{d}\longrightarrow \n_{t_1, \cds, t_k}, \q \mbox{as $n\to \infty$. }
$$
That is, the finite dimensional distributions of $\sT(\m^n)$ converge to those of $\n$, and as $\{\sT(\m^n)\}_{n\ge
1}$ is tight by part (2), we  conclude that $\sT(\m^n)\limw \n$ in $\sP(\hC_T)$. This, together with
(\ref{estnu}), further shows that $W_2(\sT(\m^n), \sT(\m))\to 0$, as $n \to\infty$, proving the continuity of $\sT$, whence (3).
The proof is now complete.
\qed

 As a consequence of Theorem \ref{compact}, we have the following existence result for SDE (\ref{SDE}).
\begin{prop}
\label{existence} Let Assumption \ref{Assum2} hold. Then SDE (\ref{SDE}) has at least one solution in the sense of Definition \ref{sol}.
\end{prop}

{\it Proof.}
The proof follows from Theorem \ref{compact} and a generalization of the Schauder Fixed Point Theorem by Cauty (see \cite{Cauty}, or a recent generalization
%of this theorem made recently by Cauty himself the reader is referred to
\cite{Cauty2}). To do this we must check:
%In order to be able to apply Cauty's result we need that
(i) $\sE$ is a convex subset of a Hausdorff topological linear space, (ii)  $\sT$ is continuous and $\sT(\sE)\subseteq \sE$; and (iii) $\sT(\sE)\subset K$, for some compact $K$
in $\sP_2(\hC_T)$.

To imbed ${\sE}$ into a Hausdorff topological linear space, we borrow the argument of Li-Min  \cite{LiMin}.
Let $\sM_1(\hC_T)$ be the space of all bounded signed Borel measures $\n(\cd)$ on $\hC_T$ such that $|\int_{\hC_T}\|\f\|_{\hC_T}\n(d\f)|<+\infty$, endowed with the norm:
$$
\|\n\|_1:=\sup\Big\{\Big|\int_{\hC_T}hd\n\Big|\ :\ h\in\mbox{Lip}_1(\hC_T),~|h(0)|\le 1 \Big\}. \footnote{$Lip_1(\hC_T)$ denotes the set of all real-valued Lipschitz functions over $\hC_T$ with Lipschitz constant 1.}$$
Clearly $(\sM_1(\hC_T),\|\cdot\|_1)$ is a normed (hence  Hausdorff topological) linear space. Since $\sP_2(\hC_T)\subset\sP_1(\hC_T)\subset\sM_1(\hC_T)$, and by the Kantorovich-Rubinstein formula,
\beaa
 W_1(\n^1,\n^2)=\sup\Big\{\Big|\int_{\hC_T}hd(\n^1-\n^2)\Big|\ :\ h\in\mbox{Lip}_1(\hC_T),~|h(0)|\le 1 \Big\}
 %\\& =&\sup\Big\{\Big|\int_{\hC_T}(h-h(0))d(\n^1-\n^2)\Big|\ :\ h\in\mbox{Lip}_1(\hC_T),~|h(0)|\le 1 \Big\}\\
=\|\n^1-\n^2\|_1,
\eeaa
for all $\n^1,\n^2\in \cP_1(\hC_T)$, the
 topology generated by the norm $\|\cdot\|_1$ on $\sP_2(\hC_T)$ coincides with the one generated by the 1-Wasserstein metric on $\sP_2(\hC_T)$. Thus, $\sE\subset\sP_2(\hC_T)$ is a convex subset of $\sM_1(\hC_T)$, proving (i). Further, note
 that
 %a Hausdorff topological linear space, it is obviously convex, and from the preceding theorem we know that
 $\sT:\sE\rightarrow \sP_2(\hC_T)$ is continuous under the 1-Wasserstein metric, hence also under the $\|\cdot\|_1$-norm, verifying (ii). Finally,  since $\sT(\sE)\subset\sE$, and $\sE$ is compact under the 2-Wasserstein metric, hence also under the $\|\cdot\|_1$-norm, proving (iii). We can now apply Cauty's theorem to conclude the existence of a fixed point $\n\in\sE\subset\sP_2(\hC_T)$ such that $\sT(\n) =\n.$

We note that the existence of the fixed point $\m$ amounts to saying that  SDE (\ref{XBLQ0}) has
a solution on the probability space $(\O, \cF, \hQ^0)$, with $\m=\m^{X|Y}=\hP\circ [U]^{-1}$, and $U_t=\hE^{\hP}[X_t|\cF^Y_t]$, $t\ge 0$, where $d\hP=L_Td\hQ^0$ by construction. But this in turn defines a solution of (\ref{SDE}) on the probability space $(\O, \cF, \hP)$, thanks to the Girsanov transformation. However, since under $\hP$,
$(B^1, B^2)$ constructed in (\ref{XBLQ0}) is a Brownian motion, $(\O, \cF, \hP, X, Y, B^1, B^2)$ defines a
%we can identify it  with the original measure $\hP^0$. That is, we proved the existence of the
(weak) solution of SDE (\ref{SDE}).
\qed

\section{Uniqueness}
\setcounter{equation}{0}

In this section we investigate the uniqueness of the solution to SDE (\ref{SDE}).
%We first establish some fundamental estimates for the solutions of the state processes and their ``variational" processes.
We note that the general uniqueness for the weak solution for this problem is quite difficult, we will content ourselves with a version that is relatively more amendable.

To begin with, and let $\hQ^0$ be  the reference probability measure
under which $(B^1, Y)$ is a Brownian motion. For each $u\in L^{\infty-}_{\hF^Y}(\hQ^0,[0,T])$, consider the SDE on $(\O, \cF, \hQ^0)$:
\bea
\label{XBLQ1}
\left\{\ba{lll}
dX^u_t=\si(t, X^u_{\cd\wedge t}, \m^{X^u|Y}_t, u_t)dB^1_t, \q &X^u_0=x;\ms\\
dB^2_t=dY_t-h(t, X^u_{ t})dt, & B^2_0=0; \ms\\
dL^u_t= h(t, X^u_{ t})L^u_t dY_t, & L_0=1,
\ea\right.
\eea
where $\mu_t^{X^{u}|Y}:=\hP^{u}\circ[\hE^{\hP^{u}}[X^u_t|\cF_t^Y]]^{-1}$, and $d\hP^{u}:=L_T^{u}d\hQ^0$.
We shall argue that,  under Assumption \ref{Assum1}, the solution of the SDE (\ref{XBLQ1})  is pathwisely unique.
\begin{rem}
\label{remarkuniq} {\rm  It should be clear that if  $u\in  L^{\infty-}_{\hF^Y}(\hQ^0,[0,T])$, and $(X^u, B^2, L^u)$ is a solution to (\ref{XBLQ1}) under $\hQ^0$, then $u\in  L^{\infty-}_{\hF^Y}(\hP^u,[0,T])$ (since
%$\hP^u\sim \hQ^0$
%if $u\in L^{\infty-}_{\hF^Y}(\hP, [0,T])$ for some $\hP\in\sP(\O)$, then $u\in L^{\infty-}_{\hF^Y}(\tilde\hP, [0,T])$ for any $\tilde \hP\sim \hP$ with
$\frac{d\hP^u}{d\hQ^0}\in L^p(\O)$ for all $p>1$, thanks to Assumption \ref{Assum1}),
% (here ``$\sim$" stands for ``equivalent"). Next,  we),
 and the process $(X^u, Y, B^1, B^2)$ is a solution to (\ref{controlsys}) and (\ref{observation}) on the probability space $(\O, \cF, \hP^u, \hF)$ in the sense of Definition \ref{sol}, where $\hF:=\hF^{B^1, Y}$.
%Consequently, $(\hP^u, u)\in \sU_{ad}$ by Definition \ref{admissible}.
 %We note that the conditional law $\m^{X^u|Y}$ is completely determined by the solution
 Conversely,  if $(\O,\cF,\hP^u,\hF,B^1,B^2,X,Y)$ is a weak solution of (\ref{controlsys})-(\ref{observation}), then following the argument of \S2.2, we see that  $d\hQ^0=[L^{u}_T]^{-1}d\hP^u$ defines a reference measure, where $L^u$ is defined by (\ref{barL}) or (\ref{barLexp}), and $(X, B^2, [L^u]^{-1})$ will be a solution of (\ref{XBLQ1}) with respect to the $\hQ^0$-Brownian motion $(B^1,Y)$.
%In a control problem context, such a uniqueness result amounts to saying that, given a observation process $Y$ such that  $(B^1, Y)$ is a Brownian motion under the reference measure $\hQ^0$, then for any $(\hQ^0, u)\in \sU_{ad}$, the state-observation dynamics  $(X^u, Y, L^u)$ is pathwisely unique under $\hQ^0$, as the solution to (\ref{XBLQ1}), thanks to Assumption \ref{Assum1}.
In what follows we shall call the solution to (\ref{XBLQ1}) the $\hQ^0$-dynamics of the  system (\ref{controlsys}) and (\ref{observation}).}
\qed
\end{rem}

%2) If the control process . Then we see on the one hand that, if $X^{u}$ is a strong solution on $(\Omega,\cF,\hQ^0)$ with respect to the $\hQ^0$-Brownian motion $(B^1,Y)$, then $(\Omega,\cF,\hF,\hP^{u},B^1,B^{2},X^{u},Y)$ with $\hF$ the natural filtration generated by $(B^1,Y)$ and $dP^{u}=L^{u}_Td\hQ^0$ is a weak solution of (\ref{controlsys})-(\ref{observation}).
Bearing Remark \ref{remarkuniq} in mind, let us first try to establish a result in the spirit of the Yamada-Watanabe Theorem: {\it the pathwise uniqueness  of $(\ref{XBLQ1})$  implies the uniqueness in law for the original SDEs (\ref{controlsys}) and (\ref{observation})}. To do this, we begin by noting that, given the ``regular" nature of the canonical space $\O$, a process $u\in L^{\infty-}_{\hF^Y}(\hP^u,[0,T])$ amounts to saying that (cf. e.g., \cite{SV,yong-zhou}) there exists a
progressively measurable functional ${\bf u}:[0,T]\times \hC_T\mapsto U$ such that $u_t(\o)={\bf u}(t, Y_{\cd\wedge t}(\o))$, $dtd\hP^u$-a.s., such that $u$ has all the finite moments under $\hP^u$ (hence also true under $\hQ^0\sim
\hP^u$!).
We have the following Proposition.
\begin{prop}
\label{YamaWata}
 Assume that Assumption \ref{Assum1} is in force, and that the pathwise uniqueness holds for SDE (\ref{XBLQ1}).  Let ${\bf u}:[0,T]\times \O\mapsto U$ be a given progressively measurable functional, and  $(\Omega,\cF,\hP^{i},\hF,B^{1,i},B^{2,i},X^{i},Y^{i})$, $i=1,2$, be two  (weak) solutions of (\ref{controlsys})- (\ref{observation}) corresponding to the controls  $u^{i}={\bf u}(\cd, Y^{i})$, $i=1,2$, respectively. Then, it holds that
  %$Given a measurable functional $u(\cdot)$, for  let  be an and
$$\hP^{1}\circ[(B^{1,1},B^{2,1},X^{1},Y^{1})]^{-1}=\hP^{2}\circ[(B^{1,2},B^{2,2},X^{2},Y^{2})]^{-1}.$$
\end{prop}

{\it Proof.} Following the argument of \S2.2, we define $d\hQ^{0,i}=[L^{i}_T]^{-1}d\hP^{i}$, where $L^{i}=[\bar{L}^{i}]^{-1}$ and $\bar{L}^{i}$ is the unique solution of the SDE (\ref{barL}) with respect to $(X^i, B^{1,i}, Y^i)$, $i=1,2$.
%$d\bar{L}^{i}_t=-h(t,X_t^{i})\bar{L}^{i}_t dB_t^{2,i},\, \bar{L}^{i}_0=1.$
Then, as the $\hQ^{0,i}$-dynamics, $(X^i, B^{2,i},L^{i})$ satisfies (\ref{XBLQ1}), $i=1,2$,  $\hQ^{0,i}$-a.s.
%is an $(\hF^{i},$-Brownian motion and $X^{i}$ a solution of (\ref{XBLQ1}) governed by $(B^{1,i},Y^{i})$. In particular we have, for ,
%
%\centerline{$dX^{i}_t=\si(t,X^{i}_{\cdot\wedge t},\mu^{X^{i}|Y^{i}}_t,u_t(Y^{i}_{\cdot\wedge t}))dB^{1,i}_t,\, t\in[0,T],\ X^{i}_0=x,$}
%
%
%\noindent and we observe that $L^{i}$ is the unique solution of the equation
%
%\centerline{$dL^{i}_t=h(t,X^{i}_t)L_t^{i}dY^{i}_t,\, t\in[0,T],\, L^{i}_0=1.$}
%
In particular, we recall  (\ref{KS}) that
%
%\centerline{$\displaystyle\si (t,\f_{\cdot\wedge t},\mu_t^{X^{i}|Y^{i}},v)=\hE^{\hP^{i}}\left[\si(t,\f_{\cdot\wedge t},\hE^{\hP^{i}}[X_t^{i}|\cF^{Y^{i}}_t],v)\right]=\hE^{\hQ^{0,i}}\left[L_t^{i}\si(t,\f_{\cdot\wedge t},\hE^{\hP^{i}}[X_t^{i}|\cF^{Y^{i}}_t],v)\right],$}
%
%\noindent with
%
$$U_t^{X^{i}|Y^{i}}=\hE^{\hP^{i}}[X_t^{i}|\cF^{Y^{i}}_t]= \frac{\hE^{\hQ^{0,i}}[L_t^{i}X_t^{i}|\cF^{Y^{i}}_t]}{\hE^{\hQ^{0,i}}\left[L_t^{i}|\cF^{Y^{i}}_t\right]},\qq \hQ^{0,i}\mbox{-a.s.},\, t\in[0,T].
$$

Thus, there exist two progressively measurable functionals $\Phi^{i}:[0,T]\times \O\mapsto \hR$ such that $U_t^{X^{i}|Y^{i}}=\Phi^{i}(t,Y^{i}_{\cdot\wedge t})$, $dtd\hQ^{0,i}$-a.s., $i=1,2$. We now consider an intermediate SDE on $(\O, \cF, \hQ^{0,2})$:
% We compare now the above SDE for $X^{1}$,
%
%\centerline{$dX^{1}_t=\si(t,X^{1}_{\cdot\wedge t},\Phi^{1}(t,Y^{1}_{\cdot\wedge t}),u_t(Y^{1}_{\cdot\wedge t}))dB^{1,1}_t,\, t\in[0,T],\, Q^{0,1}\mbox{-a.s.},\ X^{1}_0=x,$}
%
%\noindent with the SDE
%
\bea
\label{hatX2}
\left\{\ba{lll}
d\widehat{X}^{2}_t=\si(t,\widehat{X}^{2}_{\cdot\wedge t},\Phi^{1}(t,Y^{2}_{\cdot\wedge t}),{\bf u}(t, Y^{2}_{\cdot\wedge t}))dB^{1,2}_t,\q &\widehat{X}^{2}_0=x;  \ms \\
d\widehat{L}^{2}_t=h(t,\widehat{X}^{2}_t)\widehat{L}_t^{2}dY^{2}_t, &\widehat{L}^{2}_0=1,
\ea\right.
\qq t\in[0,T].
%, ~\hQ^{0,2}\mbox{-a.s.}
\eea
Clearly, comparing to (\ref{XBLQ1}) for $\hQ^{0,1}$-dynamics $(X^1, B^{2,1}, L^1)$, this SDE has the same coefficient $\widehat{\si}(t, \o, \varphi_{\cdot\wedge t}):=\si(t,\varphi_{\cdot\wedge t},\Phi^1(t,\o^2_{\cdot\wedge t}),{\bf u}(t,\o^2_{\cdot\wedge t}))$, and
$h(t,x)\ell$,  which is jointly measurable, uniformly Lipschitz in $\varphi$ with linear growth (in $\ell$),
 uniformly in  $(t,\o, \varphi, \ell)$, thanks to Assumption \ref{Assum1}, except that it is
driven by the $\hQ^{0,2}$-Brownian motion $(B^{1,2}, Y^2)$. Thus, by the classical SDE theory (cf. e.g., \cite{IW}) we know that
%the SDE (\ref{hatX2}) has a strong existence and is pathwisely unique. That is,
% Also the second SDE has a unique solution $\widehat{X}^2$, and
 there exists a (unique) measurable functional $\Psi:\hC_T\times\hC_T\rightarrow\hC_T\times\hC_T$ such that $(X^1,L^1)=\Psi(B^{1,1},Y^1),\, \hQ^{0,1}$-a.s., and
$(\widehat{X}^2,\widehat{L}^2)=\Psi(B^{1,2},Y^2),\, \hQ^{0,2}$-a.s. % $\Psi:\hC_T\times\hC_T\rightarrow\hC_T$ such that $X^1=\Psi(B^{1,1},Y^1)$, $\hQ^{0,1}$-a.s., and $\widehat{X}^2=\Psi(B^{1,2},Y^2)$, $\hQ^{0,2}$-a.s.
Since
% define the probability measure , we conclude, again with the argument of the pathwise uniqueness for classical SDEs, that there is a measurable functional Then it follows from \\
$\hQ^{0,1}\circ (B^{1,1},Y^1)^{-1}=\hQ^{0,2}\circ (B^{1,2},Y^2)^{-1}=\hQ^0$, the Wiener measure on $(\O, \cF)$,  we deduce that
\bea
\label{Q01=Q02}
\hQ^{0,1}\circ (B^{1,1},Y^1,X^1,L^1)^{-1}=\hQ^{0,2}\circ (B^{1,2},Y^2,\widehat{X}^2,\widehat{L}^{2})^{-1}.
\eea

We now claim that $(\h{X}^2, B^{2,2}, \h{L}^2)$ coincides with the $\hQ^{0,2}$-dynamics of (\ref{controlsys})-(\ref{observation}).
Indeed, it suffices to argue that in SDE (\ref{hatX2}),
\bea
\label{hatmu}
\Phi^1(t,Y^2_{\cdot\wedge t})=\hE^{\widehat{\hP}^{2}}[\widehat{X}_t^2|\hF_t^{Y^2}]=U_t^{\widehat{X}^2|Y^2},\qq \hQ^{0,2}\mbox{-a.s.},
\eea
where $d\widehat{\hP}^2:=\widehat{L}^2d\hQ^{0,2}$. To see this, we note that,
%
% immediate consequence of this equality of the laws is that,
for all $t\in[0,T]$ and any bounded Borel measurable function $f:\hC_T\rightarrow \hR$, it follows from (\ref{Q01=Q02}) and the
definition of $U_t^{X|Y}$ that
\beaa
 &&\hE^{\widehat{\hP}^{2}}[f(Y^2_{\cdot\wedge t})\Phi^1(t,Y^2_{\cdot\wedge t})]=\hE^{\hQ^{0,2}}[\widehat{L}_t^2 f(Y^2_{\cdot\wedge t})\Phi^1(t,Y^2_{\cdot\wedge t})]=\hE^{\hQ^{0,1}}[L_t^1 f(Y^1_{\cdot\wedge t})\Phi^1(t,Y^1_{\cdot\wedge t})]\\
&=&\hE^{\hP^{1}}[f(Y^1_{\cdot\wedge t})U_t^{X^1|Y^1}] =\hE^{\hP^{1}}[f(Y^1_{\cdot\wedge t})X_t^1]=\hE^{\hQ^{0,1}}[L_t^1f(Y^1_{\cdot\wedge t})X^1_t]=\hE^{\hQ^{0,2}}[\widehat{L}_t^2f(Y^2_{\cdot\wedge t})\widehat{X}^2_t]\\
&=&\hE^{\h{\hP}^{2}}[f(Y^2_{\cdot\wedge t})\widehat{X}_t^2]=\hE^{\h{\hP}^{2}}[f(Y^2_{\cdot\wedge t})U_t^{\widehat{X}^2|Y^2}],
\eeaa
proving (\ref{hatmu}), whence the claim.

Now, by pathwise uniqueness of SDE (\ref{XBLQ1}), we conclude that $(X^2,L^2)=(\h{X}^2,\h{L}^2)$, $\hQ^{0,2}$-a.s.
Thus (\ref{Q01=Q02}) implies that $ \hQ^{0,1}\circ [(B^{1,1},Y^1,X^1,L^1)]^{-1}
%=Q^{0,2}\circ [(B^{1,2},Y^2,\widehat{X}^2,\widehat{L}^{2})]^{-1}
=\hQ^{0,2}\circ [(B^{1,2},Y^2,X^2,L^{2})]^{-1}$, and consequently,
$\hQ^{0,1}\circ [(B^{1,1},B^{2,1},X^1,Y^1)]^{-1}=\hQ^{0,2}\circ [(B^{1,2},B^{2,2},X^2,Y^2)]^{-1}$. This proves
 the uniqueness in law for the system (\ref{controlsys})-(\ref{observation}).
 \qed

We now turn our attention to the main result of this section: the pathwise uniqueness of (\ref{XBLQ1}). We shall  establish
some fundamental estimates which will
be useful in our future discussions. Since all controlled dynamics are constructed via the reference probability
space $(\O, \cF, \hQ^0)$, we shall consider only their $\hQ^0$-dynamics, namely the solution to (\ref{XBLQ1}).
%Consider two controls $u,v\in \sU_{ad}$ under the reference measure $\hQ^0$.
 Recall the space $L^p(\hQ^0; L^2([0,T]))$, $p>1$, and the norm $\|\cd\|_{p,2,\hQ^0}$
defined by (\ref{LpqPnorm}).
We have the following important result.
\begin{prop}
\label{est1}
Assume that Assumption \ref{Assum1} is in force. Let  $u, v\in \sU_{ad}$ be given. Then, for any $p>2$, there exists a constant $C_p>0$, such that the following estimates hold:

\bea
\label{7.25}
&&{\rm (i)}\ \hE^{\hQ^0}\Big[\sup\limits_{0\leq s\leq T}(|X_s^{{u}}-X_s^{{v}}|^2+
|L_s^{{u}}-L_s^{{v}}|^2+|X_s^{{u}}L_s^{{u}}
-X_s^{{v}} L_s^{{v}}|^2)\Big]\leq C\|u-v\|^{2}_{2,2,\hQ^0};\ \ \ \ \ \\
\label{Xpest}
&&{\rm (ii)}\ \hE^{\hQ^0}\Big[\sup\limits_{0\leq s\leq T}|X_s^{{u}}-X_s^{{v}}|^p\Big]\leq C_p\|u-v\|^{p}_{p,2,\hQ^0}.
\eea
\end{prop}
{\it Proof.} We split the proof into several steps. Throughout this proof we let $C>0$ be a generic constant, depending only
on the bounds and Lipschitz constants of the coefficients and the time duration $T>0$, and it is allowed to vary from line to line.

\ms
{\it Step 1} ({\it Estimate for $X$}). First let us denote, for any $u\in \sU_{ad}$,
\bea
\label{sigmamuu}
 && \si^u(t, \f_{\cd\wedge t}, \m^u_t) \dfnn  \int_{\hR}\si(t, \f_{\cd\wedge t}, y, u_t)\m^u_t(dy), ~~(t, \f)\in [0,T]\times \hC_T,
 \eea
and
$ \m^u_t \dfnn \m^{X^u|Y}\circ P^{-1}_t=\hP^u\circ (\hE^{\hP^u}[X^u_t|\cF^Y_t])^{-1}$, $t\ge 0$. Then, we have
\bea
\label{siuv}
&&|\si^u(t, X^u_{\cd\wedge t}, \m^u_t)-\si^v(t, X^v_{\cd\wedge t}, \m^v_t)|\nonumber\\
&=& \Big|\int_{\hR}\si(t, X^u_{\cd\wedge t}, y, u_t)\m^u_t(dy)-\int_{\hR}\si(t, X^v_{\cd\wedge t}, y, v_t)\m^v_t(dy)\Big|\\
&\le & C\Big\{ |u_t-v_t|+\sup_{0\le s\le t}|X^u_s-X^v_s|+\Big|\int_{\hR}\si(t, X^v_{\cd\wedge t}, y, v_t)[\m^u_t(dy)-
%\int_{\hR}\si(t, X^v_{\cd\wedge t}, y, v_t)
\m^v_t(dy)]\Big|\Big\}.\nonumber
\eea
Next, let us denote $S^u_t=\hE^{\hQ^0}[L^u_tX^u_t|\cF^Y_t]$ and $S^{u,0}_t=\hE^{\hQ^0}[L^u_t|\cF^Y_t]$, and define
$S^{v}_t$, $S^{v,0}_t$ in a similar way. By (\ref{Bayes}) and the fact that $d\hP^u=L^u_Td\hQ^0$,  we see that
\bea
\label{I1I2}
&&\Big|\int_{\hR}\si(t, X^v_{\cd\wedge t}, y, v_t)[\m^u(dy)-\m^v(dy)]\Big|\nonumber\\
&=& \Big|\hE^u[\si(t, \f_{\cd\wedge t},\hE^u[X^u_t|\cF^Y_t], u)]-\hE^v[\si(t, \f_{\cd\wedge t},\hE^v[X^v_t|\cF^Y_t], u)]\big|_{\f=X^v, u=v_t}\Big|\\
&=& \Big|\hE^{\hQ^0}\Big\{L^u_t \si\big(t, \f_{\cd\wedge t},\frac{\hE^{\hQ^0}[L^u_tX^u_t|\cF^Y_t]}{\hE^{\hQ^0}[L^u_t|\cF^Y_t]}, u\big)
-L^v_t \si\big(t, \f_{\cd\wedge t},\frac{\hE^{\hQ^0}[L^v_tX^v_t|\cF^Y_t]}{\hE^{\hQ^0}[L^v_t|\cF^Y_t]}, u\big)\Big\}\big|_{\f=X^v, u=v_t}\Big|\nonumber\\
&\le & I_1+I_2, \nonumber
\eea
where (noting the definition of $S^u$, $S^{u,0}$ and the fact that they are both $\hF^Y$-adapted)
\beaa
\label{I1}
I_1&=&\Big|\hE^{\hQ^0}\Big\{L^u_t \si\big(t, \f_{\cd\wedge t},\frac{S^u_t}{S^{u,0}_t}, u\big)
-L^v_t \si\big(t, \f_{\cd\wedge t},\frac{S^u_t}{S^{v,0}_t}, u\big)\Big\}\big|_{\f=X^v, u=v_t}\Big|\\
&=&\Big|\hE^{\hQ^0}\Big\{S^{u,0}_t \si\big(t, \f_{\cd\wedge t},\frac{S^u_t}{S^{u,0}_t}, u\big)
-S^{v,0}_t \si\big(t, \f_{\cd\wedge t},\frac{S^u_t}{S^{v,0}_t}, u\big)\Big\}\big|_{\f=X^v, u=v_t}\Big|;
\eeaa
and
\beaa
I_2&=&\Big|\hE^{\hQ^0}\Big\{L^v_t \big[\si\big(t, \f_{\cd\wedge t},\frac{S^u_t}{S^{v,0}_t}, u\big)
-\si\big(t, \f_{\cd\wedge t},\frac{S^v_t}{S^{v,0}_t}, u\big)\big]\Big\}\big|_{\f=X^v, u=v_t}\Big|\\
&=&\Big|\hE^{\hQ^0}\Big\{S^{v,0}_t \big[\si\big(t, \f_{\cd\wedge t},\frac{S^u_t}{S^{v,0}_t}, u\big)
-\si\big(t, \f_{\cd\wedge t},\frac{S^v_t}{S^{v,0}_t}, u\big)\big]\Big\}\big|_{\f=X^v, u=v_t}\Big|.
\eeaa
Clearly, we have
\bea
\label{I2est}
I_2\le C\hE^{\hQ^0}\Big\{S^{v,0}_t\frac{|S^{u}_t-S^{v}_t|}{S^{v,0}_t}\Big\}\le C\hE^{\hQ^0}\left[|L^u_tX^u_t-L^v_tX^v_t|\right].
\eea
%\beaa
%\Big|\hE^0\Big\{\hE^0[L^u_t|\cF^Y_t] \si(t, \f_{\cd\wedge t},\frac{\hE^0[L^u_tX^u_t|\cF^Y_t]}{\hE^0[L^u_t|\cF^Y_t]}, u)\Big]\Big|_{\f=X, u=u_t}\nonumber\eeaa
To estimate $I_1$, we write $\hat \si(t,\o,\f_{\cd\wedge t}, y,z)=y\si\big(t,\f_{\cd\wedge t}, \frac{S^u_t(\o)}{y}, z\big)$. Since
\bea
\label{pahatsi}
 \pa_y\hat \si(t,\o, \f_{\cd\wedge t}, y,z)=\si\Big(t,\f_{\cd\wedge t}, \frac{S^u_t(\o)}{y}, z\Big)-\frac{S^u_t(\o)}{y}\pa_y\si
\Big(t,\f_{\cd\wedge t}, \frac{S^u_t(\o)}{y}, z\Big),
\eea
we see that $y\mapsto \pa_y\hat\si(t,\f_{\cd\wedge t}, y, z) $ is uniformly bounded thanks to Assumption
\ref{Assum1}-(iv). Thus we have
\bea
\label{I1est}
I_1
\le C\|\pa_y\hat\si\|_{\infty}\hE^{\hQ^0}|S^{u,0}_t-S^{v,0}_t|\le C\hE^{\hQ^0}|L^u_t-L^v_t|.
\eea
%we note that for any $x_1,\ x_2,\ y,$\ we define $\tilde{\si}(x):=\si(x,\varphi, u),$ then we have:
%\begin{equation}\label{7.26}|x_1\tilde{\si}(\frac{y}{x_1})-x_2\tilde{\si}(\frac{y}{x_2})|=|\tilde{\si}(\frac{y}{\theta_y})
%-\frac{y}{\theta_y}\tilde{\si}_x(\frac{y}{\theta_y})|\cdot|x_1-x_2|\leq C|x_1-x_2|,\end{equation}
%for some $\theta_y$\ between $x_1$\  and $x_2$.
Now note that (\ref{XBLQ1}) implies that
$
X_t^{{u}}-X_t^{{v}}=\int_0^t[\si^u(s, X^u_{\cd\wedge s}, \m^{u}_s)-\si^v(s, X^v_{\cd\wedge s}, \m^{v}_s)]dB^1_s$.
%\mathcal{F}_s^Y], \varphi_{\cd\wedge  s}, \tilde{u})]\big|_{\varphi=X^{\tilde{u}};\atop \tilde{u}=\tilde{u}_s}\right.\\
%& &-E^{\tilde{v}}[\si(E^{\tilde{v}}\left.[X_s^{\tilde{v}}|\mathcal{F}_s^Y], \varphi_{\cd\wedge  s}, \tilde{v})]\big|_{\varphi=X^{\tilde{v}};\atop \tilde{v}=\tilde{v}_s}\right)dB_s^1,\ \ t\geq 0.
%\eeaa
%
Combining (\ref{siuv})--(\ref{I1est}), we see that
%, we get the following inequality with the help of (\ref{7.26}):\\
\bea
\label{7.27}
\hE^{\hQ^0} \Big[\sup_{0\leq s\leq t}|X_s^{u}-X_s^{v}|^p\Big]&\le&C\hE^{\hQ^0}\Big\{\Big[\int_0^t[\sup_{r\in[0,s]}|X_r^{u}-X_r^{v}|^2+|u_s-v_s|^2\\
&&\qq+ (\hE^{\hQ^0}|L_s^{u}-L_s^{v}|)^2+(\hE^{\hQ^0}|L_s^{u} X_s^{u}-L_s^{v}X_s^v|)^2]ds\Big]^{p/2}\Big\}.\nonumber
\eea
Applying the Gronwall inequality we obtain that
\bea
\label{7.28}
\hE^{\hQ^0} \Big[\sup_{0\leq s\leq t}|X_s^{u}-X_s^{v}|^p\Big]&\leq& C\hE^{\hQ^0}\Big\{\Big[\int_0^t\big[|u_s-v_s|^2+\hE^{\hQ^0}[|L_s^{u}-L_s^{v}|^2]
\nonumber\\
&& +\hE^{\hQ^0}[|L_s^{u}X_s^{u}-L_s^{v} X_s^v|^2]\big]ds\Big]^{p/2}\Big\}.
\eea

\ms
{\it Step 2} ({\it Estimate for $L$}). We first note that, for $t\in[0,T]$,
\bea
\label{Lh}
 &&|L_t^{u}h(t, X_{t}^{u})-L_t^{v} h(t, X_{ t}^{v})| =\Big|L_t^{u}h\Big(t, \frac{L^u_tX_{t}^u}{L^u_t}\Big)-L_t^{v} h\Big(t, \frac{L^v_tX_{t}^v}{L^v_t}\Big)\Big|\nonumber\\
 &\le& \Big|L_t^{u}h\Big(t, \frac{L^u_tX_{t}^u}{L^u_t}\Big)-L_t^{u} h\Big(t, \frac{L^v_t X_{t}^v}{L^u_t}\Big)\Big|+
 \Big|L_t^{u}h\Big(t, \frac{L^v_tX_{t}^v}{L^u_t}\Big)-L_t^{v} h\Big(t, \frac{L^v_tX_{t}^v}{L^v_t}\Big)\Big|
 \\
 &\le &C |L_t^{u}X_{ t}^{u}-L_t^{v}X_{t}^{v}|+\Big|L_t^{u}h\Big(t, \frac{L^v_tX_{t}^v}{L^u_t}\Big)-L_t^{v} h\Big(t, \frac{L^v_tX_{t}^v}{L^v_t}\Big)\Big|.\nonumber
 \eea
To estimate the second term above we define, as before, $\hat h(t,\o, x)\dfnn xh\big(t, \frac{ L^v_t(\o)X^v_t(\o)}x\big)$. Then, similar to (\ref{pahatsi}), one shows that $x\mapsto \pa_x\hat h(t,\o, x)$ is uniformly bounded, thanks to Assumption
\ref{Assum1}-(v). Consequently, we have
\bea
\label{pahath}
\Big|L_t^{u}h\Big(t, \frac{L^v_tX_{t}^v}{L^u_t}\Big)-L_t^{v} h\Big(t, \frac{L^v_tX_{t}^v}{L^v_t}\Big)\Big|\le
 \|\pa_x\hat h\|_\infty |L_t^{u}-L_t^{v}|.
 \eea
%  \nonumber\\
% &&
% |h(t, X_{  t}^{u})-h(t, X_{ t}^{v})|\\
% && \qq \qq \qq \le C|L_t^{u}-L_t^{v}|+CL_t^{v}  |X_t^{u}-X_t^{v}|
%% &\leq C|L_t^{u}-L_t^{{v}}|+C|L_t^{u}X_{ t}^{u}-L_t^{v} X_{ t}^{v}|+C|X_{ t}^{u}| |L_t^{u}-L_t^{v}|\\
% \leq C(|L_t^{u}-L_t^{v}|+). \nonumber
%\eea
Now, combining (\ref{Lh}) and (\ref{pahath}) we obtain
\bea
\label{7.29}
|L_t^{u}h(t, X_{t}^{u})-L_t^{v} h(t, X_{ t}^{v})|
 %\leq C|L_t^{u}-L_t^{v}|+CL_t^{v}  |h(t, X_{  t}^{u})-h(t, X_{ t}^{v})|\\
% && \qq \qq \qq \le C|L_t^{u}-L_t^{v}|+CL_t^{v}  |X_t^{u}-X_t^{v}|
% &\leq C|L_t^{u}-L_t^{{v}}|+C|L_t^{u}X_{ t}^{u}-L_t^{v} X_{ t}^{v}|+C|X_{ t}^{u}| |L_t^{u}-L_t^{v}|\\
 \leq C(|L_t^{u}-L_t^{v}|+|L_t^{u}X_{ t}^{u}-L_t^{v}X_{t}^{v}|).
 %\nonumber
\eea
Therefore, noting that $L_t^{{u}}=1+\int_0^t h(s, X_{ s}^{{u}})L_s^{{u}}dY_s$, we deduce from (\ref{7.29}) and Gronwall's
inequality that
\begin{equation}
\label{7.31}
\hE^{\hQ^0}[\sup_{0\leq s\leq t}|L_s^{u}-L_s^{v}|^2]\leq
C\hE^{\hQ^0}[\int_0^t|L_s^{u} X_{ s}^{u}-L_s^{v}X_{s}^{v}|^2ds],\ \hQ^0\mbox{-a.s.},\ 0\leq t\leq T.
 \end{equation}

 \bs

{\it Step 3} ({\it Estimate for $L_tX_t$}). It is clear from (\ref{7.28}) and (\ref{7.31}) that it suffices to find the estimate of
$L^u_tX^u_t-L^v_tX^v_t$ in terms of $u-v$. To see this we note that
\bea
\label{LX}
L_t^{u} X_t^{u}=x+\int_0^t L_s^{u} X_s^{u}h(s,X_{s}^{u})dY_s+\int_0^t L_s^{u}\hE^{\hP^u}[\si (s, \varphi_{\cd\wedge  s},
\hE^{\hP^u}[X_s^{u}|\mathcal{F}_s^Y],  v)]\big|_{\varphi=X^{u}\atop v=u_s}dB_s^1.
\eea
Now define $\tilde h(t,x)\dfnn xh(t,x)$. Then it is easily seen that as $h$ satisfies Assumption \ref{Assum1}-(vi),
$\tilde h$ satisfies Assumption \ref{Assum1}-(v). Thus, similar to (\ref{7.29}) we have
\bea
\label{7.32}
|L_s^u X_s^u h(s, X_{s}^u)-L_s^v X_s^{v}h(s, X_{s}^{v})|&=&|L_s^u\tilde h(s,X_{s}^u)- L_{s}^{v}\tilde h(s,X_s^{v})|\nonumber\\
&\leq& C(|L_s^{u}-L_s^{v}|+|L_s^u X_{  s}^{u}-L_s^v X_{s}^{u}|).
\eea
On the other hand, for any $u\in \sU_{ad}$, recalling (\ref{sigmamuu}) for the notations $\si^u$ and $\m^u$,
we have,
\beaa
\label{7.33}
\Delta^{u,v}_t &\dfnn& \Big|L_s^u \hE^{\hP^u}[\si (s,  \varphi_{\cd\wedge  s}, \hE^{\hP^u}[X_s^{u}|\mathcal{F}_s^Y], z)]\big|_{\varphi=X^{u};\atop z=u_s}
-L_s^v\hE^{\hP^v} [\si (s, \varphi_{\cd\wedge  s}, \hE^{\hP^v}[X_s^{v}|\mathcal{F}_s^Y], z)]\big|_{\varphi=X^{v}\atop z=v_s}\Big| \nonumber\ms \\
&=& \big|L_t^u \si^u (t, X^{u}_{\cd\wedge  t}, \m^u_t)
-L_t^v \si^v (t, X^{v}_{\cd\wedge  t}, \m^v_t)\big|.
\eeaa
Then, following a similar argument as in Step 1 we have
\beaa
\label{Deltauv}
\Delta^{u,v}_t
 &\leq& CL_t^{v}(\hE^{\hQ^0}[|L_t^{u}-L_t^{v}|]+\hE^{\hQ^0}[|X_t^{u} L_t^{u}-X_t^{v}L_t^v|])  \\
 && +C(|L_t^{u}-L_t^{v}|+|L_t^u X_t^{u}-L_t^v X_t^{v}|)+CL_t^v |u_t-v_t|.\nonumber
\eeaa
Squaring both sides above and then taking the expectations we easily deduce that
\begin{equation}
\label{7.34}
\hE^{\hQ^0}[|\Delta^{u,v}_t|^2]\leq C(\hE^{\hQ^0}[|L_s^u-L_s^{v}|^2]+\hE^{\hQ^0}[|X_t^{u} L_t^{u}-X_t^v L_t^{v}|^2])+C\hE^{\hQ^0}[(L_t^{v})^2 |u_t-v_t|^2].
\end{equation}
Now, combining (\ref{LX})-- (\ref{7.34}), for $p>2$ we can find $C_p>0$ such that
\bea
\label{7.35}
&&\hE^{\hQ^0}\Big[\sup_{0\leq s\leq t}|L_s^uX_s^u  -L_s^vX_s^v |^2\Big] \nonumber\\
&\leq& C\hE^{\hQ^0}\Big[\int_0^t |L_s^u X_s^u h(s,X_{ s}^{u})-L_s^v X_s^v h(s,X_{ s}^v)|^2ds\Big]+C\hE^{\hQ^0}\int_0^t|\Delta^{u,v}_s|^2ds\\
&\leq& C_p\Big\{\hE^{\hQ^0}\Big[\Big(\int_0^t |u_s-v_s|^2ds\Big)^{p/2}\Big]\Big\}^{2/p}+C\hE^{\hQ^0}\int_0^t|L_s^{u}-L_s^{v}|^2ds \nonumber\\
&& +C\hE^{\hQ^0}\int_0^t|L_s^u X_s^{u}-L_s^v X_s^{v}|^2ds. \nonumber
\eea
Hence, applying Gronwall's inequality we obtain
\bea
\label{7.36}
\hE^{\hQ^0}\Big[\sup_{0\leq s\leq t} |L_s^uX_s^u  -L_s^vX_s^v |^2\Big]\leq C_p\|u-v\|^2_{p,2,\hQ^0}+
%(E^0[(\int_{0}^{t}|\tilde{u}_s-\tilde{v}_s|^2ds)^p])^{\frac{1}{p}}+
C\hE^{\hQ^0}\int_{0}^{t}|L_s^u-L_s^v|^2ds.
\eea

Combining  (\ref{7.36}) with (\ref{7.31}) and applying the Gronwall inequality again, we conclude that
\bea
\label{7.37}
\hE^{\hQ^0}\Big\{\sup_{0\leq s\leq t}|L_s^u-L_s^v|^2\Big\}\leq C_p\|u-v\|^2_{p,2,\hQ^0}.
%E^0[(\int_{0}^{t}|\tilde{u}_s-{\tilde{v}}_s|^2ds)^p])^\frac{1}{p}.
\eea

This, together with  (\ref{7.28}) and (\ref{7.36}), implies (\ref{7.25}). (\ref{Xpest}) then follows easily from (\ref{7.25}) and
(\ref{7.27}), proving the proposition.
\qed

A direct consequence of Proposition \ref{est1} is the following uniqueness result.
% for the solution to (\ref{XBLQ1}).
\begin{cor}
\label{uniqueness}
Assume that Assumption \ref{Assum1} holds. Then the solution to SDE (\ref{XBLQ1}) is pathwisely unique.
\end{cor}

{\it Proof.} Setting $u=v$ in Proposition \ref{est1} we obtain the result.
\qed

\section{A Stochastic Control Problem with Partial Observation}
\setcounter{equation}{0}

%\subsection{\large{Formulation of the Problem}}
We are now ready to study the stochastic control problem with partial observation. We first note that
in theory for each $(\hP^u, u)\in\sU_{ad}$ our state-observation dynamics $(X^u, Y^u)$ lives on probability space
$(\O, \cF, \hP^u)$, which varies with control $u$. We shall consider their $\hQ^0$-dynamics so that
our analysis can be carried out on a common probability space, thanks to Assumption \ref{Assump2}.
Therefore, in what follows, for each $(\hP^u, u)\in \sU_{ad}$ we consider only the $\hQ^0$-dynamics $(X^u, Y, L^u)$, which
satisfies the following SDE:
\bea
\label{7.18a}
\left\{
\begin{array}{lll}
\dis dX_t^u=\si^u(t, X^u_{\cd\wedge t},\m^u_t)dB_t^1,\q & X^u_0=x; \ms \\
\dis dB^{2,u}_t=dY_t-h(t, X^u_{t})dt, & B^{2, u}_0=0; \ms\\
\dis dL^u_t=h(t, X^u_t)L^u_tdY_t, & L^u_0=1,\q t\geq 0, 
\end{array}
\right.
\eea
where $(B^1, Y)$ is a $\hQ^0$-Brownian motion, $d\hP^u=L_T^ud\hQ^0$, and $\mu_t^{X^{u}|Y}=\hP^{u}\circ[\hE^{\hP^{u}}[X_t|\cF_t^Y]]^{-1}$. For simplicity, we denote $\hE^u[\cdot]\dfnn \hE^{\hP^u}[\cdot]$ and
$\hE^0[\cd]\dfnn \hE^{\hQ^0}[\cd]$.

\begin{rem}{\rm
A convenient and practical way to identify admissible control is to simply consider the space $L^{\infty-}_{\hF^Y}(\hQ^0;[0,T])$
(cf. Definition \ref{admissible}),
which is independently well-defined, thanks to Assumption \ref{Assump2}. It is easy to check that, under Assumption \ref{Assum1},
 $u\in L^{\infty-}_{\hF^Y}(\hQ^0;[0,T])$ if and only if $u\in L^{\infty-}_{\hF^Y}(\hP^u;[0,T])$.
 %, where $\hP^u$ is the probability measure constructed via (\ref{7.18a}).
 Therefore in what follows by $u\in\sU_{ad}$ we  mean that $u\in L^{\infty-}_{\hF^Y}(\hQ^0;[0,T])$.
\qed}
\end{rem}
%
%We {\color{red}We call such pair $(u, \hQ^0)$ the {\it natural admissible control}. Next, for a natural admissible control $(u, \hQ^0)$, consider  the following SDE:}
%}

We recall that  for $u\in\sU_{ad}$ and $\m\in\sP_2(\hC_T)$, the coefficient $\si^u$ in (\ref{7.18a}) is defined
by (\ref{sigmamuu}). Thus we can write the cost functional as
%We define the following cost functional:
\begin{equation}
\label{7.20}
J(u)\dfnn \hE^0\Big\{ \Phi(X^u_T, \m^u_T)+
%E^v[\Phi(E^v[X_T^v|\mathcal{F}_T^Y],\varphi)]\Big|_{\varphi=X_T^v}]\\
\int_0^T f^u(s, X_s^u, \m^u_s)ds\Big\}.
%|\mathcal{F}_s^Y],\varphi_{\cd\wedge s},v)]\Big|_{\varphi=X^v\atop v=v_s}ds].
%\end{array}
\end{equation}
An admissible control $u^*\in\sU_{ad}$ is said to be optimal if
\begin{equation}\label{7.21}
J(u^*)=\inf_{u\in\sU_{ad}}J(u).
\end{equation}
We remark that the cost functional $J(\cd)$ involves the law of the conditional expectation of the solution in a nonlinear way.
Therefore, such a control problem is intrinsically ``{\it time-inconsistent}" and, thus, the dynamic programming approach  in general does not apply.
For this reason, we shall consider only the necessary condition of the optimal solution, that is, Pontryagin's
 Maximum Principle.
%\qed}
%\end{rem}

 To this end, we let  $u^*\in\sU_{ad}$ be
an {\it optimal control}, and consider the convex variations of $u^*$:
\begin{equation}
\label{7.22}
u_t^{\theta, v}:= u^*_t+\theta(v_t-u^*_t),\q t\in[0,T], \q 0<\theta<1,\q v\in\sU_{ad}.
\end{equation}
Here, we assume that $u^*, v\in L^{\infty-}_{\hF^Y}(\hQ^0; [0,T])$. Since $U$ is convex, $u^{\th, v}_t\in U$, for all $t\in[0,T]$, $ v\in\sU_{ad}$, and $\th\in(0,1)$. We denote $(X^{\theta, v}, Y, L^{\th, v})$ to be the corresponding $\hQ^0$-dynamics that satisfies (\ref{7.18a}), with control $u^{\theta, v}$. Applying Proposition \ref{est1} ((\ref{7.25}) and (\ref{Xpest})) and noting that $Y$ is a Brownian motion under $\hQ^0$,  we get,  for $p>2$,
\bea
\label{limXth}
\lim_{\th\to0} \hE^0\Big[\sup_{0\leq t\leq T}|X_t^{\th, v}-X_t^{u^*}|^2\Big]&\le& C_p  \lim_{\th\to0}\| u^{\th, v}-u^*\|^2_{p,2,\hQ^0}
=0; \\
\label{limLth}
\lim_{\th\to0} \hE^0\Big[\sup_{0\leq t\leq T}|L_t^{\th, v}-L_t^{u^*}|^2\Big]&=&0.
\eea

In the rest of the section we shall derive, heuristically, the ``variational equations" which play a fundamental role
in the study of Maximum Principle. The complete proof will be given in the next section.
% of the pair $(X^{\th, v}, L^{\th, v})$ for any $v\in\sU_{ad}$ around the optimal  control $u^*$,
For notational simplicity we shall denote $u=u^*$, the optimal control, from now on, bearing in mind that all
discussions will be carried out for the $\hQ^0$-dynamics, therefore on the same probability space.

%Let $u,v\in\sU_{ad}$,
%and let $(X^u, L^u)$ and $(X^v, L^v)$ denote the corresponding solutions  We denote $\d \xi^{u,v}\dfnn \xi^u-\xi^v$,
%for $\xi=X, L$, respectively, .

Now for $u^1,u^2\in\sU_{ad}$, let  $(X^1, L^1)$ and $(X^2, L^2)$ denote the corresponding solutions of (\ref{7.18a}). We define
$\d X=\d X^{1,2}=\d X^{u^1,u^2}\dfnn X^{u^1}-X^{u^2}$ and $\d L=\d L^{1,2}=\d L^{u^1,u^2}\dfnn L^{u^1}-L^{u^2}$, and will
often drop the superscript ``$^{1,2}$" if the context is clear. Then $\d X$ and $\d L$ satisfy the equations:
\bea
\label{deltaX}
\left\{\ba{lll}
\dis \d X_t = \int_0^t [\si^{u^1}(s, X^1_{\cd\wedge s}, \m^1_s)-\si^{u^2}(s, X^2_{\cd\wedge s}, \m^2_s)]dB^1_s;\ms\\
\dis \d L_t=\int_0^t [L^1_sh(s, X^1_s)-L^2_sh(s, X^2_s)]dY_s.
\ea\right.
\eea
As before, let $U^i_t\dfnn \hE^{u^i}[X^i_t|\cF^Y_t]$ and $\m^i_t=\hP^{u^i}\circ [U^i_t]^{-1}$,  $t\ge0$,  $i=1,2$. We can easily check that
\bea
\label{deltasigma}
&& \si^{u^1}(t, X^1_{\cd\wedge t}, \m^1_t)-\si^{u^2}(t, X^2_{\cd\wedge t}, \m^2_t) \nonumber\\
&=& \hE^0\Big\{L^1_t \si(t, \f^1_{\cd\wedge t}, U^1_t, z^1)-L^2_t \si(t, \f^2_{\cd\wedge t}, U^2_t, z^2)\Big\}\Big|_{\f^1=X^1, \f^2=X^2;z^1=u^1_t, z^2=u^2_t}\nonumber\\
&=& \hE^0\Big\{\d L^{1,2}_t \si(t, \f^1_{\cd\wedge t},U^1_t, z^1)\\
&&+L^2_t \Big[\int_0^1D_\f\si(t, \f^2_{\cd\wedge t}+\l (\f^1_{\cd\wedge t}-\f^2_{\cd\wedge t}), U^1_t, z^1)(\f^1_{\cd\wedge t}-\f^2_{\cd\wedge t})d\l\nonumber\\
&&+\int_0^1\pa_y\si(t, \f^2_{\cd\wedge t}, U^2_t+\l(U^1_t-U^2_t), z^1)d\l \cd(U^1_t-U^2_t)\nonumber\\
&&+\int_0^1\pa_z\si(t, \f^2_{\cd\wedge t}, U^2_t, z^2+\l(z^1-z^2))d\l \cd(z^1-z^2)\Big]\Big\}\Big|_{\f^1=X^1, \f^2=X^2; z^1=u^1_t, z^2=u^2_t}. \nonumber
\eea
Now let $u^1=u^{\th, v}$ and $u^2=u^*=u$, and denote
$$ \d_\th X\dfnn \d_\th X^{u, v}=\frac{X^{\th, v}-X^u}{\th}, \q \d_\th L\dfnn \d_\th L^{u, v}=\frac{L^{\th, v}-L^u}{\th}, \q
\d_\th U\dfnn \d_\th U^{u,v}=\frac{U^{\th, v}-U^u}{\th}.
$$
Combining (\ref{deltaX}) and (\ref{deltasigma}) we have
\bea
\label{dthX}
\d_\th X_t
%&& \si^{u^1}(t, X^1_{\cd\wedge t}, u^1_t)-\si^{u^2}(t, X^2_{\cd\wedge t}, u^2_t) \nonumber\\
%&=& \hE^0\Big\{L^1_t \si(t, \f^1_{\cd\wedge t}, U^1_t, z^1)-L^2_t \si(t, \f^2_{\cd\wedge t}, U^2_t, z^2)\Big\}\Big|_{\f^1=X^1, \f^2=X^2;z^1=u^1_t, z^2=u^2_t}\nonumber\\
%&=&  \hE^0\Big\{(L^u_t -L^v_t) \si(t, \f^1_{\cd\wedge t}, U^u_t, z^1)
%%\nonumber\\&&
%+L^v_t [\si(t, \f^1_{\cd\wedge t}, U^u_t, z^1)-\si(t, \f^2_{\cd\wedge t}, U^v_t, z^2)\Big\}
%\Big|_{\f^1=X^u, \f^2=X^v\atop z^1=u_t, z^2=v_t}
%\nonumber \\
&=& \int_0^t\Big\{\hE^0\{\d_\th L_s \cd\si(s, \f^1_{\cd\wedge s},U^{\th, v}_s, z^1)\}\Big|_{\f^1=X^{\th, v}, z^1=u^{\th, v}_s}+
 [D\si]^{\th, u, v}_s(\d_\th X_{\cd\wedge s})\ms\nonumber\\
%\Big|_{\f^1=X^u; z^1=u_t}\\
&&\q+\hE^0\{B^{\th,u,v}(s, \f^2_{\cd\wedge s}, z^1) \d_\th U_s\}
\Big|_{\f^2=X^{u}; \atop z^1=u^{\th,v}_s}
%\Big\}\Big|_{ \f^2=X^v; z^1=u_t}
+C^{\th, u,v}_\si(s)(v_s-u_s)\Big\}dB^1_s,
%\Big]\Big\}\Big|_{\f^1=X^1, \f^2=X^2; z^1=u^1_t, z^2=u^2_t}.
 %\nonumber
\eea
where
\bea
\label{ABC}
&& [D\si]^{\th, u, v}_t(\psi)=
\hE^0\Big\{L^u_t\int_0^1D_\f\si(t, \f^2_{\cd\wedge t}+\l (\f^1_{\cd\wedge t}-\f^2_{\cd\wedge t}), U^{\th, v}_t, z^1)(\psi)d\l\Big\}
\Big|_{\f^1=X^{\th,v}, \f^2=X^u, \atop
z^1=u^{\th, v}_t\ \ \ \ \ \ \ },\nonumber\\
&&B^{\th, u,v}(t, \f^2_{\cd\wedge t}, z^1) =L^u_t\int_0^1\pa_y\si(t, \f^2_{\cd\wedge t}, U^u_t+\l(U^{\th,v}_t-U^u_t), z^1)d\l,
\\
&&C^{\th, u,v}_\si(t) =\hE^0\Big\{L^u_t \int_0^1\pa_z\si(t, \f^2_{\cd\wedge t}, U^u_t, z^2+\l(z^1-z^2))d\l \Big\}\Big|_{\f^2=X^u; z^1=u^{\th, v}_t, z^2=u_t}. \nonumber
\eea
Here the integral involving the Fr\'echet derivative $D_\f \si$ is in the sense of Bochner.
Noting that $U^{\th, v}_t=\frac{\hE^0[L^{\th,v}_tX^{\th,v}_t|\cF^Y_t]}{\hE^0[L^{\th,v}_t|\cF^Y_t]}$ and $U^u_t=\frac{\hE^0[L^u_tX^u_t|\cF^Y_t]}{\hE^0[L^u_t|\cF^Y_t]}$, we can easily check that
\bea
\label{deltaU}
\d_\th U_t&=&\frac{\hE^0[L^u_t|\cF^Y_t]\hE^0[L^{\th,v}_tX^{\th,v}_t|\cF^Y_t]-\hE^0[L^{\th,v}_t|\cF^Y_t]\hE^0[L^u_tX^u_t|\cF^Y_t]}{\th\hE^0[L^{\th,v}_t|\cF^Y_t]\hE^0[L^u_t|\cF^Y_t]}\\
&=& \frac{\hE^0[L^u_t|\cF^Y_t]\hE^0[\d_\th L_tX^{\th,v}_t+L^u_t\d_\th X_t|\cF^Y_t]-\hE^0[\d_\th L_t|\cF^Y_t]
\hE^0[L^u_tX^u_t|\cF^Y_t]}{\hE^0[L^{\th,v}_t|\cF^Y_t]\hE^0[L^u_t|\cF^Y_t]}\nonumber\\
&=&\frac{\hE^0[\d_\th L_tX^{\th,v}_t+L^u_t\d_\th X_t|\cF^Y_t]}{\hE^0[L^{\th,v}_t|\cF^Y_t]}-
\frac{\hE^0[\d_\th L_t|\cF^Y_t]}{\hE^0[L^{\th,v}_t|\cF^Y_t]}U^u_t.\nonumber
\eea
Now, sending $\th\to 0$, and assuming that
\bea
\label{KRuv}
K_t=K^{u,v}_t\dfnn \lim_{\th\to 0} \d_\th X^{u,v}_t; \qq R_t=R^{u,v}_t\dfnn\lim_{\th\to 0}\d_\th L^{u,v}_t
\eea
both exist in $L^2(\hQ^0)$, then it follows from (\ref{deltaX})-(\ref{deltaU}) we have, at least formally,
 \bea
\label{Keq0}
    K_t&=&\int_{0}^{t}\Big\{\hE^0[R_s\si(s,\f_{\cd\wedge s}, U^u_s, z)]\Big|_{\f=X^u,z=u_s}+ [D\si]^{u, v}_s(K_{
    \cd\wedge s})
    \nonumber\\
    & &+\hE^0\Big[B^{u,v}(s, \f_{\cd\wedge s}, z)\Big(\frac{\hE^0[R_sX^u_s+L^u_sK_s|\cF^Y_s]}{\hE^0[L^{u}_s|\cF^Y_s]}-
\frac{\hE^0[R_s|\cF^Y_s]}{\hE^0[L^{u}_s|\cF^Y_s]}U^u_s\Big)\Big]\Big|_{\f=X^u;\atop z=u_s}\\
   & &+ C^{u,v}_\si(s)(v_s-u_s)\Big\}dB_s^1, \nonumber
\eea
where
\bea
\label{DBC}
 [D\si]^{u, v}_t(\psi)&\dfnn&\hE^0\{L^u_t  D_\f\si(t, \f_{\cd\wedge t}, U^{u}_t, z)( \psi) \}
\Big|_{\f=X^u; z=u_t},\nonumber\\
 B^{u,v}(t, \f_{\cd\wedge t}, z) &\dfnn&L^u_t\pa_y\si(t, \f_{\cd\wedge t}, U^u_t, z), \\
C^{u,v}_\si(t) &\dfnn&\hE^0\Big\{L^u_t \pa_z\si(t, \f_{\cd\wedge t}, U^u_t, z)\Big\}\Big|_{\f=X^u; z=u_t}. \nonumber
\eea
Observing also that $U^u_t$ is $\cF^Y_t$-measurable, we have
\bea
\label{Buv}
&&\hE^0\Big[B^{u,v}(s, \f_{\cd\wedge s}, z)\Big(\frac{\hE^0[R_sX^u_s+L^u_sK_s|\cF^Y_s]}{\hE^0[L^{u}_s|\cF^Y_s]}-
\frac{\hE^0[R_s|\cF^Y_s]}{\hE^0[L^{u}_s|\cF^Y_s]}U^u_s\Big)\Big]\Big|_{\f=X^u;\atop z=u_s}\nonumber\\
&=&\hE^u\Big[\pa_y\si(s, \f_{\cd\wedge s}, U^u_s, z)\hE^u\{(L^u_s)^{-1}R_s[X^u_s-U^u_s]+K_s|\cF^Y_s\}
%{\hE^0[L^{u}_s|\cF^Y_s]}-
%\frac{\hE^0[R_s|\cF^Y_s]}{\hE^0[L^{u}_s|\cF^Y_s]}U^u_s\Big)
\Big]\Big|_{\f=X^u;\atop z=u_s}\\
&=&\hE^u\Big[(L^u_s)^{-1}\pa_y\si(s, \f_{\cd\wedge s}, U^u_s, z)\{R_s[X^u_s-U^u_s]+L^u_sK_s\}
\Big]\Big|_{\f=X^u;\atop z=u_s}\nonumber\\
&=&\hE^0\Big[\pa_y\si(s, \f_{\cd\wedge s}, U^u_s, z)(R_sX^u_s+L^u_sK_s)-U^u_s\pa_y\si(s, \f_{\cd\wedge s}, U^u_s, z)R_s\Big]\Big|_{\f=X^u;\atop z=u_s}.
\nonumber
\eea
Consequently, if we define
\bea
\label{Psi}
 \Psi(t, \f_{\cd\wedge t}, x, y,z)\dfnn\sigma(t, \f_{\cd\wedge t}, y,z)+\pa_y\sigma(t, \f_{\cd\wedge t}, y,z)(x-y),
\eea
then we can rewrite (\ref{Keq0}) as
\bea
\label{Keq}
K_t&=&\int_{0}^{t}\Big\{\hE^0\Big[\Psi(s,\varphi_{\cd\wedge s}, X^u_s, U^u_s,z)R_s
+\pa_y\sigma(s,\varphi_{\cd\wedge s},U^u_s, z)L_s^uK_s\Big]\Big|_{\varphi=X^{u}; z=u_{s}}\\
    %E^0[L_s^u\sigma_z'(E^u[X_s^u|{\cal{F}}_s^Y],\varphi_{\cd\wedge s},u)]\big|_{\varphi=X^u;\atop u=u_s}
&& \q + [D\si]^{u, v}_s(K_{\cd\wedge s})+C^{u,v}_\si(s)(v_s-u_s)\Big\}dB_s^1.\nonumber
\eea
 Similarly, we can formally write down the SDE for $R$:
\bea
\label{Req}
    R_t =\int_{0}^{t}[R_sh(s,X_{s}^u)+L_s^u\pa_xh(s,X_{s}^u)K_{s}]dY_s,
\q   t\geq0.
\eea

The following theorem is regarding the well-posedness of the SDEs (\ref{Keq}) and (\ref{Req}).
\begin{thm}\label{Theorem 5.2}
Assume that Assumption \ref{Assum1} is in force, and let $u, v\in L^{\infty-}_{\hF^Y}(\hQ^0; [0,T])$ be given. Then,
there is a unique solution $(K, R)\in \scL^{\infty-}_{\hF}(\hQ^0;\hC_T^2)$ to SDEs (\ref{Keq}) and (\ref{Req}).
\end{thm}

{\it Proof.} Let $u, v\in  L^{\infty-}_{\hF^Y}(\hQ^0; [0,T])$ be given. We define $F^1_t(K, R)$ and $F_t^2(K, R)$, $t\in [0,T]$,
to be the right hand side of (\ref{Keq}) and (\ref{Req}), respectively.

We first observe that $F_t^1(0,0)=\int_0^t C^{u,v}_\si(s)(v_s-u_s)dB^1_s$, and $F_t^2(0,0)\equiv 0$, $t\in[0,T]$. Then, for any $p>2$, it holds that
\bea
\label{F1est}
\hE^u\Big[\sup_{0\le s\le t}|F_s^1(0,0)|^p\Big]\le C_p\hE^u\Big[\Big(\int_0^t|v_s-u_s|^2ds\Big)^{p/2}\Big], \qq t\in [0,T].
\eea
Now let $(K^i,R^i)\in \scL^{\infty-}_\hF(\hQ^0;\hC_T)$, $i=1,2$. We define
$\wt K^i\dfnn F_1(K^i, R^i)$,  $\wt R^i\dfnn F_1(K^i, R^i)$, $i=1,2$, and
$\bar K\dfnn K^1-K^2$, $\bar R\dfnn R^1-R^2$, $\hat K\dfnn \wt K^1-\wt K^2$, and $\hat R\dfnn\wt R^1-\wt R^2$.
Then, noting that $\si$, $\pa_y\si$, $y\pa_y\si$,  and $\pa_z\si$ are all bounded, thanks to Assumption \ref{Assum1}, we see that
\beaa
\label{Psiest}
|\Psi(t, \f_{\cd\wedge t}, x,y,z)|\le C(1+|x|), \q (t, x,y,z)\in [0,T]\times \hR^3, ~\f\in\hC_T,
\eeaa
where, and in what follows, $C>0$ is some generic constant which is allowed to vary from line to line. It then follows that
\bea
\label{Psiest2}
&&\Big|\hE^0[\Psi(t, \f_{\cd\wedge t}, X^u_t, U^u_t,z)\bar R_s+\pa_y\si(t, \f_{\cd\wedge t}, U^u_t, z)L^u_t\bar K_t]\Big|
\nonumber\\
&& \le C\hE^0[(1+|X^u_t|)|\bar R_t|+|L^u_t \bar K_t|] \le C\Big[ \hE^0[|\bar K_t|^2+|\bar R_t|^2]\Big]^{1/2}.
\eea
Furthermore, since $D_\f\si$ is also bounded, we have
$|[D\si]^{u,v}_t(\psi)|\le C\sup_{0\le s\le t}|\psi(s)|$, for $\psi\in\hC_T$.
Then from the definition of $\hat K$ and (\ref{Psiest2}) we have, for any $p\ge 2$ and $t\in [0,T]$,
\bea
\label{hatK}
\hE^0\Big[\sup_{0\le s\le t}|\hat K_s|^{2p}\Big] \le C_p\int_0^t\Big(\hE^0[|\bar R_s|^2+|\bar K_s|^2]\Big)^pds
+C_p\int_0^t\hE^0\Big[\sup_{0\le r\le s}|\bar K_r|^{2p}\Big]ds.
\eea
On the other hand, the boundedness of $h$ and $\pa_x h$ implies that, recalling the definition of $\hat R$, for $p\ge 2$ and
$t\in [0,T]$,
\bea
\label{hatR}
\Big(\hE^0\Big[\sup_{s\le t}|\hat R_s|^p\Big)^2&\le & C_p\int_0^t \hE^0[|\bar R_s|^p]^2ds+C_p\int_0^t\hE^0[|L_s^u\bar K_s|^p]^2ds\\
&\le & C_p\int_0^t (\hE^0[|\bar R_s|^p])^2ds+C_p\int_0^t\hE^0[|\bar K_s|^{2p}]ds. \nonumber
\eea
Combining (\ref{hatK}) and (\ref{hatR}) we have, for $t\in[0,T]$,
\beaa
\hE^0\Big[\sup_{0\le s\le t}|\hat K_s|^{2p}\Big]+\Big(\hE^0\Big[\sup_{0\le s\le t}|\hat R_s|^p]\Big)^2\le
C_p\int_0^t\Big(\hE^0\Big[\sup_{0\le r\le s}|\bar K_r|^{2p}\Big]+\Big(\hE^0\Big[\sup_{0\le r\le s}|\bar R_r|^p\Big]^2\Big)ds.
\eeaa
This, together with (\ref{F1est}), enables us to apply  standard SDE arguments to deduce that there is a unique solution
$(K, R)\in  \scL^{\infty-}_{\hF}(\hP;\hC_T)$ of (\ref{Keq}) and (\ref{Req}), such that for all $p\ge 2$, it holds that
 \bea
 \label{KRest}
 \hE^0\big[\|K\|^{2p}_{\hC_T}\big]+\hE^0\big[\|R\|^{2p}_{\hC_T}\big]\le C_p\|v_s-u_s\|^2_{p,2,\hQ^0}.
 \eea
 We leave it to the interested reader, and this completes the proof.
 \qed

\section{Variational Equations}
\setcounter{equation}{0}

In this section we validate the heuristic arguments in the previous section and derive the variational equation
of the optimal trajectory rigorously.
Recall the processes $\d_\th X=\d_\th X^{u,v}$, $\d_\th L=\d_\th L^{u,v}$, and $(K,R)$ defined in the previous
section.  Denote
\begin{equation}\label{7.44}
    \eta_t^\theta\dfnn\d_\theta X_t-K_t, \qq
    \tilde{\eta}_t^\theta\dfnn \d_\th L_t-R_t,\qq  t\in[0,T].
\end{equation}
Our main purpose of this section is to prove the following result.
\begin{prop}
\label{lemma7.2}
Let $(\hP^u, u)=(\hP^{u^*}, u^*)\in \sU_{ad}$  be an optimal control, $(X^u, L^u)$ be the
corresponding solution of (\ref{7.18a}), and let $U^u_t=\hE^u[X^u_t|\cF^Y_t]$, $t\ge 0$.
For any $v\in\sU_{ad}$, let $(K, R) = (K^{u,v}, R^{u,v})$ be the solution of the  linear equations (\ref{Keq}) and (\ref{Req}).
Then, for all $p>1$, it holds that
\begin{equation}\label{7.42}
 \lim_{\th\to0}\hE^0[\|\eta^\th\|^p_{\hC_T}]= \lim_{\theta\rightarrow 0}\hE^0\Big[\sup_{s\in[0,T]} \Big|\frac{X_s^{\th, v}-X_s^u}{\theta}-K_s\Big|^p\Big]=0;
\end{equation}
\begin{equation}\label{7.43}
 \lim_{\th\to0}\hE^0[\|\tilde\eta^\th\|^p_{\hC_T}]=   \lim\limits_{\theta\rightarrow 0}\hE^0\Big[\sup_{s\in[0,T]}\Big|\frac{L_s^{\th, v}-L_s^u}{\theta}-R_s\Big|^p\Big]=0.
\end{equation}
\end{prop}

The proof of Proposition \ref{lemma7.2} is quite lengthy, we shall split it into two parts.

\ms
[{\it Proof of (\ref{7.43})}]. This part is relatively easy.
%Propositio are to show that both $\eta^\th$ and $\tilde\eta^\th$ converges to $0$ in $L^2(\O;\hC([0,T])$.
We note that with
a direct calculation using the equations (\ref{deltaX}) and (\ref{Req}) it is readily seen that $\tilde \eta^\th$ satisfies the
following SDE:
\bea
\label{7.45}
  \tilde{\eta}_t^\theta&=& \int_{0}^{t}\tilde{\eta}_r^{\theta}h(r,X_{ r}^{\theta,v})dY_r+\int_{0}^{t}L_r^u\int_{0}^{1}\pa_xh(r,X_{ r}^u+\lambda\theta(\eta_{ r}^{\theta}+K_{ r}))\eta_{ r}^{\theta }d\lambda d Y_r\nonumber\\
  &&+I_t^{1,\theta}+I_t^{2,\theta},
  \eea
where
\begin{eqnarray*}
I_t^{1,\theta}&=& \int_{0}^{t}R_r(h(r,X_{ r}^{\th, v})-h(r,X_{ r}^u))dY_r;\\
I_t^{2,\theta}&=&\int_{0}^{t}L_r^u\int_{0}^{1}\pa_xh(r,X_{ r}^u+\lambda\theta(\eta_{ r}^{\theta}+K_{ r}))K_{ r} d\lambda dY_r-\int_{0}^{t}L_r^u\pa_xh(r,X_{ r}^u)K_{ r}dY_r.
\end{eqnarray*}
We claim that, for all $p>1$,
%Then, t
\begin{equation}\label{7.46}
  \lim\limits_{\theta\rightarrow0}\hE^u[\sup_{t\in[0,T]}|I_{t}^{1,\theta}|^p]=0,\ \ \
   \lim_{\theta\rightarrow0}\hE^u[\sup_{t\in[0,T]}|I_{t}^{2,\theta}|^p]=0.
\end{equation}
Indeed, note that $dY_t=dB^2_t-h(t,X^u_t)dt$, and $B^2$ is a $\hP^u$-Brownian motion.
Proposition \ref{est1}, together with the bounded and continuity of  $h$ and $\pa_x h$, leads to that, for all $p\ge 2$,
\beaa
\hE^u\Big\{\sup_{t\in[0,T]}|I^{1, \th}_t|^p\Big\}&=&\hE^0\Big\{L^u_T\sup_{t\in[0,T]}\Big|\int_0^tR_s[h(s,X^{\th,v}_s)-h(s,X^u_s)]dY_s\Big|^p
\Big\}\\
&\le &2\hE^u\Big\{\sup_{t\in[0,T]}\Big|\int_0^t R_s[h(s,X^{\th,v}_s)-h(s,X^u_s)]dB^2_s\Big|^p\Big\}\\
&&+2\hE^0\Big\{L^u_T\sup_{t\in[0,T]}\Big|\int_0^tR_s[h(s,X^{\th, v}_s)-h(s,X^u_s)]h(s,X^u_s)ds\Big|^p\Big\}\\
&\le & C_p\hE^0\Big\{L^u_T\int_0^T R^p_s(|X^{\th,v}_s-X^u_s|^p\wedge 1)ds\Big\}\\
&\le & C_p\Big\{\neg\hE^0\big[(L^u_T)^3\big]\neg\Big\}^{\frac13}
\Big\{\hE^0\big[\neg\neg\sup_{s\in[0,T]}|R_s|^{3p}\big]\neg\Big\}^{\frac13}
\Big\{\hE^0 \big[\neg\neg\sup_{s\in[0,T]}(|X^{\th,v}_s-X^u_s|^2\neg\wedge \neg1)\big]\neg\Big\}^{\frac13}\\
&\le & C_p\|u-u^{\th,v}\|_{p,2,\hQ^0}^{\frac23}\le C|\th|^{\frac23},
\eeaa
where we used the following estimate for
%We also would like to point out that in both estimates above we have used the fact that for
any function $f\in L^\infty(\hR)$ bounded by $C_0\ge 1$:
\bea
\label{3pest}
|f(x)-f(x')|^{3p}\le (2C_0 (|f(x)-f(x')|\wedge 1))^{3p}\le (2C_0)^{3p}(|f(x)-f(x')|^2\wedge 1),~ \forall p\ge 2.
\eea
Similarly, we have
\beaa
&&\hE^u\Big\{\sup_{t\in[0,T]}|I^{2, \th}_t|^p\Big\}\\
&=&  \hE^0\Big\{L^u_T\sup_{t\in[0,T]}\Big|\int_0^tL^u_rK_r\Big[\int_0^1\big[
\pa_xh(r, X^u_r+\l\th(\eta^\th_r+K_r))-\pa_xh(r,X^u_r)]d\l\Big]dY_r\Big|^p \Big\}\\
 &\le& C_p \hE^0\Big\{L^u_T\int_0^T|L^u_r|^p|K_r|^p\Big[\int_0^1\big|
\pa_xh(r, X^u_r+\l\th(\eta^\th_r+K_r))-\pa_xh(r,X^u_r)\big|d\l\Big]^pdr\Big\}\\
&\le& C_p\hE^0\Big\{\int_0^T\Big[\int_0^1\big[|\pa_xh(r, X^u_r+\l\th(\eta^\th_r+K_r))-\pa_xh(r,X^u_r)|^2\wedge 1\big]d\l\Big] dr\Big\}^{1/3}.
\eeaa
Here in the above the second inequality follows from (\ref{3pest}) applied to $\partial_x h$, the H\"{o}lder inequality, and the fact that $L^{u},K\in \scL^{\infty-}_{\hF}(\hQ^0;\hC_T)$ (see Theorem \ref{Theorem 5.2}), and the last inequality follows from the $L^p$-estimate (\ref{KRest}).
Now, from (\ref{7.25}), (\ref{Keq}),  and (\ref{Req})
we see that
$$ \hE^0\Big\{\sup_{t\in[0,T]}\big(|\eta^\th_t|^2+|K_t|^2\big)\Big\}\le C, \qq  \th\in(0,1).
$$
Hence, since $\th [\|\eta^{\th}\|_{\hC_T}+\|K\|_{\hC_T}]\to 0$, in probability $\hQ^0$,
as $\th\to0$,  the continuity of $\pa_x h$ and the Bounded Convergence Theorem then imply (\ref{7.46}),
 proving the claim. Recalling (\ref{7.45}), we see that (\ref{7.43}) follows from (\ref{7.46}), provided
(\ref{7.42}) holds, which we now substantiate.

[{\it Proof of (\ref{7.42})}]. This part is more involved. We first rewrite (\ref{dthX}) as follows
\bea
\label{dthX1}
\d_\th X_t
&=& \int_0^t\Big\{\hE^0\{(\tilde\eta^\th_s+R_s)\si(s, \f_{\cd\wedge s},U^{\th,v}_s, z)\}\Big|_{\f=X^{\th,v}, \atop
z=u^{\th,v}_s}+
 [D\si]^{\th, u, v}_s(\eta^\th_{\cd\wedge s}+K_{\cd\wedge s})\ms\nonumber\\
%\Big|_{\f^1=X^u; z^1=u_t}\\
&&\q+\hE^0\{B^{\th,u,v}(s, \f_{\cd\wedge s}, z) \d_\th U_s\}
\Big|_{\f=X^{u}; \atop z=u^{\th,v}_s}
%\Big\}\Big|_{ \f^2=X^v; z^1=u_t}
+C^{\th, u,v}_\si(s)(v_s-u_s)\Big\}dB^1_s.
%\Big]\Big\}\Big|_{\f^1=X^1, \f^2=X^2; z^1=u^1_t, z^2=u^2_t}.
 %\nonumber
\eea
Here $[D\si]^{\th, u,v}$, $B^{\th,u,v}$, and $C^{\th, u,v}$ are defined by (\ref{ABC}). Furthermore, in light of
 (\ref{deltaU}), we can also write:
 %$\d_\th U$ inas
\bea
\label{deltaU1}
\d_\th U_t
%&=&\frac{\hE^0[L^u_t|\cF^Y_t]\hE^0[L^{\th,v}_tX^{\th,v}_t|\cF^Y_t]-\hE^0[L^{\th,v}_t|\cF^Y_t]\hE^0[L^u_tX^u_t|\cF^Y_t]}{\th\hE^0[L^{\th,v}_t|\cF^Y_t]\hE^0[L^u_t|\cF^Y_t]}\\
%&=& \frac{\hE^0[L^u_t|\cF^Y_t]\hE^0[\d_\th L_tX^{\th,v}_t+L^u_t\d_\th X_t|\cF^Y_t]-\hE^0[\d_\th L_t|\cF^Y]
%\hE^0[L^u_tX^u_t|\cF^Y_t]}{\hE^0[L^{\th,v}_t|\cF^Y_t]\hE^0[L^u_t|\cF^Y_t]}\nonumber\\
=\frac{\hE^0[(\tilde\eta^\th_t+R_t)X^{\th,v}_t+L^u_t(\eta^\th_t+K_t)|\cF^Y_t]}{\hE^0[L^{\th,v}_t|\cF^Y_t]}-
\frac{\hE^0[(\tilde\eta^\th_t+R_t)|\cF^Y_t]}{\hE^0[L^{\th,v}_t|\cF^Y_t]}U^u_t.\nonumber
\eea
Plugging this into (\ref{dthX1}) we have
\bea
\label{dthX2}
\d_\th X_t
&=& \int_0^t\Big\{\hE^0\{\tilde\eta^\th_s\si(s, \f_{\cd\wedge s},U^{\th,v}_s, z)\}\Big|_{\f=X^{\th,v}, \atop
z=u^{\th,v}_s}+
 [D\si]^{\th, u, v}_s(\eta^\th_{\cd\wedge s})\ms\nonumber\\
%\Big|_{\f^1=X^u; z^1=u_t}\\
&&\q+\hE^0\Big\{B^{\th,u,v}(s, \f_{\cd\wedge s}, z) \Big[\frac{\hE^0[\tilde\eta^\th_sX^{\th,v}_s+L^u_s\eta^\th_s|\cF^Y_s]}{\hE^0[L^{\th,v}_s|\cF^Y_s]}-
\frac{\hE^0[\tilde\eta^\th_s|\cF^Y_s]}{\hE^0[L^{\th,v}_s|\cF^Y_s]}U^u_s\Big]\Big\}
\Big|_{\f=X^{u}; \atop z=u^{\th,v}_s}\Big\}dB^1_s\nonumber\\
&&+ \int_0^t\Big\{\hE^0\{R_s\si(s, \f_{\cd\wedge s},U^{\th}_s, z)\}\Big|_{\f=X^{\th,v}, \atop
z=u^{\th,v}_s}+ [D\si]^{\th, u, v}_s(K_{\cd\wedge s})\ms\nonumber\\
&&\q+\hE^0\Big\{B^{\th,u,v}(s, \f_{\cd\wedge s}, z) \Big[\frac{\hE^0[R_sX^{\th,v}_s+L^u_sK_s|\cF^Y_s]}{\hE^0[L^{\th,v}_s|\cF^Y_s]}-
\frac{\hE^0[R_s|\cF^Y_t]}{\hE^0[L^{\th,v}_s|\cF^Y_s]}U^u_s\Big]\Big\}
\Big|_{\f=X^{\th,v}; \atop z=u^{\th,v}_s}\nonumber\\%\nonumber
&&\q+C^{\th, u,v}_\si(s)(v_s-u_s)\Big\}dB^1_s.\nonumber
\eea
Now, recalling (\ref{Keq}) (or more conveniently, (\ref{Keq0})) we have
 \bea
\label{eta}
   \eta_{t}^{\theta}&=& \d_\th X_t- K_t=\int_0^t\Big\{\hE^0\{\tilde\eta^\th_s\si(s, \f_{\cd\wedge s},U^{\th,v}_s, z)\}\Big|_{\f=X^{\th,v}, \atop
z=u^{\th,v}_s}+
 [D\si]^{\th, u, v}_s(\eta^\th_{\cd\wedge s})\ms\nonumber\\
&&+\hE^0\Big\{B^{\th,u,v}(s, \f_{\cd\wedge s}, z) \Big[\frac{\hE^0[\tilde\eta^\th_sX^{\th,v}_s+L^u_s\eta^\th_s|\cF^Y_s]}{\hE^0[L^{\th,v}_s|\cF^Y_s]}-
\frac{\hE^0[\tilde\eta^\th_s|\cF^Y_s]}{\hE^0[L^{\th,v}_s|\cF^Y_s]}U^u_s\Big]\Big\}
\Big|_{\f=X^{u}; \atop z=u^{\th,v}_s}\Big\}dB^1_s\nonumber\\
&&+ I^{3, \th, 1}_t+ I^{3, \th, 2}_t+ I^{3, \th, 3}_t+ I^{3, \th, 4}_t,
\eea
where, for $t\in[0,T]$,
\bea
\label{I3th1}
I^{3,\th,1}_t&\dfnn&\int_{0}^{t}\hE^0\big\{R_s\big[\si(s,\f^1_{\cd\wedge s}, U^{\th,v}_s, z^1)-\si(s,\f^2_{\cd\wedge s}, U^u_s, z^2)\big]\big\}\Big|_{\f^1=X^{\th,v},z^1=u^{\th,v}_s\atop \f^2=X^u,z^2=u_s}dB^1_s;\nonumber\\
I^{3,\th,2}_t&\dfnn&\int_{0}^{t}\hE^0\big\{[D\si]^{\th,u, v}_s(K_{\cd\wedge s})-
 [D\si]^{u, v}_s(K_{\cd\wedge s})\big\}dB^1_s;\\
 \nonumber
\eea
\bea
 I^{3,\th,3}_t   &\dfnn&\int_{0}^{t}\Big\{\hE^0\Big\{B^{\th,u,v}(s, \f_{\cd\wedge s}, z) \Big(\frac{\hE^0[R_sX^{\th,v}_s+L^u_sK_s|\cF^Y_s]}{\hE^0[L^{\th,v}_s|\cF^Y_s]}-
\frac{\hE^0[R_s|\cF^Y_t]}{\hE^0[L^{\th,v}_s|\cF^Y_s]}U^u_s\Big)\Big\}\Big|_{\f=X^{u}; \atop z=u^{\th,v}_s}\nonumber\\
&&-\hE^0\Big\{B^{u,v}(s, \f_{\cd\wedge s}, z)\Big(\frac{\hE^0[R_sX^u_s+L^u_sK_s|\cF^Y_s]}{\hE^0[L^{u}_s|\cF^Y_s]}-
\frac{\hE^0[R_s|\cF^Y_s]}{\hE^0[L^{u}_s|\cF^Y_s]}U^u_s\Big)\Big\}\Big|_{\f=X^u;\atop z=u_s} \Big\}dB^1_s  \nonumber\\
I^{3,\th,4}_t   &\dfnn&\int_{0}^{t}\hE^0[C^{\th, u,v}_\si(s)(v_s-u_s)- C^{u,v}_\si(s)(v_s-u_s)]dB_s^1. \nonumber
\eea

We have the following lemma.
\begin{lem}
\label{I3thconv}
Suppose that Assumption \ref{Assum1} holds. Then, for all $p>1$,
\bea
\label{I3conv}
\lim_{\th\to0}\hE^0\Big\{\sup_{0\le t\le T}|I^{3,\th,i}_t|^p\Big\}=0, \qq i=1,\cds, 4.
\eea
\end{lem}

{\it Proof.}
We first recall that $U^{\th,v}_s\dfnn \hE^{\th,v}[X_{s}^{\th,v}|\cF_{s}^{Y}]$ and $U^u_s\dfnn\hE^{u}[X_{s}^{u}|\mathcal{F}_{s}^{Y}]$.
Using the Kallianpur-Strieble formula we have
\bea
\label{J12th}
\hE^0\int_0^T|U^{\th,v}_s-U^u_s|^pds&\le&C_p\Big\{\hE^{0}\int_0^T\Big|\frac{\hE^{0}[L^{\th,v}_sX_{s}^{\th,v}|\cF_{s}^{Y}]}{\hE^0[L^{\th,v}_s|\cF^Y_s]}-\frac{\hE^{0}[L^u_sX_{s}^{u}|\mathcal{F}_{s}^{Y}]}{\hE^0[L^{\th,v}_s|\cF^Y_s]}\Big|^pds\nonumber\\
%&&+\hE^{0}\int_0^T\Big|\frac{\hE^{0}[L^u_sX_{s}^{\th}|\mathcal{F}_{s}^{Y}]}{\hE^0[L^\th_s|\cF^Y_s]}-
%\frac{\hE^{0}[L^u_sX_{s}^{\th}|\mathcal{F}_{s}^{Y}]}{\hE^0[L^u_s|\cF^Y_s]}\Big|^2ds\\
&&+\hE^{0}\int_0^T\Big|\frac{\hE^{0}[L^u_sX_{s}^{u}|\mathcal{F}_{s}^{Y}]}{\hE^0[L^{\th,v}_s|\cF^Y_s]}-
\frac{\hE^{0}[L^u_sX_{s}^{u}|\mathcal{F}_{s}^{Y}]}{\hE^0[L^u_s|\cF^Y_s]}\Big|^pds\Big\}\\
&\dfnn& C_p\{J^{1}_\th+J^{2}_\th\}. \nonumber
\eea
%where $J^{3,2,i}_\th$ are defined in an obvious way.
We now estimate $J^1_\th$ and $J^2_\th$ respectively.
First note that, for any $p>1$, we can find a constant $C_p>0$ such that for any $\th \in(0,1)$ and $u\in \sU_{ad}$,
$$ \hE^0[(L^{\th,v}_s)^p]+\hE^0[(L^{\th,v}_s)^{-p}]+\hE^0[(L^u_s)^p] \le C_p.
$$
Thus, applying the H\"older and Jensen inequalities as well as Proposition \ref{est1}, we have, for any $p>1$,
and $\th\in(0,1)$,
\bea
\label{dL1}
&&\hE^0\int_0^T\Big|\frac{\hE^0[L^{\th,v}_sX^{\th,v}_s|\cF^Y_s]-\hE^0[L^u_sX^u_s|\cF^Y_s]}{\hE^0[L^{\th,v}_s|\cF^Y_s]}\Big|^pds \le
\int_0^T\hE^0\Big\{\frac{|L^{\th,v}_sX^{\th,v}_s-L^u_sX^u_s|^p}{\hE^0[
L^{\th,v}_s|\cF^Y_s]^p}\Big\}ds\nonumber\\
&\le& \int_0^T\Big\{\{\hE^0|L^{\th,v}_sX^{\th,v}_s-L^u_sX^u_s|^2\}^{1/2}\cd \Big\{\hE^0\Big[\frac{|L^{\th,v}_sX^{\th,v}_s-L^u_sX^u_s|^{2p-2}}{\hE^0[L^{\th,v}_s|\cF^Y_s]^{2p}}\Big]\Big\}^{1/2}\Big\}ds\\
&\le& \int_0^t\{\hE^0|L^{\th,v}_sX^{\th,v}_s-L^u_sX^u_s|^2\}^{1/2}\cd \Big\{\hE^0[|L^{\th,v}_sX^{\th,v}_s-L^u_sX^u_s|^{2p-2}]
\hE^0\big[[L^{\th,v}_s]^{-2p}|\cF^Y_s\big]\Big\}^{1/2}ds\nonumber\\
&\le& C_p\th\|u-v\|_{2,2,\hQ^0}.\nonumber
\eea
Similarly, one can also argue that, for any $p>1$, the following estimates hold:
\bea
\label{dL2}
\hE^0\int_0^T\Big| \frac1{\hE^0[L^{\th,v}_s|\cF^Y_s]}-\frac1{\hE^0[L^u_s|\cF^Y_s]}
%\Big]\hE^0[L_s^{u}X^u_s|\cF^Y_s]
\Big|^pds \le C_p\th\|u-v\|_{2,2,\hQ^0}, \q\th\in(0,1).
\eea
Clearly, (\ref{dL1}) and (\ref{dL2}) imply that   $J^1_\th+J^2_\th\le C_p\th\|u-v\|_{2,2,\hQ^0}$, for some constant $C_p>0$, depending only
on $p$, the Lipschitz constant of the coefficients, and $T$. Therefore we have
%together with Dominated Convergence Theorem, imply (\ref{Uth}).
\bea
\label{Uth}
\hE^0\int_0^T|U^{\th,v}_s-U^u_s|^pds \le C_p\th\|u-v\|_{2,2,\hQ^0}\to 0, \qq \mbox{as $\th\to0$.}
\eea

We can now prove (\ref{I3conv}) for  $i=1,\cds, 4$. First, by Burkholder-Davis-Gundy inequality we have
\beaa
\label{I3th1}
\hE^0[\sup_{0\le t\le T}|I^{3,\th,1}_t|^2]\le C\int_0^T\hE^0\Big|\hE^0\big\{R_s\big[\si(s,\f^1_{\cd\wedge s}, U^{\th,v}_s, z^1)-\si(s,\f^2_{\cd\wedge s}, U^u_s, z^2)\big]\big\}\Big|_{\f^1=X^{\th,v},z^1=u^{\th,v}_s\atop \f^2=X^u,z^2=u_s}\Big|^2ds.
\eeaa
Since $\si$ is bounded and Lipschitz continuous in $(\f, y, z)$, it follows from Proposition \ref{est1} and (\ref{Uth})
that $\lim_{\th\to0} \hE^0[\sup_{0\le t\le T}| I^{3,\th, 1}_t|^2]=0$.
By the similar arguments using the continuity of $D_\f\si$ and that of $\pa_z\si$, respectively,  it is not hard to show that,
for all $p>1$,
$$\lim_{\th\to 0}\hE^0[\sup_{0\le t\le T}| I^{3,\th, 2}_t|^p]=0; \qq \lim_{\th\to0}\hE^0[\sup_{0\le t\le T}| I^{3,\th, 4}_t|^p]=0.
$$

It remains to prove the convergence of $I^{3,\th,3}$. To this end, we note that, for any $p>1$,
\bea
\label{est2}
\hE^0\Big[\sup_{s\in[0,T]}\big(|R_s|^p+|K_s|^p\big)\Big]&\le& C_p,
\eea
and by (\ref{Uth}) we have, for $p>1$,
\bea
\label{B3th}
\lim_{\th\to0}\hE^0\int_{0}^{T}\Big|\hE^{0}\big\{\big|B^{\th,u,v}(s,\f_{\cd\wedge s},z)
-B^{u,v}(s, \varphi_{\cd\wedge s}, z^1)\big|^2\big\}\big|_{\varphi=X^{u},z=u^{\th,v}_s\atop
z^1=u_s\ \ \ \ \ \ \ }\Big|^pds=0.
\eea
This, together with (\ref{dL2}), (\ref{Uth}), an estimate similar to (\ref{dL1}), and Proposition \ref{est1}, yields that
$\lim_{\th\to 0}\hE^0[\sup_{0\le t\le T}| I^{3,\th, 3}_t|^2]=0$,
proving the lemma.
\qed

We now continue the proof of (\ref{7.42}). First we rewrite (\ref{eta})  as
 \bea
\label{eta1}
   \eta_{t}^{\theta}&=&\int_0^t\Big\{\hE^0\{\tilde\eta^\th_s\si(s, \f_{\cd\wedge s},U^{\th,v}_s, z)\}\Big|_{\f=X^{\th,v}, \atop
z=u^{\th,v}_s}+
 [D\si]^{\th, u, v}_s(\eta^\th_{\cd\wedge s})\ms\nonumber\\
%\Big|_{\f^1=X^u; z^1=u_t}\\
&&+\hE^0\Big\{B^{\th,u,v}(s, \f_{\cd\wedge s}, z) \Big[\frac{\hE^0[\tilde\eta^\th_sX^{\th,v}_s+L^u_s\eta^\th_s|\cF^Y_s]}{\hE^0[L^{u}_s|\cF^Y_s]}-
\frac{\hE^0[\tilde\eta^\th_s|\cF^Y_s]}{\hE^0[L^{u}_s|\cF^Y_s]}U^u_s\Big]\Big\}
\Big|_{\f=X^{u}; \atop z=u^{\th,v}_s}\Big\}dB^1_s\nonumber\\
&&+I^{3,\th,0}_t+\sum_{i=1}^4 I^{3, \th, i}_t,
\eea
where
\beaa
I^{3,\th,0}_t&\dfnn&\int_0^t\hE^0\Big\{B^{\th,u,v}(s, \f_{\cd\wedge s}, z) \Big[\frac{\hE^0[\tilde\eta^\th_sX^{\th,v}_s+L^u_s\eta^\th_s|\cF^Y_s]}{\hE^0[L^{\th,v}_s|\cF^Y_s]}-
\frac{\hE^0[\tilde\eta^\th_s|\cF^Y_s]}{\hE^0[L^{\th,v}_s|\cF^Y_s]}U^u_s\Big]\Big\}
\Big|_{\f=X^{u}; \atop z=u^{\th,v}_s}\nonumber\\
&&-B^{\th,u,v}(s, \f_{\cd\wedge s}, z) \Big[\frac{\hE^0[\tilde\eta^\th_sX^{\th,v}_s+L^u_s\eta^\th_s|\cF^Y_s]}{\hE^0[L^{u}_s|\cF^Y_s]}-
\frac{\hE^0[\tilde\eta^\th_s|\cF^Y_s]}{\hE^0[L^{u}_s|\cF^Y_s]}U^u_s\Big]\Big\}
\Big|_{\f=X^{u}; \atop z=u^{\th,v}_s}\Big\}dB^1_s\nonumber\\
\eeaa
We note that with the same argument as before one shows that $\lim_{\th\to0}\hE^0[\sup_{0\le t\le T}|I^{3,\th,0}_t|^2]=0$. On the other hand, similar to (\ref{Buv}) one can argue that
\beaa
\label{Buv1}
&&\hE^0\Big[B^{\th,u,v}(s, \f_{\cd\wedge s}, z)\Big(\frac{\hE^0[\tilde\eta^\th_sX^{\th,v}_s+L^u_s\eta^\th_s|\cF^Y_s]}{\hE^0[L^{u}_s|\cF^Y_s]}-
\frac{\hE^0[\tilde\eta^\th_s|\cF^Y_s]}{\hE^0[L^{u}_s|\cF^Y_s]}U^u_s\Big)\Big]\Big|_{\f=X^u;\atop z=u^{\th,v}_s}\nonumber\\
&=&\hE^0\Big[\int_0^1\pa_y\si(s, \f_{\cd\wedge s}, U^u_s+\l (U^{\th,v}_s-U^u_s), z)d\l \cd (\tilde\eta^\th_sX^{\th,v}_s+L^u_s\eta^\th_s-U^u_s
\tilde\eta^\th_s)\Big]\Big|_{\f=X^u;\atop z=u^{\th,v}_s}.
\eeaa
 Consequently, we have
 \beaa
\label{eta2}
\eta_{t}^{\theta}&=&\int_0^t\Big\{\hE^0\{\a^{1,\th}_s(\f^1_{\cd\wedge s},\f^2_{\cd\wedge s}, z)\tilde\eta^\th_s\}\Big|_{\f^1=X^{\theta,v},\f^2=X^{u},\atop z=u^{\th,v}_s\ \ \ \ \ \ \ \ }+\hE^0\{\a^{2,\th}_s(\f^2_{\cd\wedge s}, z)\tilde\eta^\th_s\}\Big|_{\f^2=X^{u},\atop z=u^{\th,v}_s}\Big\}dB^1_s\\
& &+\int_0^t\Big\{\hE^0\{\b^\th_s(\f^2_{\cd\wedge s}, z)\eta^\th_s\}\Big|_{\f^2=X^u \atop z=u^{\th,v}_s}+[D\si]^{\th, u, v}_s(\eta^\th_{\cd\wedge s})\Big\}dB^1_s+I^{3,\th}_t,
\eeaa
where $I^{3,\th}_t=\sum_{i=0}^4 I^{3, \th, i}_t$, and
 \beaa
&&\a^{1,\th}_s(\f^1_{\cd\wedge s},\f^2_{\cd\wedge s}, z)\dfnn\int_0^1D_\varphi\sigma(s,\f^2_{\cd\wedge s}+\lambda(\f^1_{\cd\wedge s}-\f^2_{\cd\wedge s}),U^{\theta,v}_s,z)(\f^1_{\cd\wedge s}-\f^2_{\cd\wedge s})d\lambda,\\
&&\a^{2,\th}_s(\f^2_{\cd\wedge s}, z)\dfnn \sigma(s,\f^2_{\cd\wedge s},U^{\theta,v}_s,z)+\int_0^1\partial_y\sigma(s,\f^2_{\cd\wedge s},U_s^u+\lambda(U^{\theta,v}_s-U_s^u),z)d\lambda(U^{\theta,v}_s-U_s^u);\\
&&\b^\th_s(\f^2_{\cd\wedge s}, z)\dfnn L_s^u\int_0^1 \partial_y\sigma(s,\f^2_{\cd\wedge s},U_s^u+\lambda(U^{\theta,v}_s-U_s^u),z)d\lambda.
\eeaa

\noindent Notice that
$$|\a^{1,\th}_s(\f^1_{\cd\wedge s},\f^2_{\cd\wedge s}, z)| +|\a^{2,\th}_s(\f^2_{\cd\wedge s}, z)|\le C(1+|\f^1_{\cd\wedge s}|+|\f^2_{\cd\wedge s}|+|U_s^{\th, v}|+|U^u_s|),\ |\b_s^\th(\f_{\cd\wedge s}, z)|\le CL^u_s.
$$
 Now by the Burkholder and Cauchy-Schwartz inequalities we have, for all $p\ge 2$, $t\in[0,T]$,
 \beaa
 \hE^0\Big[\sup_{s\in[0,t]}|\eta^\th_s|^{2p}\Big]\le C_p \Big\{\hE^0[\|I^{3,\th}\|^{2p}_{\hC_T}]+\hE^0\Big\{\Big[\int_0^t\Big(\hE^0[|\eta^\th_s|^2
 +|\tilde\eta^\th_s|^2]+\sup_{r\in[0,s]}|\eta^\th_s|^2\Big)ds\Big]^p\Big\}\Big\},  \eeaa
 and from Gronwall's inequality one has
 \bea
 \label{etaest}
 \hE^0\Big[\sup_{s\in[0,t]}|\eta^\th_s|^{2p}\Big]\le C_p\Big\{\hE^0\Big[\|I^{3,\th}\|^{2p}_{\hC_T}+\int_0^t\Big(\hE^0[|\tilde\eta^\th_s|^p]
 \Big)^2ds\Big\}, \q t\in[0,T].
 \eea
 On the other hand, setting $I^\th_t \dfnn I^{1, \th}_t+I^{2,\th}_t$, $t\in[0,T]$, we have from (\ref{7.45}) that, for $p\ge 2$,
 \beaa
 \hE^0\Big[\sup_{s\in[0,t]}|\tilde\eta^\th_s|^p\Big]\le C_p\Big\{\hE^0[\|I^\th\|^p_{\hC_T}]+\int_0^t\hE^0[|\tilde\eta^\th_s|^p]ds
 +\int_0^t\big(\hE^0[|\eta^\th_s|^{2p}]\big)^{1/2}ds\Big\}, ~t\in[0,T].
 \eeaa
Then Gronwall's inequality leads to that
\bea
\label{tetaest}
\big(\hE^0\Big[\sup_{s\in[0,t]}|\tilde\eta^\th_s|^p\Big]\Big)^2\le C_p\Big\{\big(\hE^0\|I^\th\|^p_{\hC_T}\big)^2+
\int_0^t\hE^0[|\eta^\th_s|^{2p}]ds\Big\}, \q t\in[0,T].
\eea
Combining (\ref{etaest}), (\ref{tetaest}), applying (\ref{7.46}) and Lemma \ref{I3thconv} as well as the Gronwall inequality, we can easily
deduce  (\ref{7.42}) by sending $\th\to0$. Consequently, (\ref{7.43}) holds as well.
\qed

From Proposition 6.1, (\ref{deltaU}) and the above development we also obtain the following corollary.

\begin{cor}
We assume that Assumption \ref{Assum1} holds. Then, for all $p>1$,
$$\lim_{\th\to 0}\hE^0[\|\delta_\theta U-\overline{V}\|_{\hC_T}^p]=\lim_{\th\to 0}\hE^0[\sup_{0\le s\le T}|\frac{U_s^{\theta,v}-U_s^u}{\theta}-\overline{V}_s|^p]=0,
$$
where
$$\overline{V}_t\dfnn \frac{\hE^0[R_tX_t^u+L_t^uK_t|{\cal{F}}_t^Y]}{\hE^0[L_t^u|{\cal{F}}_t^Y]}-\frac{\hE^0[R_t|{\cal{F}}_t^Y]}{\hE^0[L_t^u|{\cal{F}}_t^Y]}U_t^u,\ t\in [0, T]. $$
\end{cor}

\section{ Stochastic Maximum Principle}
\setcounter{equation}{0}

We are now ready to study the Stochastic Maximum Principle. The main task will be to determine the appropriate {\it adjoint equation}, which we expect to be a backward stochastic differential equation of Mean-field type. We begin with
a simple analysis. Suppose that $u=u^*$ is an optimal control, and for any $ v\in\sU_{ad}$, we define $u^{\th, v}$
by (\ref{7.22}). Then we have
\bea
\label{7.23}
0&\le &\frac{J(u^{\theta, v})-J(u)}{\th}\nonumber\\
&=&\frac1\th ~\hE^0\Big\{\hE^0[L^{\theta, v}_T\Phi(x, U^{\theta, v}_T)]|_{x=X^\th_T}-\hE^0[L^u_T\Phi(x, U^u_T)]|_{x=X^u_T}\\
&&+\int_0^T\big[\hE^0[L^{\theta, v}_sf(s, \f_{\cd\wedge s}, U^{\theta, v}_s, z)]|_{\f=X^{\theta, v},\atop z=u^{\th,v}_s}-\hE^0[L^u_sf(s, \f_{\cd\wedge s}, U^u_s, z)]|_{\f=X^u,\atop z=u_s}\big]ds\Big\}.\nonumber
\eea
Now, repeating the same analysis as that in Proposition \ref{est1},  then sending $\th\to 0$, it follows from
Propositions \ref{est1},  \ref{lemma7.2} and  the continuity of the functions $\Phi$ and $f$ that
\bea
\label{terminal}
0&\leq& \hE^0[K_{T}\xi]+\hE^0[R_T\Th]
%\hE^{0}[L_{T}^{u}\pa_z\Phi(x,U^u_T)]|_{x=X_{T}^{u}}+K_TL^u_T\hE^0[\pa_y \Phi(X_{T}^{u}, y)]|_{y=U^u_T}
%\nonumber\\
%&&\qq
%\hE^{0}[\Phi(X_{T}^{u}, y)]\big|_{y=U^u_T}+R_T\hE^0[\pa_y \Phi(X_{T}^{u}, y)]|_{y=U^u_T}(X^u_T-U^u_T)
+ \hE^{0}\Big\{\int_0^T\Big\{\hE^0[R_sf(s,\f_{\cd\wedge s}, U^u_s,z)]|_{\f=X^{u}, z=u_s}\nonumber\\
&&\qq+\hE^0[\pa_yf(s,\f_{\cd\wedge s}, U^u_s, z)(X^u_s-U^u_s)R_s+L^u_sK_s]|_{\f=X^{u}, z=u_s}\\
&&\qq+\hE^0[L^u_s D_\f f(s,\f_{\cd\wedge s}, U^u_s,z)(\psi_{\cd\wedge s})]|_{\f=X^{u}, z=u_s,\psi=K}\nonumber\\
&&\qq +\hE^0[L^u_s \pa_z f(s,\f_{\cd\wedge s}, U^u_s,z)]|_{\f=X^{u}, z=u_s}(v_s-u_s)\Big\}ds\Big\},\nonumber
\eea
where
\bea
\label{pTQT}
\xi&\dfnn&\hE^{0}[L_{T}^{u}\pa_x\Phi(x,U^u_T)]|_{x=X_{T}^{u}}+L^u_T\hE^0[\pa_y \Phi(X_{T}^{u}, y)]|_{y=U^u_T},\nonumber\\
\Th&\dfnn&\hE^{0}[\Phi(X_{T}^{u}, y)]\big|_{y=U^u_T}+(X^u_T-U^u_T)\hE^0[\pa_y \Phi(X_{T}^{u}, y)]|_{y=U^u_T}.
\eea

We now consider the adjoint equations that take the following form of backward SDEs on the reference space
$(\O, \cF, \hQ^0)$:
\bea
\label{BSDE}
\left\{\ba{lll}dp_t= -\a_tdt+d\G_t+q_tdB^1_t+\widetilde q_t dY_t, \qq p_T=\xi,\\
dQ_t= -\b_tdt+d\Si_t+M_tdB^1_t+\widetilde M_t dY_t, \qq Q_T=\Th.
\ea\right.
\eea
Here the coefficients $\a, \b$ as well as the two bounded
variation processes $\G$ and $\Si$ are to be determined. Applying It\^o's formula and recalling the variational equations
(\ref{Keq}) and (\ref{Req}), we can easily derive (denote $U^u_t=\hE^{u}[X_{t}^{u}|{\mathcal{F}}_{t}^{Y}]$,
$t\in[0,T]$)
%, then a direct computation shows that
\bea
\label{pKQR}
&& \hE^0[\xi K_T]+\hE^0[\Th R_T] \nonumber\\
&=&\int_0^T\Big\{-\hE^0[K_s\a_s]-\hE^0[R_s\b_s]+\hE^0\Big[q_s\hE^0[R_s\si(s, \f_{\cd\wedge s}, U^u_s, z)]\big|_{\f=X^u, z=u_s}\Big]
\nonumber\\
&&+\hE^{0}\Big[q_{s}\hE^0\Big[\pa_y\si(s,\f_{z\cd\wedge s}, U^{u}_s,z)[(X^{u}_s-U^u_s)R_s+L^u_sK_s]\Big]
\Big|_{\f=X^u, z=u_s}\Big]\\
&&+\hE^{0}\big[q_s[D\si]^{u,v}_s(K_{\cd\wedge s}) +q_s C^{u,v}_\si(s)(v_s-u_s)+\wt M_sR_sh(s, X^u_s)
+\wt M_{s}K_sL^u_s\pa_xh(s, X^u_s)]\Big\}ds \nonumber\\
&&+\hE^0\Big\{\int_0^T[K_sd\G_s+R_sd\Si_s]\Big\},\nonumber
\eea
where $[D\si]^{u,v}$ and $C^{u,v}$ are defined by (\ref{DBC}).

By Fubini's Theorem we see that
\bea
\label{Fubini}
\left\{\ba{llll}
\hE^0\big[q_s\hE^0[R_s\si(s, \f{\cd\wedge s}, U^u_s, z)]\big|_{\f=X^u, z=u_s}\big]=
\hE^{0}\big[R_s\hE^0[q_{s}\si(s, X_{\cd\wedge s}, y, u_s)]\big|_{y=U^{u}_s}\big]; \ms\\
\hE^{0}\Big[q_{s}\hE^0\Big[\pa_y\si(s,\f_{z\cd\wedge s}, U^{u}_s,z)[(X^{u}_s-U^u_s)R_s+L^u_sK_s]\Big]\Big|_{\f=X^u, z=u_s}\Big]\\
\qq=\hE^{0}\Big[\hE^0\big[q_s\pa_y\si(s,X_{z\cd\wedge s}, y,u_s)]\big|_{y=U^u_s}[(X^{u}_s-U^u_s)R_s+L^u_sK_s]\Big].
\ea\right.
\eea
%{\mathcal{F}}_{s}^{Y}],\varphi,u):={L_{s}^{u}}\cdot\sigma'_x(E^{u}[X_{s}^{u}|{\mathcal F}_{s}^{Y}],\varphi,u);\\
Furthermore, in light of definition of $[D\si]^{u,v}$ ((\ref{DBC})), if we denote, for fixed $(t,\f, z)$,
\bea
\label{mu0}
\m^0_\si(t, \f_{\cd\wedge t}, z)(\cd)\dfnn \hE^0[L^u_t D_\f\si(t, \f_{\cd\wedge t}, U^u_t, z)](\cd)\in \sM[0,T],
\eea
where $\sM[0,T]$ denotes all the Borel measures on $[0,T]$, then we can write
\bea
\label{DsiK}
[D\si]^{u,v}_t(K_{\cd\wedge t})=\hE^0\big[L^u_t D_\f\si(t, \f_{\cd\wedge t}, U^u_t, z)(\psi)]\big|_{\f=X^u, z=u_t, \atop
\psi=K_{\cd\wedge t}\ \ \ \ }=\int_0^tK_r\m^0_\si(r, X^u_{\cd\wedge r}, u_r)(dr).
\eea
%Fr\'echet derivative as a measure

%{\color{red}
Let us now argue that a similar Fubini Theorem argument holds for the random measure $\m^0_\si(t, X^u_{\cd\wedge t}, u_t)
(\cd)$.
First, for a given process $q\in L^2_\hF(\hQ^0; [0,T])$, consider the following finite variation (FV) process (in fact, under Assumption \ref{Assum1}, {\it integrable variation} (IV) process):
%.:
\bea
\label{A}
A^\si_t\dfnn \int_0^T\int_0^{t\wedge s}q_s\m^0_\si(s, X_{\cd\wedge s}^u, u_s)(dr)ds, \qq t\in[0,T].
\eea
%It is easy to see that $(t,\o)\mapsto A_t(\o)\dfnn \n([0,t])(\o)$ is a
%n the sense that, for any bounded measurable function $\f\in L^0([0,T])$,
%$$ \int_0^T\f(r)\n(dr)=\int_0^T\f(r)\int_r^Tq_s\m^0_\si(s, X_{\cd\wedge s}, u_s)(dr)ds.
It is easy to check, as a (randomized) signed measure on $[0,T]$, it holds $\hQ^0$-almost surely that $dA^\si_t=\int_t^Tq_s\m^0_\si
(s, X^u_{\cd\wedge s}, u_s)(dt)ds$. We note that being a ``raw FV" process, the process $A^\si$ is not $\hF$-adapted. We now consider its {\it dual predictable projection}:
\bea
\label{dpA}
^p\neg\Big(\int_t^Tq_s\m^0_\si(s, X^u_{\cd\wedge s}, u_s)(dt)ds\Big)\dfnn d[^p\neg\neg A^\si_t], \q t\in[0,T].
 \eea
We remark that $d[^p\neg\neg A_t]$ is a predicable random measure that can be formally understood as
 \beaa
\label{dpA1}
d[ ^p\neg\neg A^\si_t]=\hE^0[dA^\si_t|\cF_{t-}]=\hE^0\Big[\int_t^Tq_s\m^0_\si(s, X^u_{\cd\wedge s}, u_s)(dt)ds\Big|\cF_{t-}\Big], \q t\in[0,T].
 \eeaa

Using the definition of dual predicable projection and (\ref{DsiK}), we see  that, for the continuous process $K\in
L^2_\hF(\hQ^0;\hC_T)$,
\bea
\label{Duvsi}
\int_0^T\hE^0[q_s[D\si]^{u,v}_s(K_{\cd\wedge s})]ds&=&\int_0^T\hE^0\Big[q_s\int_0^sK_r\m^0_\si(r,
X^u_{\cd\wedge r}, u_r)(dr)\Big]ds\nonumber\\
&=& \hE^0\Big[\int_0^TK_r dA^\si_r\Big]=\hE^0\Big[\int_0^T K_r d[^p\neg\neg A^\si_r]\Big]\\
&=& \hE^0\Big[\int_0^TK_r \, ^p\neg\Big(\int_r^Tq_s\m^0_\si(s, X^u_{\cd\wedge s}, u_s)(dr)ds \Big)\Big].
%\\&=& \hE^0\Big\{\int_0^T\hE^0\Big\{\int_r^Tq_s\m^0_\si(s, X^u_{\cd\wedge s}, u_s)(dr)ds\Big|\cF_r\Big\} K_r\Big\}.
\nonumber
\eea
Similarly, we  denote $A^f_t \dfnn \int_0^T\int_0^{t\wedge s}\m^0_f(s, X^u_{\cd\wedge s}, u_s)(dr)ds$, $t\in[0,T]$;
and denote its dual predicable projection by $^p\neg\Big(\int_t^T\m^0_f(s, X^u_{\cd\wedge s}, u_s)(dt)ds\Big) = d[^p\neg\neg  A^f_t]$, $t\in[0,T]$.
% \eea
% \bea
%  \int_0^T\int_0^{t\wedge s}q_s\m^0_\si(s, X_{\cd\wedge s}, u_s)(dr)ds=\int_0^t \int_r^Tq_s\m^0_\si(s, X_{\cd\wedge s}, u_s)(dr)ds,  ~t\in[0,T].
% \eea

We now plug (\ref{Fubini}) and (\ref{Duvsi}) into (\ref{pKQR}) to get:
\bea
\label{pKQR1}
&& \hE^0[\xi K_T]+\hE^0[\Th R_T] \nonumber\\
&=&\hE^0\Big\{\int_0^T\Big\{K_s\Big[-\a_s+L^u_s\hE^0\big[q_{s}\pa_y\si(s,X^u_{\cd\wedge s}, y,u_s)\big]\big|_{y=U^{u}_s}
+M_{s}L^u_s\pa_xh(s, X^u_s)\Big]\nonumber\\
&& +R_s\Big[-\b_s+\hE^0[q_{s}\si(s, X_{\cd\wedge s}, y, u_s)]\big|_{y=U^{u}_s}+\wt M_sh(s, X^u_s)\Big] +q_s C^{u,v}_s(v_s-u_s)\nonumber\\
&& +R_s\hE^0\big[q_s\pa_y\si(s,X^u_{\cd\wedge s}, y,u_s)\big]\big|_{y=U^{u}_s}(X^{u}_s-U^u_s)\Big\}ds+\int_0^TK_s d[^p\neg A^\si_s]\Big\}\\
&&+\hE^0\Big\{\int_0^T[K_sd\G_s+R_sd\Si_s]\Big\},\nonumber\\
&=& \hE^0\Big\{\int_0^T[-K_s\hat \a_s-R_s\hat\b_s+q_s C^{u,v}_\si(s)(v_s-u_s)]ds+K_s d[^p\neg A^\si_s]+[K_sd\G_s+R_sd\Si_s]\Big\},\nonumber
\eea
where
\bea
\label{alpha}
\left\{\ba{llllll}
\hat\a_t \dfnn \a_t-L^u_t\hE^0\big[q_t\pa_y\si(t,X^u_{\cd\wedge t}, y,u_t)\big]\big|_{y=U^{u}_t}
-\wt M_tL^u_t\pa_xh(t, X^u_t); \ms\\
\hat\b_t\dfnn \b_t-\hE^0[q_{t}\si(t, X_{\cd\wedge t}, y, u_t)]\big|_{y=U^{u}_t}-\wt M_th(t, X^u_t)\\
\qq-\hE^0\big[q_t\pa_y\si(t,X^u_{\cd\wedge t}, y,u_t)\big|_{y=U^{u}_t}(X^{u}_t-U^u_t).
%\ms \\
%\dis \hD_\si[q, X^u, u](t, dt)=\hE^0\Big[\int_t^Tq_r\m^0_\si(r, X^u_{\cd\wedge r}, u_r)(dt)dr\Big|\cF_t\Big].
\ea\right.
\eea
Combining (\ref{terminal}) and (\ref{pKQR1}) and using the processes $dA^\si$, $dA^f$ and their dual predicable projections,
we have
\bea
\label{terminal1}
0&\leq&
\hE^0\Big\{\int_0^T[-K_s\hat \a_s-R_s\hat\b_s+q_s C^{u,v}_\si(s)(v_s-u_s)]ds+\int_0^TK_sd[^p\neg\neg A^\si_s]\Big\}\\
&&+\hE^0\Big\{\int_0^T\Big[R_s\big[\hE^0[f(s, X_{\cd\wedge s}, y, u_s)]\big|_{y=U^{u}_s}+\hE^0\big[\pa_y f(s,X^u_{\cd\wedge s}, y,u_s)\big]\big|_{y=U^{u}_s}(X^{u}_s-U^u_s)\big]\nonumber\\
&&+L^u_sK_s\hE^0\big[\pa_y f(s,X^u_{\cd\wedge s}, y,u_s)\big]\big|_{y=U^{u}_s}+C^{u,v}_f(s)
(v_s-u_s)\Big]ds+\int_0^TK_s d[^p\neg\neg A^f_s]\Big\}\nonumber\\
&&+\hE^0\Big\{\int_0^T[K_sd\G_s+R_sd\Si_s]\Big\},\nonumber
\eea
where
%\bea
%\label{Cfuv}
$C^{u,v}_f(s)\dfnn \hE^0[L^u_s \pa_z f(s,\f_{\cd\wedge s}, U^u_s,z)]|_{\f=X^{u}, z=u_s}$.
%\ms\\
%\dis \hD_f[ X^u, u](t, dt)\dfnn \hE^0\Big[\int_t^T\m^0_f(t, X^u_{\cd\wedge r}, u_r)(dt)dr\Big|\cF_t\Big].
%\ea\right.\eea
%Here $\m^0_f(t,\f_{\cd\wedge t},z)(\cd)\in \sM[0,T]$ is defined in the same way as (\ref{mu0}).
Now, if we set $\Si_t=0$, and
\bea
\label{abcd}
\hat\a_t&=& L^u_t\hE^0\big[\pa_y f(t,X^u_{\cd\wedge t}, y,u_t)\big]\big|_{y=U^{u}_t}\nonumber\\
\hat\b_t&=& \hE^0[f(t, X_{\cd\wedge t}, y, u_t)]\big|_{y=U^{u}_t}+\hE^0\big[\pa_y f(t,X^u_{\cd\wedge t}, y,u_t)\big]\big|_{y=U^{u}_t}(X^{u}_t-U^u_t)\\
d\G_t&=&-d[^p\neg\neg A^\si_t] -d[^p\neg\neg A^f_t], \nonumber
\eea
then (\ref{terminal1}) becomes
\bea
\label{SMP0}
0&\le& \hE^0\Big\{\int_0^T [q_tC^{u,v}_\si(s)+C^{u,v}_f(s)](v_s-u_s)ds\Big\}, \q v\in\sU_{ad}.
\eea
From this we should be able to derive the maximum principle, provided that the adjoint equation (\ref{BSDE}) with coefficients
$\a$, $\b$, and
$\G$ determined by (\ref{alpha}) and (\ref{abcd}) is well-defined.
\begin{rem}
\label{hDsi}
{\rm 1) We remark that the process $\G$ in (\ref{abcd}) should be considered as a mapping from the space $L^2_\hF([0,T]\times\O)\times L^2_\hF(\O;\hC_T)\times
L^2_\hF([0,T]\times\O; U)$ to $\sM_\hF([0,T])$, the space of all the random measures on $[0,T]$, such that

(i) $(t,\o)\mapsto \m(t, \o, A)$ is $\hF$-progressively measurable, for all $A\in\sB([0,T])$;

(ii) $\m(t,\o, \cd)\in\sM([0,T])$ is a finite Borel measure on $[0, T]$.

\ms
2) Assumption \ref{Assum1}-(iii) implies that the random measure $\hD_\si[q,X^u, u](t, dt)$ satisfies the following estimate:
for any $q\in L^2_\hF([0,T]\times \O)$ and $u\in\sU_{ad}$,
\bea
\label{hDsiest}
\hE^0\Big[\int_0^T|d\,^p\neg\neg A^\si_t|\Big]&=&
%\Big|\hE^0\Big\{\int_0^T\hD_\si[q, X^u, u](t, dt)\Big\}\Big|=
\hE^0\Big\{\int_0^T\Big|\,^p\neg\Big(\int_t^Tq_s\m^0_\si(s, X^u_{\cd\wedge s}, u_s)(dt)ds\Big) \Big|\Big\}\nonumber\\
&\le& \hE^0\Big\{\int_0^T\int_0^s|q_s||\m^0_\si(s, X^u_{\cd\wedge s }, u_s)(dt)|ds\Big\}\\
&\le&\hE^0\Big\{\int_0^T|q_s|\int_0^s\ell(s,dt)ds\Big\}\le C\hE^0\Big\{\int_0^T|q_s|ds\Big\}\le C\|q\|_{2,2,\hQ^0}.\nonumber
\eea
The same estimate holds for $\hD_f[X^u, u](t, dt)$ as well.

\ms
3) Clearly, the processes $A^\si$ and $A^f$ are originated from the Fr\'echet derivatives of $\si$ and $f$, respectively, with respect to
the path $\f_{\cd\wedge t}$. If $\si$ and $f$ are  of Markovian type, then they will be absolutely
continuous with respect to the Lebesgue measure.
\qed}
\end{rem}

We shall now validate all the arguments presented above. To begin with, we note that the choice of $\a$, $\b$, and
$\G$ via by (\ref{alpha}) and (\ref{abcd}), together with the terminal condition $(\xi, \Th)$ by (\ref{pTQT}),  amounts to saying
that the processes $(p, q, \tilde q)$ and $(Q, M, \tilde M)$ solve the BSDE:
\bea
\label{BSDE1}
\left\{\ba{lllllll}
dp_t=-L^u_t\Big\{\hE^0\big[\pa_y f(t,X^u_{\cd\wedge t}, y,u_t)\big]\big|_{y=U^{u}_t}+\hE^0\big[q_t\pa_y\si(t,X^u_{\cd\wedge t}, y,u_t)\big]\big|_{y=U^{u}_t}\ms\\
\qq\q+\wt M_t\pa_xh(t, X^u_t)\Big\}dt-d\,^p\neg\neg A^\si_t - d\,^p\neg\neg A^f_t+q_tdB^1_t+\wt q_tdY_t\ms\\
dQ_t=-\Big\{\hE^0[q_{t}\si(t, X^u_{\cd\wedge t}, y, u_t)]\big|_{y=U^{u}_t}-\wt M_th(t, X^u_t)\ms\\
\qq\q+\hE^0\big[q_t\pa_y\si(t,X^u_{\cd\wedge t}, y,u_t)]\big|_{y=U^{u}_t}(X^{u}_t-U^u_t)\ms\\
\qq\q+\hE^0[f(t, X_{\cd\wedge t}, y, u_t)]\big|_{y=U^{u}_t}+\hE^0\big[\pa_y f(t,X^u_{\cd\wedge t}, y,u_t)\big]\big|_{y=U^{u}_t}(X^{u}_t-U^u_t)\Big\}dt\ms\\
\qq\q+M_tdB^1_t+\wt M_tdY_t,\ms\\
p_T=\xi, \q Q_T=\Th.
\ea\right.
\eea
Now if we denote $\eta=(p, Q)^T$, $W=(B^1, Y)^T$, $\Xi=\Big[\ba{ll}q&\tilde q\\ M&\tilde M\ea\Big]$, then we can rewrite
(\ref{BSDE1}) in a more abstract (vector) form:
\bea
\label{BSDE2}
\left\{\ba{ll}
d\eta_t=-\{A_t+\hE^0[G_t\Xi_tg(t,y)]\big|_{y=U^u_t}+H_t\Xi_t h_t\}dt-\G(\Xi)(t,dt)-\G_0(t,dt)+\Xi_tdW_t, \ms\\
\eta_T=\Upsilon,
\ea\right.
\eea
where $\Upsilon\in L^2_{\hF^W_T}(\O;\hQ^0)$;
$A, G, H$ and $h$ are bounded, vector or matrix-valued $\hF^W$-adapted processes with appropriate dimensions,
$g$ is an $\hR^2$-valued progressively measurable random field, and $U$ is an $\hF^Y$-adapted process. Moreover,
the $\hR^2$-valued finite variation processes $\G(\Xi)(t,dt)$ and $\G_0(t,dt)$ take the form:
\bea
\label{Gamma}
\G(\Xi)(t,dt)=~^p\neg\Big(\int_t^T\Xi_r\m^1_r(dt)dr\Big), \q \G_0(t, dt)=~^p\neg\Big(\int_t^T\m^2_r(dt)dr\Big),
\eea
where $r\mapsto \m^i_r(\cd)$, $i=1,2$, are $\sM[0,T]$-valued measurable random processes satisfying, as
measures with respect to the total variation norm,
\bea
\label{muest}
|\m^1_r(dt)|+|\m^2_r(dt)|\le \ell(r,dt), \q r\in[0,T], ~\hQ^0\mbox{a.s.}
\eea
We note that $\G(\Xi)(dt)$ and $\G_0(dt)$ are representing $d[^p\neg A^\si_t]$ and $[^p\neg A^f_t]$
%\hD_\si[q, X^u,u](t,dt)$ and $\hD_f[X^u,u](t,dt)$
in (\ref{BSDE1}), respectively, and can be substantiated by (\ref{A}) and (\ref{dpA}). Furthermore, by Assumption
\ref{Assum1}, they both satisfy (\ref{muest}). To the best of our knowledge,
BSDE (\ref{BSDE2}) is beyond all the existing frameworks of BSDEs, and we shall give a brief proof for its well-posedness.
\begin{thm}
\label{wellposed}
Assume that the Assumption \ref{Assum1} is in force. Then, the BSDE (\ref{BSDE2}) has a unique solution $(\eta, \Xi)$.
%p, q, \tilde q)$, $(Q, M,\tilde M)$.
\end{thm}

{\it Proof.} The proof is more or less standard, we shall only point out a key estimate. For any given $\widetilde{\Xi}^i\in  L^2_{\hF^W}([0, T]\times\O;\hR^4)$, obviously we have a unique solution $(\eta^i, \Xi^i)$ of (\ref{BSDE2}), $i=1, 2$, respectively, i.e.,
$$
\label{BSDE3}
\left\{\ba{ll}
d\eta^i_t=-\{A_t+\hE^0[G_t\widetilde{\Xi}^i_tg(t,y)]\big|_{y=U^u_t}+H_t\widetilde{\Xi}^i_t h_t\}dt-\G(\widetilde{\Xi}^i)(t,dt)-\G_0(t,dt)+\Xi^i_tdW_t, \ms\\
\eta_T^i=\Upsilon.
\ea\right.
$$
We define $\widehat{\xi}=\xi^1-\xi^2$, $\xi^i=\eta^i, \Xi^i$, $i=1, 2$, respectively. $\widehat{\widetilde{\Xi}}=\widetilde{\Xi}^1-\widetilde{\Xi}^2$. Noting the linearity of BSDE (\ref{BSDE2}) we see that $\widehat{\eta}$ satisfies:
\bea
\label{BSDE3}
\widehat{\eta}_t=\int_t^T\Big\{\hE^0[G_s\widehat{\widetilde{\Xi}}_s g(s,y)]\big|_{y=U^u_s}+H_s\widehat{\widetilde{\Xi}}_s h_s\Big\}ds+\int_t^T\G(\widehat{\widetilde{\Xi}})(s,ds)-M^T_t,
\eea
where $M^T_t\dfnn \int_t^T\widehat{\Xi}_sdW_s$. Therefore,
\beaa
|\widehat{\eta}_t+M_t^T|^2\le 2\Big\{\Big|\int_t^T\Big\{\hE^0[G_s\widehat{\widetilde{\Xi}}_s g(s,y)]\big|_{y=U^u_s}+H_s\widehat{\widetilde{\Xi}}_s h_s\Big\}ds\Big|^2+\Big|\int_t^T\G(\widehat{\widetilde{\Xi}})(s,ds)\Big|^2\Big\}.
\eeaa
Taking expectation on both sides above and noting that
$ \hE^0[\widehat{\eta}_t M^T_t]=0$ and
$$\hE^0\Big\{\Big|\int_t^T\Big\{\hE^0[G_s\widehat{\widetilde{\Xi}}_sg(s,y)]\big|_{y=U^u_s}+H_s\widehat{\widetilde{\Xi}}_s h_s\Big\}ds\Big|^2\Big\}\le C(T-t)\hE^0\Big[\int_t^T|\widehat{\widetilde{\Xi}}_s|^2ds\Big],$$
we have
\bea
\label{BSDEest}
 \hE^0[|\widehat{\eta}_t|^2] +\hE^0\Big[\int_t^T|\widehat{\Xi}_s|^2ds\Big]\le C(T-t)\hE^0\Big[\int_t^T|\widehat{\widetilde{\Xi}}_s|^2ds\Big] +\hE^0\Big\{\Big|\int_t^T\G(\widehat{\widetilde{\Xi}})(s,ds)\Big|^2\Big\}.
\eea
To estimate the term involving $\G(\widehat{\widetilde{\Xi}})$ we note that (recall (\ref{Gamma})) if a square-integrable process $V$ is increasing and continuous,
then so is its dual predictable projection $^pV$. Thus, by the definition of $^pV$ we have
\beaa
&&\hE^0\Big[\Big|\int_t^T d[^pV_s]\Big|^2\Big]=2 \hE^0\Big[\int_t^T(^pV_s-{}^pV_t)d[^pV_s]\Big]=
2 \hE^0\Big[\int_t^T(^pV_s-{}^pV_t)dV_s\Big]\\
&\le& 2\hE^0[(^pV_T-{}^pV_t)(V_T-V_t)]\le 2\Big(\hE^0\Big[\Big|\int_t^T d[^pV_s]\Big|^2\Big]\Big)^{1/2}
\Big(\hE^0\Big[\Big|\int_t^T dV_s\Big|^2\Big]\Big)^{1/2}.
\eeaa
That is,
\bea
\label{dpAest}
\hE^0\Big[\Big|\int_t^T d[^pV_s]\Big|^2\Big]\le 4\hE^0\Big[\Big|\int_t^T dV_s\Big|^2\Big].
\eea
Applying this to $V_t\dfnn \int_0^T\int_0^{t\wedge r}|\widehat{\widetilde{\Xi}}_r||\m^1_r(ds)|dr$, $t\in[0,T]$, we have
\beaa
\hE^0\Big[\Big|\int_t^T \G(\widehat{\widetilde{\Xi}})(s,ds)\Big|^2\Big]
&\le&\hE^0\Big[\Big|\int_t^T{}^p\Big(\int_s^T |\widehat{\widetilde{\Xi}}_r||\m^1_r(ds)|dr\Big)\Big|^2\Big]\le 4\hE^0\Big[\Big|\int_t^T\int_s^T |\widehat{\widetilde{\Xi}}_r| |\m^1_r(ds)|dr\Big|^2\Big]\\
&\le&4\hE^0\Big[\Big|\int_t^T \int_s^T |\widehat{\widetilde{\Xi}}_r|\ell(r,ds)dr\Big|^2\Big\}\\
&\le& C \hE^0\Big[\Big|\int_t^T |\widehat{\widetilde{\Xi}}_r|dr\Big|^2\Big\}\le C(T-t)\hE^0\Big[\int_0^T|\widehat{\widetilde{\Xi}}_s|^2ds\Big],
\eeaa
and therefore (\ref{BSDEest}) becomes
\bea
\label{BSDEest1}
\hE^0[|\widehat{\eta}_t|^2] +\hE^0\Big[\int_t^T|\widehat{\Xi}_s|^2ds\Big]\le C(T-t)\hE^0\Big[\int_t^T|\widehat{\widetilde{\Xi}}_s|^2ds\Big] .
\eea
With this estimate, and following the standard argument one shows that BSDE (\ref{BSDE1}) is well-posed on $[T-\d, T]$ for some (uniform)
$\d>0$. Iterating the argument one can then obtain the well-posedness on $[0,T]$. We leave the details to the interested reader.
\qed

We are now ready to prove the main result of this paper. Let us define the {\it Hamiltonian}:
for $(\f, \m)\in \hC_T\times \sP(\hC_T)$, and $k: [0, T]\times \Omega\rightarrow \hR$ adapted process, $(t, \o, z)\in[0,T]\times \O\times \hR$,
\bea
\label{Hamilton}
\sH(t, \o, \f_{\cd\wedge t}, \m, z; k)\dfnn k_t(\o)\cdot \sigma(t, \f_{\cd\wedge t}, \m,z)+f(t, \f_{\cd\wedge t}, \m, z).
\eea
We have the following theorem.
\begin{thm}[Stochastic Maximum Principle]
\label{theorem7.5}
%(SMP in integral form) For any $v\in\mathcal{U}$,
Assume that the Assumptions \ref{Assum1} and \ref{Assum2} hold. Assume further that the mapping $z\mapsto \sH(t, \f_{\cd\wedge t}, \m, z)$ is
convex. Let $u=u^*\in\sU_{ad}$ be an optimal control and $X^u$ the corresponding trajectory. Then,
for $dt\times d\hQ^0$-a.e. $(t,\o)\in[0,T]\times\O$ it holds that
\bea
\label{7.57}
\sH(t, \o, X_{\cd\wedge t}^u, \m^u_t, u_t; q_t)=\inf_{v\in {U}}\sH(t, \o, X_{\cd\wedge t}^u, \m^u_t, v; q_t),
\eea
 where $(p,q, \tilde q)$ and $(Q, M, \tilde M)$ constitute the unique solution of the BSDE (\ref{BSDE1}).
\end{thm}

{\it Proof.} We first recall from (\ref{DBC})  that
\beaa
C^{u,v}_f(t)&=& \hE^0[L^u_t \pa_z f(t,\f_{\cd\wedge t}, U^u_t,z)]|_{\f=X^{u}, z=u_t}=\pa_zf(t,X^u_{\cd\wedge t}, \m^u_t,u_t);\\
C^{u,v}_\si(t) &=&\hE^0\Big\{L^u_t \pa_z\si(t, \f_{\cd\wedge t}, U^u_t, z) \Big]\Big\}\Big|_{\f=X^u; z=u_t}=\pa_z\si(t,X^u_{\cd\wedge t}, \m^u_t,u_t).
\eeaa
Then (\ref{SMP0}) implies that
\bea
\label{7.55}
%\sH(t, E^v[X_t^v|{\mathcal{F}}_t^Y], \varphi, k) &\dfnn& q_t\cdot E^v[\sigma(E^v[X_t^v|\mathcal{F}_t^Y],\varphi,k)]+E^v[f(E^v[X_t^v|\mathcal{F}_t^Y],\varphi,k)].\\
0&\le& \hE^0\Big[\int_0^T [q_tC^{u,v}_\si(t)+C^{u,v}_f(t)](v_t-u_t) dt\Big]\\
&=&\hE^0\Big[\int_0^T\pa_z\sH(t, \o, X_{\cd\wedge t}^u, \m^u_t, u_t; q_t)(v_t-u_t) dt\Big]. \nonumber
\eea
Therefore for $dt\times d\hQ^0$-a.e. $(t,\o)\in[0,T]\times\O$, and any $v\in U$,  it holds that
\bea
\label{7.56}
\pa_z\sH(t, \o, X_{\cd\wedge t}^u, \m^u_t, u_t; q_t)(v-u_t)\ge 0.
\eea
Now, for any $v\in {U}$, one has, $dt\times d\hQ^0$-a.e.  on $ [0,T]\times\Omega$,
\beaa
&&\sH(t, \o, X_{\cd\wedge t}^u, \m^u_t, v; q_t)-\sH(t, \o, X_{\cd\wedge t}^u, \m^u_t, u_t; q_t)\\
&=&\int_0^1 \pa_z\sH(t, \o, X_{\cd\wedge t}^u, \m^u_t, u_t+\l(v-u_t); q_t)(v-u_t)d\l\\
&=&\int_0^1\Big[\pa_z\sH(t, \o, X_{\cd\wedge t}^u, \m^u_t, u_t+\l(v-u_t); q_t)-\pa_z\sH(t, \o, X_{\cd\wedge t}^u, \m^u_t, u_t; q_t)\Big](v-u_t)
d\l\\
&& +\pa_z\sH(t, \o, X_{\cd\wedge t}^u, \m^u_t, u_t; q_t)(v-u_t)\geq 0,.
\eeaa
Here the first integral on the right hand side above is nonnegative due to the convexity of $\sH$ in variable $z$, and the last term is non-negative because of (\ref{7.56}). The identity (\ref{7.57}) now follows immediately.
\qed

\begin{rem}
{\rm In stochastic control literature the inequality (\ref{7.55}) is sometimes referred to as {\it Stochastic Maximum Principle in integral form}, which
in many applications is useful, as it does not require the  convexity assumption on the Hamiltonian $\sH$.
\qed}
\end{rem}

\bs
\noindent{\bf Acknowledgment.}  We would like to thank the anonymous referee for his/her very careful reading of the manuscript
and many incisive and constructive questions and suggestions, which helped us to make the paper a much better product.


\begin{thebibliography}{1}
\bibitem{Bens}
Bensoussan, A.  (1992), {\sl Stochastic Control of Partially Observable Systems}, Cambridge University Press.
\bibitem{BDL} Buckdahn, R., Djehiche, B., and Li, J.,  (2011), {\it  A general Stochastic Maximum Principle for SDEs of Mean-field Type}. Appl. Math. Optim. {\bf 64}, no. 2, 197--216.

\bibitem{BDLP} Buckdahn, R., Djehiche, B., Li, J., and Peng, S. (2009), {\it Mean-field Backward Stochastic Differential Equations: A Limit Approach.} Ann. Probab. {\bf 37},   no. 4, 1524-1565.


\bibitem{BLP} Buckdahn, R., Li, J., and Peng, S. (2009), {\it Mean-field Backward Stochastic Differential Equations and Related Partial Differential Equations. Stochastic Process}. Appl. {\bf 119}, no. 10, 3133--3154.

\bibitem{BMZ1}
Buckdahn, R., Ma, J., and Zhang, J. (2015), {\it Pathwise Taylor Expansions for Random Fields on
Multiple Dimensional Paths},  Stochastic Process. Appl. {\bf 125},  no. 7, 2820--2855.
%preprint, arXiv:1310.0517v1 [math.PR] 1 Oct 2013.

\bibitem{BMZ2}
Buckdahn, R., Keller, C., Ma, J., and Zhang, J. {\it Pathwise Viscosity Solutions of SPDEs and Forward PPDEs}, preprint.

\bibitem{CD1} Carnoma R. and Delarue, F.  (2012), {\it Optimal control of McKean-Vlasov stochastic dynamics}, Technical Report.

\bibitem{CD2} Carnoma R., Delarue, F., and  Lachapelle, A.  (2013), {\it Control of MaKean-Vlasov versus Mean Field Games}, Math. Financ. Econ. {\bf 7}, no. 2, 131-166.

\bibitem{CD3} Carnoma R. and Delarue, F.  (2013), {\it Probabilistic Analysis of Mean-Field Games},  SIAM J. Control Optim. {\bf 51}, no. 4, 2705-2734.

\bibitem{CD4} Carmona, R. and Delarue, F. (2013), {\it Mean Field Forward-Backward Stochastic Differential Equations},  Electron. Commun. Probab. {\bf 18},  no. 68, 15 pp.

\bibitem{CD5} Carmona, R. and  Delarue, F., (2015), {\it Forward-backward stochastic differential equations and controlled
McKean-Vlasov dynamics}. Ann. Probab. {\bf 43}, no. 5, 2647-2700.

\bibitem{CZ}
Carmona, R. and Zhu, X. (2014), {\it A probabilistic approach to mean field games with major and minor players}, arXiv: 1409.7141v [math.PR].

\bibitem{Cauty}
Cauty, R. (2001), {\it Solution du probl\`{e}me de point fixe de Schauder}, Fundamenta Matematicae, {\bf 170}, 231-246.

\bibitem{Cauty2}
Cauty, R.  (2012), {\it Une g\'en\'eralisation de la conjecture de point fixe de Schauder}, arXiv: 1201.2586 [math.AT].

\bibitem{Dupire}
Dupire, B. (2009), {\it Functional It\^o calculus},  papers.ssrn.com.

\bibitem{EK}
Ethier, S., and Kurtz, T. (1986), {\sl Markov Processes, Characterization and Convergence}, John Williams \& Sons Inc.


%\bibitem{GT}
%Gilbarg, D. and Trudinger, N., {\sl Elliptic Partial Differential Equations of Second Order}, %Springer-Verlag, Berlin-Heidelberg-New York, 3rd printing 1998.

\bibitem{HMC}
Huang, M., Malham\'{e}, R., and Caines, P. (2006), {\it Large population stochastic dynamic games: closed-loop
McKean-Vlasov systems and the Nash certainty equivalence principle}, Commun. Inf. Syst. {\bf 6} (3),
221-252.

\bibitem{IW}
Ikeda, N. and Watanabe, S. (1981), {\sl Stochastic differential equations and diffusion processes}, North-Holland, Amsterdam.

\bibitem{Krylov}
Krylov, N. V. (1999), {\it An analytic approach to SPDEs}. {\sl Stochastic partial differential equations: six perspectives}, 185-242, Math. Surveys Monogr., 64, Amer. Math. Soc., Providence, RI.

\bibitem{LasryLions}
Lasry, J. M. and Lions, P.L. (2007), {\it Mean Field Games}, Japanese Journal of Mathematics, {\bf 2} (1), Mar.

\bibitem{Li}
Li, J.  (2012), {\it Stochastic maximum principle in the mean-field controls}, Automatica J. IFAC. {\bf 48}, no. 2, 366-373.

\bibitem{LiMin}
Li, J. and Min, H. (2016), {\it Weak solutions of mean-field stochastic differential equations and application to zero-sum stochastic differential games}, to appear in SIAM Control and Optimization.

\bibitem{SV}
	Strook, D., and Varadhan, S.V.S. (1979), {\sl Multidimensional Diffusion Processes}, Springer-Verlag, New York.
	
		
\bibitem{Villani}
Villani, C. (2003), {\sl Topics in optimal transportations}. Graduate Studies in Mathematics, {\bf 58}, AMS, Providence, RI.
	
\bibitem{yong-zhou}
	Yong, J. and Zhou, X. (1999), {\sl Stochastic Controls: Hamiltonian Systems and HJB Equations},
	New York: Springer-Verlag.
	
\bibitem{Zeit}
Zeitouni, O. (1986), {\it On the reference probability approach to the equations of nonlinear filtering}, Stochastics,
vol. 19,  133-149.
\end{thebibliography}
  \end{document}